\documentclass[12pt]{article}
\pdfoutput=1    
\usepackage{amscd,amssymb,stmaryrd}
\usepackage{amsmath}
\usepackage{amscd,amssymb,stmaryrd}
\usepackage{graphicx}
\usepackage{subfigure,fancybox}
\usepackage{amscd,amstext}
\usepackage{amsmath}
\usepackage{mathtools}
\usepackage{setspace}
\usepackage{amssymb}
\usepackage{slashed}
\usepackage{pstricks,pst-node,pst-plot,pst-coil}
\usepackage{graphicx}
\usepackage{ntheorem}
\usepackage{upgreek}
\usepackage{enumerate}
\usepackage{graphicx}
\usepackage{cite}
\usepackage{tikz}
\usepackage{authblk}
\usepackage{txfonts}
\usepackage[colorlinks,linkcolor=black,anchorcolor=black,citecolor=black,CJKbookmarks=True]{hyperref}
\theorembodyfont{\normalfont}
\newtheorem{thm}{Theorem}[section]
\newtheorem{lem}{Lemma}[section]
\newtheorem*{pf}{PROOF:}
\newtheorem{prop}{Proposition}[section]
\newtheorem{defi}{Definition}[section]

\newtheorem{remark}{Remark}[section]
\newtheorem{cor}{Corollary}[section]

\numberwithin{equation}{section}

\title{\large Formation of shifted shock for the 3D compressible Euler equations with time-dependent damping}
\author[1]{\normalsize Chen Zhendong\footnote{Acknowledgement:
This is part of the Ph. D thesis of the author written under the supervision of Professor Zhouping Xin at the Institute of Mathematical Sciences of The Chinese University of Hong Kong. }}

\affil[1]{Institute of Mathematical Science, The Chinese University of Hong Kong, Shatin, NT, Hong Kong.}
\date{}
\begin{document}           
\def\ra{\rightarrow}\def\Ra{\Rightarrow}\def\lra{\Longleftrightarrow}\def\fai{\varphi}
\def\s{\slashed}\def\Ga{\Gamma}\def\til{\tilde}
\def\de{\Delta}\def\fe{\phi}\def\sg{\slashed{g}}\def\la{\lambda}\def\R{\mathbb{R}}\def\m{\mathbb}
\def\a{\alpha}\def\p{\rho}\def\ga{\gamma}\def\lam{\lambda}\def\ta{\theta}\def\sna{\slashed{\nabla}}
\def\pa{\partial}\def\be{\beta}\def\da{\delta}\def\ep{\epsilon}\def\dc{\underline{\chi}}
\def\si{\sigma}\def\Si{\Sigma}\def\wi{\widetilde}\def\wih{\widehat}\def\beeq{\begin{eqnarray*}}
\def\eeq{\end{eqnarray*}}\def\na{\nabla}\def\lie{{\mathcal{L}\mkern-9mu/}}\def\rie{\mathcal{R}}
\def\ud{\underline}\def\les{\lesssim}\def\ka{\kappa}
\def\dl{\underline{L}}\def\du{\underline{u}}
\def\hs{\hspace*{0.5cm}}\def\bee{\begin{equation}}\def\ee{\end{equation}}\def\been{\begin{enumerate}}
\def\een{\end{enumerate}}\def\bes{\begin{split}}\def\zpi{\prescript{(Z)}{}{\pi}}\def\szpi{\prescript{(Z)}{}{\slashed{\pi}}}
\def\ees{\end{split}}\def\hra{\hookrightarrow}\def\udx{\ud{\xi}_{A}}\def\ude{\ud{\eta}_{A}}
\def\tpi{\prescript{(T)}{}{\pi}}\def\lpi{\prescript{(L)}{}{\pi}}\def\qpi{\prescript{(Q)}{}{\pi}}
\def\stpi{\prescript{(T)}{}{\slashed{\pi}}}\def\sqpi{\prescript{(Q)}{}{\slashed{\pi}}}
\def\zgpi{\prescript{(Z)}{}{\wi{\pi}}}\def\szgpi{\prescript{(Z)}{}{\wi{\slashed{\pi}}}}
\def\tgpi{\prescript{(T)}{}{\wi{\pi}}}\def\lgpi{\prescript{(L)}{}{\wi{\pi}}}\def\qgpi{\prescript{(Q)}{}{\wi{\pi}}}
\def\stgpi{\prescript{(T)}{}{\wi{\slashed{\pi}}}}\def\sqgpi{\prescript{(Q)}{}{\wi{\slashed{\pi}}}}
\def\pre #1 #2 #3{\prescript{#1}{#2}{#3}}\def\spi{\slashed{\pi}}\def\rgpi{\prescript{(R_{i})}{}{\wi{\pi}}}
\def\rpi{\prescript{(R_{i})}{}{\pi}}\def\srpi{\prescript{(R_{i})}{}{\s{\pi}}}
\def\pp{\uprho}\def\vw{\upvarpi}\def\srgpi{\pre {(R_{i})} {} {\wi{\s{\pi}}}}
\def\supnormda{L^{\infty}(\Si_{t}^{\tilde{\da}})}\def\supnormu{L^{\infty}(\Si_{t}^{u})}
\def\normda2{L^{2}(\Si_{t}^{\tilde{\da}})}\def\normu2{L^{2}(\Si_{t}^{u})}
\def\les{\lesssim}\def\sdfai{|\s{d}\fai|}\def\sd{\s{d}}\def\mari{\mathring}\def\wv{\dfrac{2(1+t')}{\frac{1}{2}\til{tr}\chi}}
\def\oli{\overline}\def\vta{\vartheta}\def\ol{\overline}\def\Llra{\Longleftrightarrow}
\def\ea+2{E_{\leq|\a|+2}}\def\dea+2{\ud{E}_{\leq|\a|+2}}\def\fa+2{F_{\leq|\a|+2}}\def\dfa+2{\ud{F}_{\leq|\a|+2}}
\def\tea+2{\wi{E}_{\leq|\a|+2}}\def\tdea+2{\ud{\wi{E}}_{\leq|\a|+2}}\def\tfa+2{\wi{F}_{\leq|\a|+2}}\def\tdfa+2{\ud{\wi{F}}_{\leq|\a|+2}}
\def\ntea{\wi{E}_{\leq|\a|+1}}\def\ntdea{\ud{\wi{E}}_{\leq|\a|+1}}\def\ntfa{\wi{F}_{\leq|\a|+1}}\def\ntdfa{\ud{\wi{F}}_{\leq|\a|+1}}
\def\ba2{b_{|\a|+2}}\def\Ka{K_{\leq|\a|+2}}\def\mumba{\mu_m^{-2\ba2}}\def\nba{b_{|\a|+1}}
\def\endpf{\hfill\raisebox{-0.3cm}{\rule{0.1mm}{3mm}\rule{0.5mm}{0.1mm}\rule[-0.6mm]{2.3mm}{0.7mm}\raisebox{3mm}[0pt][0pt]{\makebox[0pt][r]{\rule{2.9mm}{0.1mm}}}\rule[-0.6mm]{0.1mm}{3.5mm}\rule[-0.6mm]{0.5mm}{3mm}} \vspace{0.5cm}\\}

\pagestyle{myheadings} \thispagestyle{empty} \markright{}
\maketitle
\begin{abstract}
\hs In this paper, we show the shock formation to the compressible Euler equations with time-dependent damping $\frac{a\p u}{(1+t)^{\lam}}$ in three spatial dimensions without any symmetry conditions. It's well-known that for $\lam>1$, the damping is too weak to prevent the shock formation for suitably large data. However, the classical results only showed the finite existence of the solution. Follow the work by D.Christodoulou in\cite{christodoulou2007}, starting from the initial isentropic and irrotational short pulse data, we show the formation of shock is characterized by the collapse of the characteristic hypersurfaces and the vanishing of the inverse foliation density function $\mu$, at which the first derivatives of the velocity and the density blow up, and the lifespan $T_{\ast}(a,\lam)$ is exponentially large. Moreover, the damping effect will shift the time of shock formation $T_{\ast}$. The methods in the paper can also be extended to the Euler equations with general time-decay damping. 
\end{abstract}

\tableofcontents
\section{\textbf{Introduction}}\label{section1}
\hs In this paper, we will consider the following compressible Euler equations with time-dependent damping in $R^3$:
\bee\label{cee}
\left\{\bes
&\dfrac{\pa}{\pa t}\p+\nabla\cdot(\p v)=0,\\
&\dfrac{\pa}{\pa t}(\p v)+\nabla\left(\p v\otimes v+p I_3\right)+
\dfrac{a}{(1+t)^{\lam}}\p v=0,\\
&\dfrac{\pa}{\pa t}S+v\cdot\nabla S=0,
\end{split}\right.
\ee
where $\p$, $v=(v^1,v^2,v^3)$, $p$, $S$ represent the density, the velocity, the pressure and the entropy of the flow, respectively. $I_3$ is a $3\times 3 $ identical matrix and $a,\lam$ are the damping constants to describe the decay-rate of the damping concerning time. These equations describe the motion of a perfect fluid which follow from the conservation of mass, momentum and energy respectively. The enthalpy $h$ of a thermodynamic system is defined as the sum of its internal energy and the product of its pressure and volume:
 \bee
 h=e+pV.
 \ee
 Here $e$ and $V=\frac{1}{\p}$ represent the internal energy and the volume of the system, respectively. It follows from the relation of thermodynamics that
 \bee\label{enthalpy}
 dh=TdS+Vdp,
 \ee
 where $T$ is the absolute temperature. The sound speed $\eta$ is defined as $\eta=\sqrt{\dfrac{\pa p}{\pa\p}}$. The equation of state is given by $p=p(\p,S)$ and we assume that:
   \bee\label{nonlinearcondition}
   p\ \text{is not linear in}\ \dfrac{1}{\p}.
   \ee
  This assumption will be clarified later. We additional require $\lam>\frac{3}{2}$, which will be explained in section5.1.
  \subsection{\textbf{Brief review of former results}}
\hs There many results established before and the time dependent damping makes the problem more complicated compared with the case $\lam=0$\footnote{For the former results of $\lam=0$ and $a=0$, see \cite{CZD1} section1.1.}. For $1D$ isentropic case, in\cite{PAN2016435} and \cite{PAN2016327}, Pan showed that for $0<\lam<1$, $a>0$ or $\lam=1$, $a>2$, the smooth solution exists globally for the small data(in $H^1$ sense) while for $\lam>1$, $a>0$ or $\lam=1$, $0\leq a\leq 2$, the solution will blow up in finite time. Later, such results were generalized by Chen-Li-Li-Mei-Zhang \cite{CHEN20205035} with general time-decay damping and general initial data. For the $3D$ isentropic case, Hou-Witt-Yin in \cite{hou2017global} and \cite{HOU2017949} proved that for $0<\lam\leq 1$, if the initial data is irrotational and sufficiently small, then the damping is strong enough to prevent the formation of singularity and the system admits a global solution. For $\lam>1$, the lifespan $T$ of the solution is finite for suitably large data. Later in\cite{doi:10.1080/00036811.2020.1722805}, Pan derived the similar results and showed the decay rate of the solution in the case of $0<\lam\leq 1$ and long time behavior and the optimal decay rates for the solution were studied by Ji and Mei \cite{meiming1} and \cite{meiming2}.\\
\hs However, our work which follows from Christodoulou\cite{christodoulou2014compressible} is different from the above works in the following sense:
\been[(1)]
\item In multi-dimensional case, all of them only showed the finite existence of the solution while we will give a precise description of the mechanism of the break down of the solution(shock formation).
\item To prove the finite existence of the solution, they imposed some certain conditions on the initial data while we will consider a class of short pulse data(small in low norm but large in high norm).
\een
\hs In multi-dimensional case, a major breakthrough in understanding the shock mechanism for hyperbolic systems was achieved by D.Christodoulou \cite{christodoulou2014compressible}. Starting with initial pulse data, Christodoulou showed the shock formation to compressible isentropic and irrotational Euler¡¯s equations in three space dimension. In particular, he obtained a complete geometric description of the shock formation and a detailed analysis of the behavior of the solution near shock.\\
\hs Later, Yu-Miao\cite{Ontheformationofshocks} applied Christodoulou's framework to the quasilinear wave equation in three space dimensions. By constructing a family of short pulse initial data, they showed how the solution breaks down near the first singularity(shock) and gave a sufficient condition leading to the shock formation.\\
\hs In 2018, J.Speck and J.Luk studied $2D$ Euler system without the assumption of irrotation by applying Christodoulou's framework. In\cite{luk2018shock}, they studied plane symmetric data with initial short pulse perturbation, and for such initial data, they showed the shock mechanism of the Euler system such that the first derivatives of the $v$ and $\p$ blow up, which is due to the collapse of the characteristic hypersurfaces. Then, they generalized this results into $3$D case in\cite{Specklukshockformation}.\\ 
\hs Recently, there have been many very interesting works on global well-posedness or finite blowup of smooth solutions to various nonlinear multi-dimensional wave equations by adapting the Christodoulou's framework, see\cite{yin1,yin2,yin3,Speck2022,Speck2020,Speck2017,Speck20172,Speck2016,Speck20162,Speck2014,Speck20142,yupin2019,Speck20173}.\\
\hs More recently, Buckmaster, Shkoller, Vicol obtained a series results for the shock formation to the Euler system in multi-dimensional case. In\cite{BSV2Disentropiceuler}, they considered the $2$D isentropic Euler equations under azimuthal symmetry with smooth initial datum of finite energy and nontrivial vorticity. By using the modulated self-similar variables, they obtained the point shock forms in finite time with explicitly computable blow-up time and location and the solution near shock is of cusp type. Later, in\cite{BSV3Disentropiceuler} and \cite{BSV3Dfulleuler}, they generalized the above results into $3$D isentropic and non-isentropic Euler system. \\
\hs Following the framework of Christodoulou\cite{christodoulou2014compressible} and the framework in\cite{CZD1}, we generalized the result in\cite{CZD1} into time damping case. Despite the difficulties in \cite{CZD1}, our work are different from the work \cite{CZD1} in the following sense:
\been[(1)]
\item In \cite{CZD1}, we restrict the time to be $[-2,\sigma)$ where $\sigma<0$ is a fixed small constant so that the time is finite, while in this paper, the growth of time is necessary to be considered so that the framework in\cite{CZD1} should be modified. In particular, one has to control the decay of the solution together with its derivatives(see \eqref{bs}) and the two multipliers are re-chosen which are adapted from\cite{christodoulou2014compressible}. This difficulty leads to the require of more delicate estimates in the energy estimates. Precisely, in \cite{CZD1}, to close the energy estimates, the coefficients of the energies $E_0,E_1$ only need to be small(compared with 1) while in this paper, the corresponding terms are required to be integrable in time(decay faster than $\frac{1}{(1+t)^{1^+}}$). Moreover, to obtain the high order energy estimates, the energies and fluxes need to be modified cautiously(see the definition\eqref{modifiedenergy}).
\item Consider the following covariant wave equation:
      \[
      \square_g\fai=f.
      \]
Due to the damping effect, the inhomogeneous terms $f$ become more complicated which contains the most difficult term $\frac{a}{(1+t)^{\lam}}\dl\fai$(In \cite{CZD1}, the corresponding term is $a\dl\fai$). Due to the bad decay of the derivative of $\fai$ along the incoming null-hypersurfaces(i.e. $\dl\fai$), one can not bound this term in $L^2$ norm directly which require additional estimates, while in \cite{CZD1}, one can bound this term easily. Moreover, this term will illustrate how the damping effect shift the shock formation.
\een
\subsection{\textbf{Notations}}
Through the whole paper, the following notations will be used unless stated otherwise:
\begin{itemize}
\item Latin indices $\{i,j,k,l,\cdots\}$ take the values $1, 2, 3, $ Greek indices $\{\a,\be,\ga,\cdots\}$
take the values $ 0, 1, 2, 3$ and capital letter $\{A,B,C,\cdots\}$ take values $1,2.$ Repeated indices are meant to be summed.
\item The convention $f\les h$ means that there exists a universal positive constant $C$ such that $f\leq Ch$.
\item The notations l.o.ts (lower order terms) mean the terms are of lower order. Here, the order means the number of total derivatives acting on $\fai$ and we set $\fai$ to be order $0$. For example, one can rewrite $\pa^2\fai+\pa\fai$ as $\pa^2\fai+$l.o.ts. Let $O_{\frac{b}{c},d}^{\leq a}$ be the terms of
order $\leq a$ with bound $\dfrac{(1+\ln(1+t))^{b}\da^{\frac{d}{2}}}{(1+t)^c}$. 
\item For the metric $g_{\a\be}$, $g^{\a\be}$ means its inverse such that $g_{\a\be}g^{\be\ga}=\da^{\ga}_{\be}$ with $\da^{\ga}_{\be}$ being the Kronecker symbol.
\item The box operator $\Box_g:=g^{\a\be}D^2_{\a\be}$ denotes the covariant wave operator corresponding to the spacetime metric g and $\s{\de}:=\sg^{AB}\s{D}^2_{AB}$ denotes the covariant Laplacian corresponding to $\sg$ on $S_{t,u}$, where $D$, $\s{D}$ are the Levi-Civita connections corresponds to $g$, $\sg$ respectively. We also denote $\til{D}$ to be the Levi-Civita connection corresponds to $\til{g}$.
\item For a object $q$, $\s{q}$ means its restriction(projection) on $S_{t,u}$. In particular, $\s{div}$ represents the divergence operator on $S_{t,u}$ such that $\s{div}Y:=\s{D}_AY^A$ for any $S_{t,u}$ vector field $Y$. $\s{d}$, $\lie$ represent the restriction of the standard differentiation $d$ and the lie derivative $\mathcal{L}$ on $S_{t,u}$ respectively.
\item For $q$ being a $(0,2)$ tensor and $Y,Z$ being the vector fields, set the contraction as
      \bee
      q_{YZ}:=q_{\a\be}Y^{\a}Z^{\be},
      \ee
      and similar for the other type tensors.
\item For a spacetime vector field $V$, denote the decomposition of $V$ relative to the frame $\{L,T,X\}$ as
      \bee
      V=V^LL+V^TT+V^AX_A,
      \ee
      with $V^L,V^T$ and $V^A$ being functions.
\item For a $(0,2)$ tensor $\ta$, the following decomposition holds:
      \bee
      \ta=\hat{\ta}+\frac{1}{2}tr\ta\cdot g,
      \ee
      where $\hat{\ta}$ is the trace free part of $\ta$ and $tr\ta=g^{\a\be}\ta_{\a\be}$ is the trace of $\ta$ with respect to $g$.
\end{itemize}
\subsection{\textbf{Geometric blow up for the Burgers equation with time-damping}}
Consider the following Cauchy problem for $1$D Burger's equation with time-damping
\bee\label{burgersdamp}
\left\{\bes
&\pa_t\fai+\fai\pa_x\fai=-\dfrac{a\fai}{(1+t)^{\lam}},\\
&\fai(x,t=0)=f(x).
\end{split}\right.
\ee
where $a$ and $\lam$ are the damping constants with $\lam>1$. For simplicity, we assume additionally that $f(0)=0$ and $\min \pa_xf(x)=\pa_xf(0)=-c$. Then for $a=0$, this problem returns to the standard Burger's equation and by standard characteristic method, the solution of \eqref{burgersdamp} will form a shock at time $T_{\ast}=\frac{1}{c}$ for $c>0$ and the location $x_{\ast}=0$ with $\pa_xu(0,t)\ra-\infty$ as $t\to T_{\ast}$. Define the Eikonal function as
 \bee
 \left\{\bes
 &Lu:=(\pa_t+\fai\pa_x)u=0,\\
 &u(x,0)=x.
 \end{split}\right.
 \ee
 This Eikonal function yields a new coordinates system $(t,u)$ in which $L=\dfrac{\pa}{\pa t}|_{(t,u)}$. Then, the original equation becomes
 \bee
 L\fai=-\dfrac{a\fai}{(1+t)^{\lam}}.
 \ee
 Hence, $\fai(t,u)=\fai_0(u)\exp\left\{\frac{a}{\lam-1}\left((1+t)^{1-\lam}-1\right)\right\}$ and if initial $\fai_0$ is smooth, $\fai$ together with it's derivatives will remain smooth for all time in $(t,u)$ coordinates. Define the inverse foliation density function as\footnote{Note that $\mu$ exactly characterizes the shock formation of \eqref{burgersdamp} and the Jacobian of coordinates transformation. See \cite{CZD1} section1.3.}
 \bee
 \mu=\dfrac{1}{\pa_xu}.
 \ee
 Then, $\mu$ satisfies the equation
\bee
\left\{\bes
&L\mu=\mu\pa_x\fai=\dfrac{\pa \fai}{\pa u},\\
&\mu(x,0)=1.
\end{split}\right.
\ee
It follows that
\begin{align}
\mu(t,x)&=\mu_0+\pa_xf(x)\int_{0}^t\exp\left\{\frac{a}{\lam-1}\left((1+s)^{1-\lam}-1\right)\right\} ds\\
&\geq 1-c\int_{0}^t\exp\left\{\frac{a}{\lam-1}\left((1+s)^{1-\lam}-1\right)\right\} ds.
\end{align}
\hs Therefore,
\been[(1)]
\item If $c>0$, then the damping effect is too weak to prevent the shock formation and $\mu\to 0$ as $t\to T_{\ast}$ with $T_{\ast}$ defined by
\bee
\int_0^{T_{\ast}}\exp\left\{\frac{a}{\lam-1}\left((1+s)^{1-\lam}-1\right)\right\} ds=\frac{1}{c}.
\ee Furthermore, as $t\to T_{\ast}$, $\pa_x\fai=\dfrac{1}{\mu}\dfrac{\pa \fai}{\pa u}$ blows up like $\dfrac{1}{\mu}$. Also, it can been seen that the damping constants $a$ and $\lam$ will shift $T_{\ast}$ in the following sense compared with the undamped case:
\begin{itemize}
\item for fixed $\lam>1$, if $a>0$, then the damping will delay the shock formation; while for $a<0$, the anti-damping will lead to the shock formation in advance;
\item for fixed $a$, then as $\lam$ increasing, the damping becomes weaker and the blowup time $T_{\ast}$ will close to $T_0=1$, which is the blowup time for the undamped case.
\end{itemize}
\item If $c\leq0$, then $\mu>0$ for all $t$, that is, the characteristics never intersect and we will obtain a global solution for \eqref{burgersdamp}.
    \een
\subsection{\textbf{The inhomogeneous nonlinear wave equation for the Euler system with damping}}
\hs For a irrotational flow, there exists a potential function $\fe$ such that $v=-\nabla\fe$, while for a smooth isentropic flow, the entropy will remain constant from the equation of energy. Then, it follows from the equation of momentum that up to a constant,
\bee
h=\dfrac{\pa}{\pa t}\fe-\dfrac{1}{2}|\nabla\fe|^2+\dfrac{a}{(1+t)^{\lam}}\fe.
\ee
 Then, the continuity equation can be rewritten as
\begin{align*}
&-\dfrac{\p}{\eta^2}\left[
\dfrac{\pa^2\fe}{\pa t^2}-\eta^2\de\fe-2\dfrac{\pa\fe}{\pa x^i}\dfrac{\pa ^2\fe}{\pa x^i\pa t}+\dfrac{\pa\fe}{\pa x^i}\dfrac{\pa \fe}{\pa x^j}\dfrac{
\pa^2\fe}{\pa x^i\pa x^j}+\dfrac{a}{(1+t)^{\lam}}\left(\dfrac{\pa\fe}{\pa t}-|\nabla\fe|^2-\dfrac{\lam}{(1+t)}\fe\right)\right]=0\\
\end{align*}
i.e.,
\bee\label{mass2}
g^{\a\be}\pa_{\a}\pa_{\be}\fe=\dfrac{1}{\eta^2}\dfrac{a}{(1+t)^{\lam}}\left(\dfrac{\pa\fe}{\pa t}-|\nabla\fe|^2-\dfrac{\lam}{(1+t)}\fe\right),
\ee
where $g$ is the acoustical metric given by
\bee
\left(\begin{array}{cccc}
-\eta^2+\sum_{i}(v^i)^2 & -v^1 & -v^2 & -v^3\\
-v^1 & 1 &0 &0\\
-v^2 & 0 & 1 &0\\
-v^3 & 0 & 0 &1
\end{array}\right)
\ee
with it's inverse $g^{-1}$
\bee
\dfrac{1}{\eta^2}\left(\begin{array}{cccc}
-1 & -v^1 & -v^2 & -v^3\\
-v^1 & \eta^2-(v^1)^2 &-v^1v^2 &-v^1v^3\\
-v^2 & -v^1v^2 & \eta^2-(v^2)^2 &-v^2v^3\\
-v^3 & -v^1v^3 & -v^2v^3 &\eta^2-(v^3)^2
\end{array}\right).
\ee
\hs Let $A$ be any one of $\{\pa_{i},\Omega_{ij}=x^i\pa_j-x^j\pa_i\}$. Then, acting $A$ on the equation \eqref{mass2} directly yields the following inhomogeneous covariant wave equations for $\fai$\footnote{\eqref{nonlinearwave'} is called the linearized equation of\eqref{mass2} with respect to(w.r.t) $A$}
\bee\label{nonlinearwave'}
\bes
\square_{\tilde{g}(h)}\fai&=-\dfrac{2\eta'}{\eta^2}\frac{a}{(1+t)^{\lam}}\fai\de\fe+\dfrac{1}{\p \eta}\frac{a}{(1+t)^{\lam}}\left(\dfrac{\pa\fai}{\pa t}-\frac{\pa\fe}{\pa x^i}\frac{\pa \fai}{\pa x^i}-\frac{\lam}{(1+t)}\fai\right),
\end{split}
\ee
where $\til{g}$ is the conformal acoustical metric
\bee
\til{g}=\Omega g=\dfrac{\p}{\eta} g.
\ee
We also need to derive the linearized equation with respect to $A=\dfrac{\pa }{\pa t}$. The only difference is that in this case, $Ah=\dfrac{\pa \fai}{\pa t}-\dfrac{\pa\fe}{\pa x^i}\dfrac{\pa\fai}{\pa x^i}+\dfrac{a}{(1+t)^{\lam}}\fai-\dfrac{a\lam}{(1+t)^{\lam+1}}\fe.$ Then, it follows that for $A=\frac{\pa}{\pa t}$,
\bee\label{nonlinearwave}
\square_{\tilde{g}}\fai
=-\dfrac
{2\eta'}{\eta^2}A(\dfrac{a}{(1+t)^{\lam}}\fe)\de\fe+\dfrac{1}{\p\eta}A\left[\dfrac{a}{(1+t)^{\lam}}\left(\dfrac{\pa\fe}{\pa t}-\frac{1}{2}|\nabla\fe|^2-\frac{\lam}{(1+t)}\fe\right)\right].
\ee
\hs Note that \eqref{nonlinearwave} also holds for $A\in \{\pa_{i},\Omega_{ij}=x^i\pa_j-x^j\pa_i\}$ and when $\lam=0$, \eqref{nonlinearwave} is the corresponding linearized equation in the case of damping term $a\p v$.
\begin{remark}
\hs From above two examples, one can deduce that for general damping terms $a f(t)\p v$, where $f$ is a smooth function of $t$, the corresponding linearized equation for the variation $\fai$ of the potential function $\fe$ is
\[
\square_{\tilde{g}}\fai=-\dfrac{2\eta'}{\eta^2}A(af(t)\fe)\de\fe+
\dfrac{1}{\p\eta}A\left[\frac{\pa}{\pa t}(af(t)\fe)-\frac{1}{2}|\nabla\fe|^2af(t)\right].
\]
Hence, if $f$ has some decay such as $\dfrac{1}{(1+t)^{\lam}}$, then the same arguments in this paper hold.
\end{remark}
\begin{remark}
Note that for the metric $g$, the corresponding Christoffel symbols $\Ga_{\a\be\ga}$ satisfy
\bee\label{christoffelsymbol}
\Ga_{0ij}=\Ga_{ijk}=0,
\ee
which will be used later.
\end{remark}
\begin{remark}
Let $H=-2h-\eta^2$. Then for isentropic and irrotational flows, the condition \eqref{nonlinearcondition} is equivalent to
\bee
\dfrac{d H}{dh}\neq0,
\ee
which plays a key role in the shock formation.
\end{remark}
\subsection{\textbf{Initial data construction}}
\hs Let the initial data for \eqref{cee} be given on the hyperplane $\Si_{0}$ which is isentropic and irrotational. Denote $S_{0,a}$ to be the sphere on $\Si_{0}$ centered at origin with radius $1-a$. Outside the sphere $S_{0,0}$, set
\bee
\p=\p_0,\ s=s_0,\ v=0,\ \eta=\eta_0.
\ee
Then the initial data for the nonlinear wave equation \eqref{mass2} outside $S_{0,0}$ is given by
\bee
\fe=0,\ \pa_t\fe=h_0,
\ee
where $h_0$ is the initial enthalpy. For simplicity, one can set
\bee
\p_0=1,\eta_0=1,h_0=0.
\ee
Let the initial data be only non-trivial in the $\da$-annulus region:
\bee
\Si_{0}^{\da}:=\left\{x\in\Si_{0} |1-\da\leq r(x)\leq 1\right\}.
\ee
\hs Following the work in \cite{Ontheformationofshocks}, we construct the short pulse data for \eqref{mass2} in $\Si_{0}^{\da}$ as follows, whose proof is the same as Lemma1.1 in\cite{CZD1}.
\begin{lem}
For any given seed data $(\fe_1,\fe_2)$, there exists a $\da'>0$ depending only on $(\fe_1,\fe_2)$ such that for all $\da<\da'$, there exists another function $\fe_0\in C^{\infty}((0,1]\times S^2)$ smoothly depending on $(\fe_1,\fe_2)$ with the property that if one sets the initial data of \eqref{mass2} as follows:\\
\hs for $r\geq 1$ on $\Si_{0}$, $(\fe(0,x),\pa_t\fe(0,x))=(0,0)$; while for $1-\da\leq r\leq 1$,
\bee\label{initialdata}
\fe(0,x)=\da^{\frac{5}{2}}\fe_0\left(\dfrac{1-r}{\da},\ta\right),\hs \pa_t\fe(0,x)=
\da^{\frac{3}{2}}\fe_1\left(\dfrac{1-r}{\da},\ta\right).
\ee
Then, it holds that on $\Si_{0}$,
\bee\label{pat+parfe}
|(\pa_t+\pa_r)^2\fe|\les\da^{\frac{3}{2}}.
\ee
As a result, it holds that on $\Si_{0}$, $|(\pa_t+\pa_r)\pa_{\a}\fe|\les\da^{\frac{3}{2}}$ for $\a=0,1,2,3$. 
\end{lem}
\section{\textbf{The Geometric formulation, basic structure equations and the main results}}\label{section2}
\subsection{\textbf{Frames and coordinates}}
\hs In the following, we restrict the initial data on following the annular region:
\bee
\Si_{0}^{\da^3}:=\{x\in\Si_{0}:\ 1-\da^3\leq r(x)\leq 1\},
\ee
where $r=\sqrt{(x^1)^2+(x^2)^2+(x^3)^2}$.\\
\hs Define the function $u:=1-r$ on $\Si_0$. For each value of $u$, the corresponding level set $S_{0,u}$ is a sphere of radius $1-u$ and as a consequence,
\bee
\Si_{0}^{\tilde{\da}}=\bigcup_{u\in[0,\da^3]}S_{0,u}, \ \tilde{\da}:=\da^3.
\ee
\hs Consider, in the domain of maximal solution\footnote{For the notation of "The maximal development" and related notations, see\cite{CZD1} Section2.1.}, the family of outgoing characteristic null-hypersurface emanated from $S_{0,u}$, denoted to be $C_u$. Obviously, $C_u\cap\Si_{0}=S_{0,u}$, for any $u\in[0,\tilde{\da}]$.
By the domain of dependence, the solution outside $C_{0}$ is completely trivial and it's natural to consider the solution in the region foliated by $C_u$
\bee
W_{\tilde{\da}}=\bigcup_{u\in[0,\tilde{\da}]}C_u.
\ee
 Next. the function $u$ will be extended to $W_{\tilde{\da}}$ in the following definition, which implies that its level sets are precisely the outgoing characteristic null-hypersurfaces $C_u$.
\begin{defi}
The eikonal function $u$ satisifes
\bee\label{ekinoal}
\left\{\bes
&g^{\a\be}\dfrac{\pa u}{\pa x^{\a}}\dfrac{\pa u}{\pa x^{\be}}=0,\\
& u|_{\Si_{0}}=1-r,
\end{split}\right.
\ee
and $\pa_tu>0$.
\end{defi}
\begin{defi}
The outgoing null geodesic vector field is defined as
\bee
\hat{L}=-g^{\a\be}\pa_{\a}u\pa_{\be}.
\ee
\end{defi}
\hs It follows that $\hat{L}$ is $g-$null (since $g(\hat{L},\hat{L})=g^{\a\be}\pa_{\a}u\pa_{\be}u=0$) and $g-$orthogonal to $C_u$. Using $\hat{L}$, one can define the important inverse foliation density function as follows.
\begin{defi}
The inverse foliation density function is set to be:
\bee
\mu:=\dfrac{1}{-g^{\a\be}\pa_{\a}t\pa_{\be}u}=\dfrac{1}{\hat{L}^0},
\ee
whose reciprocal measures the density of the foliation $\bigcup_{u\in[0,\tilde{\da}]}C_u\cap\Si_t$.
\end{defi}
\hs Define for each $u\in[0,\tilde{\da}]$ and fixed $t$, the min of $\mu$ on the set $S_{t,u}$ to be $\mu(t,u)$. Then, define
\bee\label{defmumu}
\mu_m^u=\min\{\inf_{u'\in[0,u]}\mu(t,u'),1\},
\ee
and
\bee
s_{\ast}=\sup\{t|t\geq0\ \text{and}\ \mu_m^{\tilde{\da}}(t)>0\}.
\ee
\hs For each $u\in[0,\tilde{\da}]$, one can define $t_{\ast}(u)$ to be the lifespan of the solution to \eqref{ekinoal} and define $t_{\ast}$ to be:
\bee
t_{\ast}=\inf_{u\in[0,\tilde{\da}]}t_{\ast}(u)=\sup\{\tau|\text{smooth solution exists for all}\ (t,u)\in[0,\tau)\times[0,\tilde{\da}]\}.
\ee
We finally restrict time on $[0,t^{\ast})$ with
\bee
t^{\ast}=\min\{t_{\ast},s_{\ast}\}.
\ee
In the following, we work on $W_{\tilde{\da}}^{\ast}\subset W_{\tilde{\da}}$, where
\bee
W_{\tilde{\da}}^{\ast}=\bigcup_{(t,u)\in[0,t^{\ast})\times[0,\tilde{\da}]}S_{t,u},
\ee
where $S_{t,u}=C_u\cap\Si_t$ is a $2-$dimensional topological sphere. For the coordinate system on $S_{t,u}$, one can define it as follows.\\
\hs Since $S_{0,0}$ is the standard Euclidean sphere and each $S_{0,u}$ is diffeomorphic to $S_{-2,0}$, then if the local coordinates $(\til{\ta}^1,\til{\ta}^2)$ are chosen on $S_{0,0}$, the diffeomorphism defines a local coordinate system on $S_{0,u}$ for every $u\in[0,\tilde{\da}]$ which is denoted to be $(\ta_1,\ta_2)$. Extend $(\ta_1,\ta_2)$ to $S_{t,u}$ for $(t,u)\in[0,t_{\ast})\times[0,\tilde{\da}]$ by requiring that (for the definition of $L$, see\eqref{defiL})
\bee
L\ta^A=0,\ \text{for}\ A=1,2.
\ee
\hs The local coordinates $(\ta^1,\ta^2)$ together with $(t,u)$ define a complete local coordinates system $(t,u,\ta^1,\ta^2)$ on $W_{\tilde{\da}}^{\ast}$, which is called the acoustical coordinates.\par
Set:
\bee
\bes
\Si_t^u&:=\{(t,u',\ta)\in\Si_t|\ 0\leq u'\leq u\}\\
C_u^t&:=\{(t',u,\ta)\in C_u|\ 0\leq t'\leq t\}\\
W_t^u&:=\bigcup_{(t',u')\in[0,t)\times[0,u]} S_{t',u'}.
\end{split}
\ee
See also the following picture for the geometric notations.
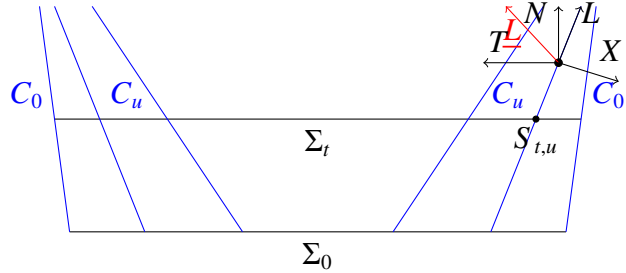
\begin{center}
\begin{tikzpicture}
\draw (-3.3,0)--node[below] {$\Si_0$}(3.3,0);
\draw (-3.5,1.5)--node[below] {$\Si_t$}(3.5,1.5);
\draw[blue] (1,0)--(3,3);
\draw[blue] (2.3,0)--node[above left] {$C_u$}(3.5,3);
\filldraw (2.9,1.5) circle (.04);
\node at (2.9,1.2) {$S_{t,u}$};
\draw[blue] (3.3,0)--node[above right] {$C_0$}(3.7,3);
\draw[blue] (-1,0)--(-3,3);
\draw[blue] (-2.3,0)--node[above right] {$C_u$}(-3.5,3);
\draw[blue] (-3.3,0)--node[above left] {$C_0$}(-3.7,3);
\filldraw (3.2,2.25) circle (.05);
\draw[->] (3.2,2.25)--node[above right] {$L$}(3.5,3);
\draw[->] (3.2,2.25)--node[above left] {$N$} (3.2,3);
\draw[red][->] (3.2,2.25)--node[left] {$\dl$}(2.5,3);
\draw[->] (3.2,2.25)--node[above left] {$T$}(2.2,2.25);
\draw[->] (3.2,2.25)--node[above right] {$X$}(4,2);
\end{tikzpicture}
\end{center}
\begin{defi}(Some important vector fileds)
\begin{itemize}
\item Define the outgoing null vector $L$ as
      \bee\label{defiL}
      L=\mu\hat{L},
      \ee
      which is $g-$orthogonal to $C_u$. Moreover, it holds that $Lt=1$.
\item Define on $W_{\tilde{\da}}^{\ast}$ a vector filed $T$ such that it is tangential to $\Si_t$ and $g-$orthogonal to $\{S_{t,u}\}$ for $u\in[0,\tilde{\da}]$, and normalized by
\bee
Tu=1.
\ee
\item Define the incoming null vector field $\dl$ as
      \bee
      \dl=\eta^{-2}\mu L+2T,
      \ee
      which is $g-$orthogonal to $S_{t,u}$.
\item Define the vector field
\bee
N=-\eta^2 g^{\a\be}\pa_{\a}t\pa_{\be},
\ee
which is $g-$orthogonal to $\Si_t$, future directed and timelike verifying $Nt=1$.
\item Let $\Lambda=[L,T]$ be the commutator of $L,T$ which is $S_{t,u}$ tangent. 
\item Define two $S_{t,u}-$tangent vector fields as $X_A=\dfrac{\pa}{\pa\ta^A}$ for $A=1,2$.
\end{itemize}
\end{defi}
\hs Then, one has the following properties of these vector fields, whose proofs are given in\cite{CZD1} Proposition2.1 and (2.24)-(2.25).
\begin{prop}\label{LTrelation}
\been[(1)]
\item It holds that
      \bee
      g(L,T)=-\mu,\quad g(L,\dl)=-2\mu,\quad g(T,T)=\kappa^2=(\eta^{-1}\mu)^2\ where\ \kappa>0.
      \ee
\item In the Cartesian coordinates, one has the following representations:
      \begin{align}
      L&=\dfrac{\pa}{\pa t}-\left(\eta\hat{T}^i+\dfrac{\pa\fe}{\pa x^i}\right)\dfrac{\pa}{\pa x^i},\ where\ \hat{T}=\kappa^{-1}T,\ |\hat{T}|=1.\label{Ldecomposition}\\
      N&=-\eta^2 g^{\a\be}\pa_{\a}t\pa_{\be}.
      \end{align}
\item In the acoustical coordinates, one has the following representations:
      \begin{itemize}
      \item $L=\dfrac{\pa}{\pa t}$.
      \item For the vector field $T$, one can decompose it as $T=\dfrac{\pa}{\pa u}-\Xi=\dfrac{\pa}{\pa u}-\Xi^A\dfrac{\pa}{\pa \ta^A}$. Moreover, it holds that $[L,\Xi]=-[L,T]=-\Lambda.$
      \end{itemize}
\een
\end{prop}

\hs Denote
\bee
\sg_{AB}=g\left(\dfrac{\pa}{\pa\ta^A},\dfrac{\pa}{\pa\ta^B}\right).
\ee
Then, $g$ has the following representations:
\bee
g=-2\mu dudt+\kappa^2 du^2+\sg_{AB}(d\ta^A+\Xi^Adu)(d\ta^B+\Xi^Bdu),
\ee
and
\bee\label{decomposeg2}
g^{\a\be}=-\dfrac{1}{2\mu}\left(L^{\a}\dl^{\be}+L^{\be}\dl^{\a}\right)+\sg^{AB}X_A^{\a}
X_B^{\be}.
\ee
\hs Consider the Jacobian of the coordinates transformation
\bee
(t,u,\ta^1,\ta^2)\ra(x_0,x_1,x_2,x_3).
\ee
Then, one has the following Fundamental theorem whose proof can be found in \cite{CZD1} Theorem 2.1.
\begin{thm}
\bee
\de=\det\dfrac{\pa(x_0,\cdots x_3)}{\pa(t,u,\ta^1,\ta^2)}=\eta^{-1}\mu\sqrt{\det\sg}.
\ee
Since $\eta\sim1$ (see Lemma\ref{hrhoeta}), the above relation implies that the diffeomorphism between the acoustical coordinates and the rectangular coordinates never degenerates \textbf{as long as $\mu>0$}, which implies that the frames $\{L,\dl,X_1,X_2\}$ and $\{L,T,X_1,X_2\}$ are equivalent to $\{\pa_{\a}\}_{\a=0,1,2,3}$ as long as $\mu>0$.
\end{thm}
\subsection{\textbf{$2^{nd}$ fundamental forms $k,\chi$ and $\ta$}}
\begin{defi}
Let $\bar{\mathcal{L}}$ and $\lie$ be the restriction of the lie derivative $\mathcal{L}$ on $\Si_t$ and $S_{t,u}$ respectively.
\been
\item The $2^{nd}$ fundamental form $k$ for the embedding $\Si_t^u\hookrightarrow W_{t}^{u}$ is defined as:
    \bee
    2\eta k=\bar{\mathcal{L}}_{N}g.
    \ee
    Moreover, $k$ can be computed as
    \bee\label{kij}
    2\eta k_{ij}=-2\pa_i\fai_j.
    \ee
\item The $2^{nd}$ fundamental form $\chi_{AB}$ for the embedding $S_{t,u}\hookrightarrow C_u^t$ is defined as:
    \bee
    2\chi_{AB}=\lie_L\sg_{AB}=2g(D_{A}L,X_B).
    \ee
\item The $2^{nd}$ fundamental form $\ta_{AB}$ for the embedding $S_{t,u}\hookrightarrow \Si_t^u$ is defined as:
    \bee
    2\eta^{-1}\mu\ta_{AB}=\lie_T\sg_{AB}=2g(D_AT,X_B).
    \ee
    Then, it follows from the definitions that $\chi=\eta(\s{k}-\ta)$ where $\s{k}$ is the restriction of $k$ on $S_{t,u}$.
\een
One can also define $\dc=\frac{1}{2}\lie_{\dl}\sg$ to be the $2^{nd}$ fundamental form for the embedding $S_{t,u}\hookrightarrow C_u^t$ associated to $\dl$.
\end{defi}
Set the following three $1-$forms $\ep,\zeta,\eta$ on $S_{t,u}$ as:
\bee
\eta^{-1}\mu\ep_A=k(X_A,T),\quad \zeta_A=g(D_AL,T),\quad \eta_A=-g(D_AT,L).
\ee
Then,
\bee\label{etaAcompute}
\zeta_A=\eta\kappa\ep_A-\kappa X_A\eta,\quad \eta_A= \zeta_A+X_A\mu.
\ee
\hs The following three structure equations will be used whose proofs are given in\cite{christodoulou2014compressible}.
\been[(1)]
\item The Gauss equation for the $2^{nd}$ fundamental form $\ta$ is:
      \bee\label{gaussta}
      \dfrac{1}{2}(tr\ta)^2-\dfrac{1}{2}|\ta|^2_{\sg}=K,
      \ee
      where $K$ is the Gauss curvature of $S_{t,u}$.
\item The Codazzi equations for the $2^{nd}$ fundamental form $\chi$ are:
      \bee\label{Codazzi1}
      \bes
      \s{D}_A\chi_{BC}-\s{D}_B\chi_{AC}&=R_{ABCL}-\mu^{-1}(\zeta_A\chi_{BC}-
      \zeta_B\chi_{AC})\\
      &=\beta_C\varepsilon_{AB}-\mu^{-1}(\zeta_A\chi_{BC}-
      \zeta_B\chi_{AC}),
      \end{split}
      \ee
      where $R$ is the Riemannian curvature tensor, $\beta$ is the $1-$form and $\varepsilon_{AB}$ is the area form of $\sqrt{\det\sg}$.\\
      \hs Contracting $\sg^{AC}$ on both sides of \eqref{Codazzi1} leads to
      \bee\label{Codazzi2}
      \s{div}\chi-\sd tr\chi=\beta^{\ast}-\mu^{-1}(\zeta\cdot\chi-\zeta\cdot
      tr\chi),
      \ee
      where $\beta^{\ast}$ is the $\sg$-dual of $\beta.$
\item $Ttr\chi$ satisfies the equation:
      \bee\label{Ttrchi}
      Ttr\chi=\s{D}^B\eta_B+\mu^{-1}\zeta\cdot\eta-\eta^{-1}L(\eta^{-1}\mu)tr\chi-\eta^{-1}\mu\ta\cdot\chi-\sg^{AB}R_{ATBL}.
      \ee
      It follows from \eqref{etaAcompute} and \eqref{Ttrchi} that $Ttr\chi-\s{\de}\mu\sim\mu\sd^2\fai+l.o.ts.$
\een
\hs Noticing
\bee\label{decomposewavenonconformal}
\square_gf=\s{\de}f-\frac{1}{2}\mu^{-1}tr\dc(Lf)-\frac{1}{2}\mu^{-1}tr\chi(\dl f)-\mu^{-1}L(\dl f)-2\mu^{-1}\zeta\cdot\sd f,
\ee
and
\bee\label{relationconfromal}
\square_{\til{g}}f=\Omega^{-1}\square_gf+\Omega^{-2}\dfrac{d\Omega}{dh}g^{\a\be}\pa_{\a}h\pa_{\be}f,
\ee
yields
\bee\label{decomposewaveeq}
\bes
\square_{\til{g}}f&=\Omega^{-1}\s{\de}f-\Omega^{-1}\mu^{-1}L(\dl f)+2\Omega^{-1}\mu^{-1}\zeta\cdot\sd f\\
&-\dfrac{1}{2}\mu^{-1}\Omega^{-1}\left((tr\dc+\dfrac{d\log\Omega}{dh}\dl h)Lf+
(tr\chi+\dfrac{d\log\Omega}{dh}Lh)\dl f\right).
\end{split}
\ee
\subsection{\textbf{Curvature tensor}}
\hs Define\footnote{$R_{\a\be\ga\da}$ is the Riemannian curvature tensor of $g$. In our case, the non-vanishing principle part (the second derivatives of $g$) of the component of the curvature tensor is $R_{0i0j}$ with its principle part given by $R_{0i0j}^{[P]}=-\dfrac{1}{2}\dfrac{dH}{dh}D^2_{ij}h.$}
\bee
\a_{AB}=R_{\a\be\ga\da}X_A^{\a}L^{\be}X_B^{\ga}L^{\da}=R_{ijkl}X_A^iL^jX_B^kL^l+R_{i0jk}X_A^iX_B^jL^k+R_{ijk0}X_A^iL^jX_B^k+
R_{0i0j}X_A^iX_B^j.
\ee
\hs Then, the principle part of $\a_{AB}$ is given by
\bee
\a_{AB}^{[P]}=R_{0i0j}^{[P]}X_A^iX_B^j=-\dfrac{1}{2}\dfrac{dH}{dh}D^2_{AB}h=-\dfrac{1}{2}\dfrac{dH}{dh}\left(\s{D}^2_{AB}h-\eta^{-1}\s{k}_{AB}Lh-\mu^{-1}
\chi_{AB}Th\right).
\ee
Denote
\bee
\a_{AB}'=\a_{AB}-\dfrac{1}{2}\mu^{-1}\dfrac{dH}{dh}\chi_{AB}Th=-\dfrac{1}{2}\dfrac{dH}{dh}\left(\s{D}^2_{AB}h-\eta^{-1}\s{k}_{AB}Lh\right)+\a_{AB}^{[N]},
\ee
with $\a_{AB}^{[N]}$ given by
\bee\label{a_AB^N}
\bes
&X_A^iX_B\fai_iL^{\a}L\fai_{\a}-L^{\a}L^{\be}X_A\fai_{\a}X_B\fai_{\be}+X_AhX_Bh\\
&+2LhX_A^iX_B\fai_i-X_AhL^{\a}X_B\fai_{\a}-X_BhL^{\a}X_A\fai_{\a},
\end{split}
\ee
which contains at most first order derivatives of $g$.
\subsection{\textbf{The rotation vector fields}}
\begin{defi}
Let $\mari{R}_i$ be the standard space rotation in the Euclidean space such that $\mari{R}_i=\ep_{ijk}x^j\frac{\pa}{\pa x^k}=\frac{1}{2}\ep_{ijk}(x^j\pa_k-x^k\pa_j)$ where $\ep_{ijk}$ is the skew-symmetric symbol. Then, define the rotation vector fields as follows:
\bee
R_i=\Pi\mari{R}_i, \ i=1,2,3,
\ee
where $\Pi$ is the projection from $\Si_t$ to $S_{t,u}$. Precisely,
\bee
\Pi_{i}^j=\da_i^j-\bar{g}_{ik}\hat{T}^k\hat{T}^j=\da_i^j-\hat{T}^i\hat{T}^j.
\ee
\end{defi}
Let $\{Z_i\}$ with $i=1,2,3,4,5$ be the set of the commutation vectorfields such that
\bee
Z_1=Q=(1+t)L,\ Z_2=T,\ Z_{i+2}=R_i,\ \text{for}\ i=1,2,3.
\ee
\subsection{\textbf{Transport equations for $\mu$ and $\chi$}}
\begin{prop}
$\mu$ and $\chi$ satisfy the following transport equations respectively\footnote{For the computation, see\cite{CZD1} Proposition2.2.}.
\begin{align}
L\mu&=m+\mu e,\label{transportmu}\\
\begin{split}\label{transportchi}
L\chi_{AB}&=(\mu^{-1}L\mu)\chi_{AB}+\chi_A^C\chi_{BC}-\a_{AB}\\
&=e\chi_{AB}+\chi_{A}^C\chi_{BC}-\a'_{AB}+a\eta^{-1}\dfrac{\hat{T}\fe}{(1+t)^{\lam}}\chi_{AB},
\end{split}\\
Ltr\chi&=(\mu^{-1}L\mu)tr\chi-|\chi|^2_{\sg}-tr\a.\label{transporttrchi}
\end{align}
with
\begin{align}
m&=\dfrac{1}{2}\dfrac{dH}{dh}Th+a\dfrac{T\fe}{(1+t)^{\lam}},\ H=-2h-\eta^2,\label{mcomputation}\\
e&=\dfrac{1}{2\eta^2}\left(\dfrac{\p}{\p'}\right)'Lh+\eta^{-1}\hat{T}^i(L\fai_i).\label{edefinition}
\end{align}
\end{prop}
\begin{remark}
\been[(1)]
\item Note that the right hand side of both equations are regular in $\mu$, that is, when $\mu\ra 0$, the quantities on the right hand side of \eqref{transportmu}, \eqref{transportchi} and \eqref{transporttrchi} do not blow up.
\item If the assumption \eqref{nonlinearcondition} violates, i.e. $\frac{dH}{dh}=0$, then $m=a\dfrac{T\fe}{(1+t)^{\lam}}$ and it can be seen later that $\mu$ never vanishes and \eqref{cee} admits a global solution.
\een
\end{remark}
\subsection{\textbf{The Main Result}}
\hs At the end of this section, we introduce the main theorem in this paper.
\begin{thm}
For the nonlinear wave equation \eqref{mass2}, 
 let the initial data $(\fe,\pa_t\fe)$ be given as in Section 1.6.\\
\hs Then, in the acoustical coordinates, the bootstrap assumptions \eqref{bs} hold on $W_{\tilde{\da}}^{\ast}$ and $\fai,\mu,\chi,g$ together with their derivatives are regular on $\Si_{t^{\ast}}$.\\
\hs In the rectangular coordinates, we have the following two possibilities:
\been[(1)]
\item If the following largeness condition is satisfied:\textbf{ $\min\pa_r\fe_1\leq-1$}, then there exists $\da_0$ such that for all $\da<\da_0$, the damping can not prevent the shock formation to the nonlinear wave equations \eqref{nonlinearwave}, which occurs before $T_{\ast}=\exp\left(\dfrac{\da^{-\frac{1}{2}}C_{a,\lam}}{2|\pa_r\fe_1|}\right)-1$, where $C_{a,\lam}$ is a constant depending on $a,\lam$ bounded as $e^{-|\frac{a}{1-\lam}|}\leq C_{a,\lam}\leq e^{|\frac{a}{1-\lam}|}$ and $C_{a,\lam}>(<)1$ if $a>(<)0$.
\item If the initial data satisfies: \textbf{$\min\pa_{r}\fe_1=c>0$}, then there exists $\da_0$ such that for all $\da<\da_0$, smooth solutions to the nonlinear wave equations \eqref{nonlinearwave} exist on $[0,+\infty)$.
\een
\hs Moreover, in case $(1)$, the inverse foliation density $\mu$ vanishes at some points on $\Si_{t^{\ast}}$ at which the outgoing characteristic null hypersurfaces $C_u$ become infinitely dense. More precisely, $\pa_i\pa_j\fe$ blows up at these points where $\fe$, $\pa\fe$ remain bounded. That is, for the compressible Euler equations, $\pa_iv^j$ and $\pa_i\p$ blow up before $t=T_{\ast}$ where $v^i$ and $\p$ remain bounded. Furthermore, the damping effect will affect the formation of shock in the following sense:
\begin{itemize}
\item for fixed $\lam$, if $a>0$, then the damping effect will delay the shock formation; while if $a<0$, then the anti-damping term will lead to the formation of shock in advance.
\item For fixed $a$, as $\lam $ decreasing, the damping effect becomes stronger which delays the formation of shock, and vice versa.
\end{itemize}
\hs In case $(2)$, we obtain a global solution to the compressible Euler equations \eqref{cee} on $[0,+\infty)$.
\end{thm}
\section{\textbf{Bootstrap Assumption and preliminary estimates}}\label{section3}
\hs Let $\fe$ be the solution to \eqref{mass2} and $\fai$ be its variations which are guaranteed by the local well-posedness theory. By the choice of the initial data, it holds that on $\Si_{0}$
\begin{align}
&|L\fai|=|(\pa_t+\pa_r)\fai|\les\da^{\frac{3}{2}},\ |T\fai|=|\pa_r\fai|\les \da^{\frac{1}{2}},\ |\sd\fai|\les\da^{\frac{3}{2}},\\
&\da^m|R_{i_n}\cdots R_{i_1}(T)^mQ^p\fai|\les \da^{\frac{3}{2}-m},
\end{align}
where $p\leq2$. Then, assume that the following bootstrap assumptions hold for all $t\in[0,t^{\ast})$.
\been[(1)]
\item  For $2\leq|\a|\leq N_{\infty}:=[\dfrac{N_{top}}{2}]+3$,
\bee\label{bs}
\max_{\a}\max_{i_1,...i_n}\da^m||R_{i_n}\cdots R_{i_1}(T)^mQ^p\fai||_{\supnormda}+||L\fai||_{\supnormda}+\da||T\fai||_{\supnormda}+(1+t)||\sd\fai||_{\supnormda}\les \da^{\frac{3}{2}}M(1+t)^{-1},
\ee
for $p+m+n=|\a|$ with $p\leq 2$;
\item
\bee\label{bs2}
||T\fe||_{\supnormda}\les \da^{\frac{1}{2}}M,
\ee
\een
where $M$ is a large constant independent of $\da$ to be determined later and $N_{top}$ is an integer measuring the totally derivatives acting on $\fai$. For $M$ large enough, it follows immediately that $|\fe|\leq|\int_{0}^t\fai_0 ds|+|O(\da^{\frac{5}{2}})|\les\da^{\frac{3}{2}}\ln(1+t)M$\footnote{Later one can see that as shock forms, $\ln(1+t)\sim\da^{-\frac{1}{2}}$, which  implies $\fe$ remains bounded(small).}.
\subsection{\textbf{Preliminary Estimate for the metric, $\mu$ and the second fundamental forms under the Bootstrap assumptions}}
Using the same arguments in \cite{CZD1}, one obtains the following three lemmas.
\begin{lem}\label{hrhoeta}
Under the bootstrap assumptions \eqref{bs} and \eqref{bs2}, the following estimates hold for sufficiently small $\da$:
\begin{equation*}
\bes
&||h||_{\supnormda}\les\dfrac{\da^{\frac{3}{2}} M}{1+t},\quad||\p-\p_0||_{\supnormda}\les\dfrac{\da^{\frac{3}{2}} M}{1+t},\quad ||\eta-\eta_0||_{\supnormda}\les\dfrac{\da^{\frac{3}{2}} M}{1+t},\\
&||L\fe||_{\supnormda}\les\dfrac{\da^{\frac{3}{2}} M}{1+t},\hs ||\sd\fe||_{\supnormda}\les\dfrac{\da^{\frac{3}{2}} M}{1+t}.
\end{split}
\end{equation*}
\hs As a consequence, $C^{-1}\leq\Omega\leq C$ for some constant $C$.
\end{lem}

\hs In the following, the bounds $\p\sim 1$ and $\eta,\eta'\sim 1$ will be used without mention.
\begin{lem}\label{esthme}
Under the bootstrap assumptions \eqref{bs} and \eqref{bs2}, the following estimates hold for sufficiently small $\da$:
\beeq
&&||L(h,\p,\eta)||_{\supnormda}+||\sd (h,\p,\eta)||_{\supnormda}+\da(1+t)^{-1}||T(h,\p,\eta)||_{\supnormda}\les\dfrac{\da^{\frac{3}{2}} M}{(1+t)^2},\\
&&||m||_{\supnormda}\les\dfrac{\da^{\frac{1}{2}} M}{1+t}, \hs||e||_{\supnormda}\les\dfrac{\da^{\frac{3}{2}} M}{(1+t)^2}.
\eeq
\end{lem}
\begin{lem}\label{estmu}
Under the same assumptions as in Lemma\ref{esthme}, it holds that for sufficiently small $\da$
\begin{align*}
||\mu-1||_{\supnormda}&\les \da^{\frac{1}{2}}M(1+\ln(1+t)),\\
||L\mu||_{\supnormda}&\les\dfrac{\da^{\frac{1}{2}}M}{1+t}.
\end{align*}
As a consequence, one can obtain $||T\fe||_{\supnormda}\les \da^{\frac{3}{2}}M^2$. which recovers the bootstrap assumption \eqref{bs2} for sufficiently small $\da$.
\end{lem}
Let $\dfrac{\sg_{AB}}{1-u+t}$ be the null $2^{nd}$ fundamental form in Minkowski space. Then set
\bee
\chi'_{AB}=\chi_{AB}-\dfrac{\sg_{AB}}{1-u+t},\quad \ta'_{AB}=\ta_{AB}+\dfrac{\sg_{AB}}{1-u+t}.
\ee
\begin{lem}\label{2ndff}
Under the bootstrap assumptions, it holds that for sufficiently small $\da$
\beeq
&||\s{k}||_{\supnormda}\les\dfrac{\da^{\frac{3}{2}} M}{(1+t)^2},&\\
&||\chi'||_{\supnormda}\les\dfrac{\da^{\frac{3}{2}} M(1+\ln (1+t))}{(1+t)^2},&\\
&||\ta'||_{\supnormda}\les\dfrac{\da^{\frac{3}{2}} M(1+\ln (1+t))}{(1+t)^2}.&
\eeq
\end{lem}
\begin{pf}
Since $\s{k}_{AB}=X_A^iX_B^jk_{ij}=-\eta^{-1}X_A^iX_B^j\pa_i\fai_j$, it follows that
\beeq
|\s{k}|^2&=&\sg^{AC}\sg^{BD}\s{k}_{AB}\s{k}_{CD}=\eta^{-2}\sg^{AC}\sg^{BD}X_A^i
X_B^j\pa_i\fai_jX_C^kX_D^l\pa_k\fai_l\\
&\leq&\eta^{-2}\bar{g}^{ik}\sg^{BD}X_B^j\pa_j\fai_iX_D^l\pa_l\fai_k=
\eta^{-2}\sum_i\sg^{BD}X_B^j\pa_j\fai_iX_D^l\pa_l\fai_i\\
&\leq&\eta^{-2}\sum_i|\sd\fai_i|^2\les \dfrac{\da^3 M^2}{(1+t)^4}.
\eeq
As a consequence, it holds that
\[
|\sd^2\fe|^2=\sg^{AC}\sg^{BD}X_A^i
X_B^j\pa_i\fai_jX_C^kX_D^l\pa_k\fai_l\les\dfrac{\da^3 M^2}{(1+t)^4}.
\]
\hs It follows from \eqref{transportchi} and $L\dfrac{\sg_{AB}}{1-u+t}=-\dfrac{2\chi'_{AB}}{1-u+t}+\dfrac{\sg_{AB}}{(1-u+t)^2}$ that
\bee\label{Lchi'AB}
L\chi'_{AB}=e\chi'_{AB}+\chi_A^{'\ C}\chi'_{BC}-\a'_{AB}+a\frac{\eta^{-1}\hat{T}\fe}{(1+t)^{\lam}}\chi'_{AB}+\left(e+a\frac{\eta^{-1}\hat{T}\fe}{(1+t)^{\lam}}\right)\dfrac{\sg_{AB}}{1-u+t}.
\ee
Therefore,
\begin{align*}
&L(|\chi'|^2_{\sg})=L\left(\sg^{AC}\sg^{BD}\chi'_{AB}\chi'_{CD}\right)\\
&=(2e+a\frac{2\eta^{-1}\hat{T}\fe}{(1+t)^{\lam}})|\chi'|^2_{\sg}-\dfrac{2}{1-u+t}|\chi'|^2-2\chi^{'\ AC}\chi_A^{'\ D}\chi'_{CD}
-2\chi'^{AB}\a'_{AB}+\dfrac{2}{1-u+t}tr\chi'\left(e+a\frac{\eta^{-1}\hat{T}\fe}{(1+t)^{\lam}}\right),
\end{align*}
which implies that
\bee\label{lt-uchi'}
L((1-u+t)^2|\chi'|_{\sg})\leq(1-u+t)^2\left((|e|+|\frac{\hat{T}\fe}{(1+t)^{\lam}}|+|\chi'|)|\chi'|+|\a'|+|\dfrac{1}{1-u+t}|(|e|+|\frac{\hat{T}\fe}{(1+t)^{\lam}}|)\right).
\ee
It remains to estimate $\a'_{AB}$. Due to the definition $\a'_{AB}=-\dfrac{1}{2}\dfrac{dH}{dh}\s{D}^2_{A,B}h+\dfrac{1}{2}\dfrac{dH}{dh}
\eta^{-1}\s{k}_{AB}Lh+\a_{AB}^{[N]}$ and \eqref{a_AB^N}, where the second term and the last term are bounded by $\dfrac{\da^{\frac{3}{2}}M}{(1+t)^3}$, it remains to estimate $\s{D}^2_{A,B}h$. It follows from Proposition\ref{relationangular} that
\bee
|\s{D}^2_{A,B}h|\les|\s{d}^2 h|+|\s{d}h|\les|\sd^2\fai_0|+|\sd(\fai_i\sd\fai_i)|+|\sd^2\fe|\les\frac{\da^{\frac{3}{2}} M}{(1+t)^3},
\ee
 which implies that $|\a'|\les\da \dfrac{\da^{\frac{3}{2}} M}{(1+t)^3}$.\\
\hs Let $P(t)$ be the set of $t'$ that $||\chi'||_{\supnormda}\leq C_0\dfrac{\da^{\frac{3}{2}} M(1+\ln(1+t))}{(1+t)^2}$ for all $t'\in[0,t]$ with sufficiently small $\da$, where $C_0$ is a constant to be determined later. Initially, $|\chi_{AB}|=|\eta(\s{k}_{AB}-\ta_{AB})|\leq|\eta(\sd\fai-\dfrac{\sg_{AB}}{2})|$. It follows from the construction of initial data and choosing $C_0$ large enough that $||\chi'||_{L^{\infty}(\Si_{0}^{\da})}<C_0\da^{\frac{3}{2}}<C_0\da^{\frac{3}{2}}M$. That is, $0\in P(t)$.\\
 \hs Let $t_0$ be the upper bound of $P(t)$. For $t\leq t_0$, it holds that $|e|+|\frac{\hat{T}\fe}{(1+t)^{\lam}}|+3|\chi'|\leq (C_0+C_1)\dfrac{\da^{\frac{3}{2}} M(1+\ln(1+t))}{(1+t)^2}$ and $|\a'|\leq C_2\dfrac{\da^{\frac{3}{2}} M(1+\ln(1+t))}{(1+t)^3}$ for some constants $C_1$ and $C_2$. Thus, it follow from \eqref{lt-uchi'} that
 \bee\label{L(t-u)^2chi}
 L((1-u+t)^2|\chi'|)\leq(1-u+t)^2\left[(C_0+C_1)\dfrac{\da^{\frac{3}{2}} M(1+\ln(1+t))}{(1+t)^2}|\chi'|+(C_2+C_1)\da^{\frac{3}{2}}M(1+t)^{-1}\right].
 \ee
 Let $x=(1-u+t)^2|\chi'|$. Then, along integral curves of $L$, one can rewrite \eqref{L(t-u)^2chi} as
\[
\dfrac{dx}{dt}\leq fx+g,
\]
where $f=(C_0+C_1)\dfrac{\da^{\frac{3}{2}} M(1+\ln(1+t))}{(1+t)^2}$ and $g=(C_2+C_1)\da^{\frac{3}{2}}M(1+t)^{-1}$. Integrating from $0$ to $t$ and noting that $(1-u+t)^2\sim (1+t)^2$ yield
\[
(1+t)^2|\chi'|_{\sg}\leq e^{(C_0+C_1)\da^{\frac{3}{2}} M}\left(||\chi'||_{L^{\infty}(\Si_{0}^{\da})}+C_3(C_1+C_2)\da^{\frac{3}{2}} M\ln (1+t)\right)
\]
Choose $C_0$ such that $C_0>4C_3(C_1+C_2)$. Then, it holds that $(1+t)^2||\chi'||_{\supnormda}<C_0\da^{\frac{3}{2}} M(1+\ln(1+t))$ for $\da^{\frac{3}{2}}<\dfrac{\log 2}{(3C_0+C_1)M}$. Hence, by a continuous argument, one shows that $P(t)$ holds for all $t\in[0,t^{\ast})$.\\
\hs The estimates for $\ta'$ follows from $|\ta'_{AB}|=|\s{k}_{AB}-\eta^{-1}\left(\chi_{AB}-\dfrac{\eta\sg_{AB}}{1-u+t}\right)
|\les\dfrac{\da^{\frac{3}{2}} M\ln (1+t)}{(1+t)^2}$.
\end{pf}
\hs Based on the above estimates, one has the following two lemmas, whose proofs are similar to Lemma3.5 and 3.6 in \cite{CZD1}.
\begin{lem}\label{LTidTi}
It holds that
\begin{align}
L\hat{T}&=(\hat{T}^j(X^A\fai_j)+X^A(\eta))X_A,\label{LTi}\\
X_A\hat{T}&=\ta_A^BX_B,\\
T\hat{T}&=-X^A(\kappa)X_A.
\end{align}
Therefore, the following bounds hold for sufficiently small $\da$:
\begin{align}
|L\hat{T}^i|&\les|\hat{T}^j||X^A\fai_j|+|X^A\eta|\les\dfrac{\da^{\frac{3}{2}} M}{1+t},\\
|\sd\hat{T}^i|&\les|\ta|\les \dfrac{M}{1+t}.\label{sdhatTi}
\end{align}
\end{lem}
\begin{lem}\label{ltfe}
The following estimates hold for sufficiently small $\da$:
\begin{align}
||LT\fe||_{\supnormda}&\les \dfrac{\da^{\frac{3}{2}} M}{(1+t)^2}.\\
||T^2\fe||_{\supnormda}&\les\da^{-\frac{1}{2}}M,\label{T^2fe}\\
||\sd T\fe||_{\supnormda}&\les\dfrac{\da^{\frac{1}{2}} M}{1+t}.
\end{align}
Meanwhile, the following estimates on $\mu$ hold:
\bee
||\sd\mu||_{\supnormda}\les \dfrac{\da^{\frac{1}{2}}M(1+\ln(1+t))}{1+t},\ ||T\mu||_{\supnormda}\les\da^{-\frac{1}{2}} M(1+\ln(1+t)).
\ee
\end{lem}
\subsection{\textbf{Accurate estimate for $\mu$ and its derivatives}}
\hs Considering the inhomogeneous covariant wave equation $\square_{\til{g}}\fai_{\a}=-\dfrac
{2\eta'}{\eta^2}A(\dfrac{a}{(1+t)^{\lam}}\fe)\de\fe+\dfrac{1}{\p\eta}A\left[\dfrac{a}{(1+t)^{\lam}}\left(\dfrac{\pa\fe}{\pa t}-\dfrac{1}{2}|\nabla\fe|^2-\dfrac{\lam}{(1+t)}\fe\right)\right]$ and noticing that \eqref{decomposewaveeq} can be written as
\bee\label{decompose}
\bes
\mu\Omega\square_{\til{g}}f&=\mu\s{\de}f-L(\dl f)-\dfrac{1}{2}\left(
(\til{tr}\chi+L\log\Omega)\dl f+(\til{tr}\dc+\dl\log\Omega) Lf\right)\\
&-2\zeta\cdot\sd f+\Omega^{-1}\dfrac{d\Omega}{dh}\sd h\cdot\sd f
\end{split}
\ee
yield the following transport equation for $\dl\fai_{\a}$:
\bee\label{transportdlfai}
\bes
L(\dl \fai_{\a})+\dfrac{1}{2}\til{tr}\chi\dl \fai_{\a}&=\underbrace{\mu\s{\de}\fai_{\a}-\frac{1}{2}L(\log\Omega)\dl\fai_{\a}-\dfrac{1}{2}(\til{tr}\dc+\dl\log\Omega)L\fai_{\a}-2\zeta\cdot\sd \fai_{\a}+\Omega^{-1}\dfrac{d\Omega}{dh}\sd{h}\cdot\sd \fai_{\a}}_{I}\\
&\underbrace{\dfrac
{2\p\mu\eta'}{\eta^3}A(\dfrac{a}{(1+t)^{\lam}}\fe)\de\fe-\dfrac{\mu}{\eta^2}A\left[\dfrac{a}{(1+t)^{\lam}}\left(\dfrac{\pa\fe}{\pa t}-\frac{1}{2}|\nabla\fe|^2-\frac{\lam}{(1+t)}\fe\right)\right]}_{II},
\end{split}
\ee
where $A$ is the variation operator such that $A\fe=\fai_{\a}$. It follows from Lemma\ref{estmu},\ref{2ndff}, the bootstrap assumptions \eqref{bs}, Proposition\ref{relationangular} and the relation $tr\dc=\eta^{-2}\mu tr\chi+2\eta^{-1}\mu tr\ta$ that $I$ is bounded by $\dfrac{\da^{\frac{3}{2}}M}{(1+t)^3}$. It remains to estimate $II$. For $\mu\de\fe$, it holds that
\bee
\bes
|\mu\de\fe|&=|\mu\da^{ij}\frac{\pa\fai_j}{\pa x^i}|=|\mu(\hat{T}^i\hat{T}^j+\sg^{AB}X_A^iX_B^j)\frac{\pa\fai_j}{\pa x^i}|\\
&\les|\hat{T}^iT\fai_i|+\mu|\sd\fai|\les \dfrac{\da^{\frac{1}{2}}M}{1+t}.
\end{split}
\ee
\begin{itemize}
\item For $A=\frac{\pa}{\pa x^i}$, it holds that $|A(\dfrac{a}{(1+t)^{\lam}}\fe)|\les\dfrac{\da^{\frac{3}{2}}M}{(1+t)^{\lam+1}}$;
\item while for $A=\frac{\pa}{\pa t}$, it holds that $|A(\dfrac{a}{(1+t)^{\lam}}\fe)|\les\dfrac{\da^{\frac{3}{2}}M\ln(1+t)}{(1+t)^{\lam+1}}$.
\end{itemize}
Thus, it holds that $|\frac{2\p\mu\eta'}{\eta^3}A(\dfrac{a}{(1+t)^{\lam}}\fe)\de\fe|\les\dfrac{\da^{\frac{3}{2}}M}{(1+t)^3}$. Note that $\mu\dfrac{\pa\fai_{\a}}{\pa t}-\mu\fai_i\dfrac{\pa\fai_{\a}}{\pa x^i}=\mu L\fai_{\a}+\eta^{2}T\fai_{\a}+\mu\fai_i\dfrac{\pa\fai_{\a}}{\pa x^i}-\mu\fai_i\dfrac{\pa\fai_{\a}}{\pa x^i}=\frac{1}{2}\eta^{2}\dl\fai_{\a}+\frac{1}{2}\mu L\fai_{\a}$, then
\begin{itemize}
\item for $A=\frac{\pa}{\pa x^i}$, it holds that
\bee\label{inhomogeneous1}
\mu A\left[\dfrac{a}{(1+t)^{\lam}}\left(\dfrac{\pa\fe}{\pa t}-\frac{1}{2}|\nabla\fe|^2-\frac{\lam}{(1+t)}\fe\right)\right]
=\dfrac{a}{(1+t)^{\lam}}\left(\frac{1}{2}\eta^2\dl\fai_i+\frac{1}{2}\mu L\fai_i-\frac{\lam}{1+t}\fai_i\right);
\ee
\item while for $A=\frac{\pa}{\pa t}$, it holds that
\bee\label{inhomogeneous2}
\bes
&\mu A\left[\dfrac{a}{(1+t)^{\lam}}\left(\dfrac{\pa\fe}{\pa t}-\frac{1}{2}|\nabla\fe|^2-\frac{\lam}{(1+t)}\fe\right)\right]\\
=&\dfrac{a}{(1+t)^{\lam}}\left(\frac{1}{2}\eta^2\dl\fai_0+\frac{1}{2}\mu L\fai_0-\frac{\lam}{1+t}\fai_0\right)-\dfrac{a\lam}{(1+t)^{\lam+1}}\left(\fai_0-\frac{1}{2}|\fai_i|^2-\frac{\lam}{1+t}\fe\right).
\end{split}
\ee
\end{itemize}
Thus, \eqref{transportdlfai} can be written as
\bee\label{transportdlfai2}
L(\dl\fai_{\a})+\dfrac{1}{2}\left(\til{tr}\chi+\dfrac{a}{(1+t)^{\lam}}\right)
\dl\fai_{\a}=O(\dfrac{\da^{\frac{3}{2}} M}{(1+t)^3}).
\ee
For convenience, we may replace $\a$ by $2a$ in the analysis.
\begin{prop}\label{accruatemu1}
For sufficiently small $\da$, it holds that
\bee\label{aestlmu}
|(1+t)A(t) L\mu(t,u,\ta)-L\mu(0,u,\ta)|\les\da^{\frac{3}{2}}M\left(1-\dfrac{1}{1+t}\right),
\ee
\bee\label{aestmu}
|\mu(t,u,\ta)-1-L\mu(0,u,\ta)\int_0^t\dfrac{1}{(1+\tau)A(\tau)}d\tau|\les\da^{\frac{3}{2}}M\left(
1+\ln(1+t)-\dfrac{1}{1+t}\right),
\ee
where
$A(t)=e^{\int_0^t\frac{a}{(1+s)^{\lam}}ds}=e^{\frac{a}{\lam-1}\left(1-\frac{1}
{(1+t)^{\lam-1}}\right)}$.
\end{prop}
\begin{pf}
\ Since $\dfrac{1}{2}\til{tr}\chi=\dfrac{1}{2}\til{tr}\chi'+\dfrac{1}{1-u+t}$, where $|\til{tr}\chi'|\les\dfrac{\da^{\frac{3}{2}} M(1+\ln (1+t))}{(1+t)^2}$ due to Lemma\ref{2ndff}, it follows from\eqref{transportdlfai2} that
\[
L(\dl\fai_{\a})+\left(\dfrac{1}{1-u+t}+\dfrac{a}{(1+t)^{\lam}}\right)\dl\fai_{\a}
=O(\dfrac{\da^{\frac{3}{2}}M}{(1+t)^3}),
\]
which implies
\bee\label{Leatu-tdlfai}
|L((1-u+t)A(t)\dl\fai_{\a})|\les\dfrac{\da^{\frac{3}{2}}M}{(1+t)^2}A(t)\les
\dfrac{\da^{\frac{3}{2}}M}{(1+t)^2},
\ee
where $A(t)=e^{\int_0^t\frac{a}{(1+s)^{\lam}}ds}=e^{\frac{a}{\lam-1}\left(1-\frac{1}
{(1+t)^{\lam-1}}\right)}$, which is uniformly bounded below and above by a constant $C(a,\lam)$.\\
\hs Integrating \eqref{Leatu-tdlfai} along integral curves of $L$ and using the bootstrap assumptions yield
\bee\label{estdlfai}
\bes
& |(1+t)A(t)\dl\fai_{\a}(t,u,\ta)- \dl\fai_{\a}(0,u,\ta)|\les\da^{\frac{3}{2}}M\left(1-\dfrac{1}{1+t}\right),\\
 & |(1+t)A(t) T\fai_{\a}(t,u,\ta)-  T\fai_{\a}(0,u,\ta)|\les\da^{\frac{3}{2}}M\left(1-\dfrac{1}{1+t}\right),\\
&|(1+t)A(t)\fai_{\a}(t,u,\ta)- \fai_{\a}(0,u,\ta)|
\les\da^{\frac{5}{2}}M\left(1-\dfrac{1}{1+t}\right).
\end{split}
\ee
 It follows from $|T\fe|\les\da^{\frac{3}{2}}M$ that
\bee\label{esttfe}
|(1+t)A(t)\dfrac{T\fe(t,u,\ta)}{(1+t)^{\lam}}-  T\fe(0,u,\ta)|\les\da^{\frac{3}{2}}M\left(1-\dfrac{1}{1+t}\right).
\ee
\hs It follows from \eqref{transportmu} that
\bee\label{lmu-2t}
\bes
(1+t)A(t)L\mu(t,u,\ta)- L\mu(0,u,\ta)&=(1+t)A(t)m(t,u,\ta)-m(0,u,\ta)\\
&+(1+t)A(t)\mu e|_{(t,u,\ta)}-\mu e|_{(0,u,\ta)}\\
&=(1+t)A(t)m(t,u,\ta)-m(0,u,\ta)+O(\da^2 M^2),
\end{split}
\ee
where the last equality is due to Lemma\ref{esthme} and \ref{estmu}. It follows from \eqref{mcomputation}, Lemma\ref{hrhoeta}
, \eqref{estdlfai} and \eqref{esttfe} that
\bee\label{m-2t}
\bes
|(1+t)A(t)m(t,u,\ta)-m(0,u,\ta)|&=|(1+\tau)A(\tau)\left[(-1-\eta'\eta)(T\fai_0-
\fai_iT\fai_i+\dfrac{T\fe}{(1+t)^{\lam}})\right]|^t_0\\
&+a\dfrac{(1+\tau)A(\tau)T\fe}{(1+\tau)^{\lam}}|^t_0\les \da^{\frac{3}{2}}M\left(1-\dfrac{1}{1+t}\right).
\end{split}
\ee
The estimates \eqref{m-2t} and \eqref{lmu-2t} complete the proof for \eqref{aestlmu}. \eqref{aestmu} will follow from \eqref{aestlmu} as follows.
\bee
\bes
\mu(t,u,\ta)&=\mu(0,u,\ta)+\int_{0}^t\dfrac{(1+\tau)A(\tau)}{(1+\tau)A(\tau)}L\mu(\tau,u,\ta)d\tau\\
&=1+\int_{0}^t\dfrac{1}{(1+\tau)A(\tau)}L\mu(0,u,\ta)+\dfrac{O(\da^{\frac{3}{2}}M)(1-
\dfrac{1}{1+\tau})}{(1+\tau)A(\tau)}d\tau\\
&=1+L\mu(0,u,\ta)\int_0^t\dfrac{1}{(1+\tau)A(\tau)}d\tau+O(\da^{\frac{3}{2}}M)\left(1+\ln(1+t)-\dfrac{1}{1+t}
\right).
\end{split}
\ee
\end{pf}
\begin{remark}
It holds that for the function $A(t)$
\bee
\dfrac{1}{C(a,\lam)}=e^{-|\frac{a}{1-\lam}|}\leq |A(t)|\leq e^{|\frac{a}{1-\lam}|}=C(a,\lam).
\ee
For $a>0$, $A(t)$ is increasing in $t$; for $a<0$, $A(t)$ is decreasing in $t$.
\end{remark}
\subsubsection{\textbf{The behavior of $\mu$ near shock and the damping effect on shock formation}}
\hs Define the shock region as: $W_{shock}=\{(t,u,\ta)|\mu(t,u,\ta)\leq\dfrac{1}{10}\}$. 
\begin{prop}\label{keymu1}
For sufficiently small $\da$ and all $(t,u,\ta)\in W_{shock}$, it holds that
\bee
L\mu(t,u,\ta)\leq-C\dfrac{1}{(1+t)A(t)}\left(\int_0^t\dfrac{1}{(1+\tau)A(\tau)}\right)^{-1},
\ee
for some positive constant $C<1$. Actually, one can refine the estimate to be
\bee
L\mu(t,u,\ta)\leq-C(a,\lam)(1+t)^{-1}(1+\ln(1+t))^{-1},
\ee
where the constant $C(a,\lam)<1$ depends on $a,\lam$.
\end{prop}
\begin{remark}
It can be show that\footnote{One can apply the same argument in \cite{CZD1} Remark3.2.} \textbf{when a shock forms, $\pa_i\pa_j\fe$ blows up, for some $i,j\in\{1,2,3\}$}. That is, for the compressible Euler equations, when a shock forms, $\nabla v$ and $\nabla\p$ blow up.
\end{remark}
\begin{remark}
Here we will illustrate that how the largeness condition leads to the shock formation in finite time and how the damping term affects the shock formation. Since $m=\dfrac{1}{2}\dfrac{dH}{dh}Th+a\dfrac{T\fe}{(1+t)^{\lam}}=(-1-\eta\eta')(T\fai_0-\fai_iT\fai_i+aT\fe)+a\dfrac{T\fe}{(1+t)^{\lam}}$, and
by the estimates \eqref{estdlfai}, \eqref{esttfe}, and $|\eta-\eta_0|+|\eta'-\eta'_0|\les\dfrac{\da^{\frac{3}{2}} M}{1+t}$, it is clear that up to the order $O(\dfrac{\da^{\frac{3}{2}}M}{1+t})$, $T\fai_{\a}(t,u,\ta)$, $\fai_{\a}(t,u,\ta)$, and $\dfrac{T\fe(t,u,\ta)}{(1+t)^{\lam}}$ can be replaced by $\dfrac{T\fai_{\a}(0,u,\ta)}{(1+t)A(t)}$, $\dfrac{\fai_{\a}(0,u,\ta)}{(1+t)A(t)}$ and $\dfrac{T\fe(0,u,\ta)}{(1+t)^{\lam+1}A(t)}$ respectively, so that
\bee
\bes
m&=-2\left(\dfrac{T\fai_0(0,u,\ta)}{(1+t)A(t)}-\dfrac{\fai_i(0,u,\ta)}{(1+t)A(t)}
\dfrac{\fai_i(0,u,\ta)}{(1+t)A(t)}+\dfrac{aT\fe(0,u,\ta)}{(1+t)^{\lam+1}A(t)}\right)+O
\left(\dfrac{\da^{\frac{3}{2}}}{1+t}\right)\\
&=\dfrac{2\da^{\frac{1}{2}}}{(1+t)A(t)}\pa_r\fe_1\left(\dfrac{1-r}{\da},\ta\right)+O
\left(\dfrac{\da^{\frac{3}{2}}}{1+t}\right).
\end{split}
\ee
Hence, $L\mu=m+\mu e=\dfrac{2\da^{\frac{1}{2}}}{(1+t)A(t)}\pa_r\fe_1\left(\dfrac{1-r}{\da},\ta\right)+O
\left(\dfrac{\da^{\frac{3}{2}}}{1+t}\right)$, which can be integrated to obtain
\bee\label{mushock}
\bes
\mu(t,u,\ta)&=\mu(0,\mu,\ta)+\int_{0}^tL\mu(\tau,u,\ta)d\tau\\
&=\mu(0,u,\ta)+\int_0^t\dfrac{2\da^{\frac{1}{2}}}{(1+\tau)A(\tau)}\pa_r\fe_1\left(\dfrac{1-r}{\da},\ta\right)+O
\left(\dfrac{\da^{\frac{3}{2}}}{1+\tau}\right)d\tau\\
&=1+2\da^{\frac{1}{2}}\pa_r\fe_1\left(\dfrac{1-r}{\da},\ta\right)\frac{1}{A(t_0)}\ln(1+t)
+O\left(\da^{\frac{3}{2}}\ln(1+t)\right),
\end{split}
\ee
where $t_0\in[0,t]$ and $A(t)$ is a constants depending on $a,\lam$ with $A(t_0)>(<)1$ for $a>(<)0$.
\been
\item If $\min\pa_s\fe_1\leq -1$, then
\bee
\mu(t,u,\ta)=1-\dfrac{2\da^{\frac{1}{2}}}{A(t_0)}\left|\pa_r\fe_1\right|\ln(1+t)
+O\left(\da^{\frac{3}{2}}\ln(1+t)\right).
\ee
It can be seen that $t$ can not be larger than $T_{\ast}=\exp{\left(\frac{
\da^{-\frac{1}{2}}A(t_0)}{2|\pa_r\fe_1|}\right)}-1$ for sufficiently small $\da $, that is, a shock forms before $t=T_{\ast}$. Furthermore, $T_{\ast}(a,\lam)$ is increasing in $a$ for fixed $\lam$ and decreasing in $\lam$ for fixed $a$, that is,
\begin{itemize}
\item For fixed $\lam$,
\been[(1)]
\item if $a>0$, then the damping effect will delay the shock formation;
\item if $a<0$, then the anti-damping term will lead to the formation of shock in advance.
\een
\item For fixed $a>0(<0)$, as $\lam $ increasing, the (anti-)damping effect becomes weaker so that $A(t)\to 1$ and the solution behaves like the solution of standard Euler system while as $\lam$ decreasing, the (anti-)damping effect becomes stronger and then the shock formation is delayed(advanced);
\end{itemize}
\item If $\min\pa_s\fe_1>0$, then it follows from \eqref{mushock} that $\mu>0$ on $t\in[0,+\infty)$, so that one can obtain a smooth solution on $[0,+\infty)$ to the system\eqref{cee}.
\een
Moreover, if $\dfrac{dH}{dh}=0$, i.e. the flow is the Chaplygin gas, then $\mu>0$ on $t\in[0,+\infty)$, which implies the system \eqref{cee} admits a global solution.
\end{remark}
\begin{pf}
By proposition \ref{accruatemu1}, it holds that
\bee
\int_0^t\dfrac{1}{(1+\tau)A(\tau)}d\tau L\mu(0,u,\ta)\leq\underbrace{\mu(t,u,\ta)-1-O\left(\da^{\frac{3}{2}}M(1+\ln(1+t)
-\dfrac{1}{1+t})\right)}_{I},
\ee
where the term $I$ in the shock region is $\leq-\frac{1}{2}$ for sufficiently small $\da$. It follows from this and Proposition\ref{accruatemu1} that
\bee
(1+t)A(t)L\mu(t,u,\ta)\leq L\mu(0,u,\ta)+O(\da^{\frac{3}{2}}M)\leq-C\left(\int_0^t\dfrac{1}{(1+\tau)A(\tau)}\right)^{-1}.
\ee
\end{pf}
\hs Define $t_1:=\inf\{t|\mu(\tau)\leq\dfrac{1}{10}\text{for all} \tau\geq t\}$. Then, the estimate \eqref{aestlmu} shows that for $t\in[t_1,s]$, $\mu(t,u,\ta)$ behaves like $\ln(1+s)-\ln(1+t)$, which implies that $\mu^{-1}T\mu$ behaves like $\dfrac{\da^{-\frac{1}{2}}M(1+\ln(1+t))}{\ln(1+s)-\ln(1+t)}$. Then, the following integral which will encounter in the energy estimate may not be bounded:
\begin{align*}
&\int_0^s\mu^{-1}T\mu\dfrac{(1+\ln(1+t))^4}{(1+t)^2}dt=\da^{-\frac{1}{2}}M\int_{0}^{\ln(1+s)}\dfrac{(1+x)^5}{(\ln(1+s)-x)e^{x}}dx\\
&\geq\int_{\frac{1}{2}\sigma}^{\sigma}\dfrac{(1+x)^5}
{(\sigma-x)e^{x}}dx\geq\dfrac{(1+\frac{1}{2}\sigma)^5}{e^{\sigma}}\int_{\frac{1}{2}\sigma}^{\sigma}\dfrac{1}
{\sigma-x}dx=\infty.
\end{align*}
However, by changing $\mu^{-1}T\mu$ into its positive part, then the above integral is integrable in time as shown in the following proposition.
\begin{prop}\label{keymu2}
For sufficiently small $\da$ and all $(t,u,\ta)\in W_{shock}$ where $t\in[t_1,s]$, it holds that
\bee\label{mu-1tmu+}
(\mu^{-1}T\mu)_+\les\dfrac{\da^{-\frac{3}{4}}M(1+\ln(1+t))}
{\sqrt{\ln(1+s)-\ln(1+t)}}.
\ee
Moreover, it holds that $\int_0^s\mu^{-1}T\mu\dfrac{(1+\ln(1+t))^4}{(1+t)^2}dt\leq C$ for some constant $C$ which is independent of $s$.
\end{prop}
\begin{pf}
Note that\footnote{See the proof of Proposition3.2 in\cite{CZD1}.}
\[
 |(\mu^{-1}T\mu)_+|\leq\sqrt{\dfrac{||T^2\mu||_{L^{\infty}([0,\tilde{\da}])}}{\inf_{[0,\tilde{\da}]}\mu}}.
\]
It remains to bound $T^2\mu$ and $\mu.$ Following the same argument in the proof of 
Proposition\ref{accruatemu1}, one can obtain
\bee\label{T^3fai}
\left|(1+t)A(t)T\fai_{\a}(t,u,\ta)-T\fai_{\a}(0,u,\ta)\right|\les\da^{\frac{3}{2}} M(1-\dfrac{1}{1+t}).
\ee
The similar argument for \eqref{T^2fe} implies 
\bee\label{T^3fe}
|T^3\fe|\les\da^{-\frac{1}{2}}M.
\ee
 Then, commuting $T^2$ with $L\mu=m+\mu e$ yields
\bee\label{lt2muequation}
LT^2\mu=T^2m+(T^2\mu)e+2(T\mu)(Te)+\mu(T^2e)+[L,T^2]\mu.
\ee
Then, it follows from \eqref{T^3fai}-\eqref{lt2muequation} and the similar argument in proving Proposition\ref{accruatemu1} that, 
\bee\label{LT^2mu}
|(1+t)A(t)LT^2\mu(t,u,\ta)-LT^2\mu(0,u,\ta)|=O(\da^{-\frac{1}{2}}M^2),
\ee
which implies
\bee
|T^2\mu|\les\da^{-\frac{3}{2}}M(1+\ln(1+t)).
\ee
 \hs It follows from  Proposition\ref{keymu1} and \ref{accruatemu1} that for all $(t,u,\ta)\in W_{shock}$
 \begin{equation*}
 L\mu(0,u,\ta)\leq(1+t)A(t)L\mu(t,u,\ta)+O(\da^{\frac{3}{2}}M(1-\dfrac{1}{1+t}))\les
-\left(\int_0^t\dfrac{1}{(1+\tau)A(\tau)}d\tau\right)^{-1}\les-\dfrac{1}{ 1+\ln(1+t)}.
\end{equation*}
Hence,
\begin{align*}
\mu(t,u,\ta)&=\mu(s,u,\ta)+\int_s^t\dfrac{(1+\tau)A(\tau)}{(1+\tau)A(\tau)}L\mu(\tau,u,\ta)d\tau\\
&=\mu(s,u,\ta)+\int_s^t\dfrac{L\mu(0,u,\ta)}{(1+\tau)A(\tau)}+\dfrac
{O(\da^{\frac{3}{2}}M(1-\frac{1}{1+\tau}))}{(1+\tau)A(\tau)}d\tau\\
&=\mu(s,u,\ta)+L\mu(0,u,\ta)\int_0^t\dfrac{1}{(1+\tau)A(\tau)}d\tau+
O(\da^{\frac{3}{2}}M(1-\frac{1}{1+\tau}))\\
&\gtrsim\dfrac{\ln(1+s)-\ln(1+t)}{1+\ln(1+t)},
\end{align*}
for sufficiently small $\da$. Therefore, $(\mu^{-1}T\mu)_+\les\dfrac{\da^{-\frac{3}{4}}M(1+\ln(1+t))}
{\sqrt{\ln(1+s)-\ln(1+t)}}.$\\
\hs As a consequence, it holds that
\begin{align*}
&\int_0^s\mu^{-1}T\mu\dfrac{(1+\ln(1+t))^4}{(1+t)^2}dt\les\da^{-\frac{3}{4}}M\int_{0}^{\ln(1+s)}\dfrac{(1+x)^5}{\sqrt{\ln(1+s)-x}e^{x}}dx\\
&=\da^{-\frac{3}{4}}M\left(\int_{\frac{1}{2}\sigma}^{\sigma}+\int_0^{\frac{1}{2}\sigma}\right)\dfrac{(1+x)^5}
{\sqrt{\sigma-x}e^{x}}dx:=
\da^{-\frac{3}{4}}M(I_1+I_2).
\end{align*}
Note that $|I_2|\leq\int_0^{\frac{1}{2}\sigma}\dfrac{(1+x)^5}{e^{x}\sqrt{x}}dx\leq \int_0^{\infty}\dfrac{(1+x)^6}{e^{2x}\sqrt{x}}dx\leq C_1$, where $C_1$ is independent of $s$. For $I_1$, it holds that
\begin{align*}
I_1&\leq\dfrac{(1+\sigma)^5}{e^{\frac{1}{2}\sigma}}\int_{\frac{1}{2}\sigma}^{\sigma}\dfrac{1}
{\sqrt{\sigma-x}}dx=\dfrac{(1+\sigma)^5}{e^{\frac{1}{2}\sigma}}\cdot\sqrt{2\sigma}\leq C_2,
\end{align*}
where $C_2$ is independent of $s$.
\end{pf}
\section{\textbf{Estimates for deformation tensors}}\label{section4}
\hs The deformation tensor associated with a vector filed $Z$ with respect to $g$ is defined as: $\zpi_{\a\be}=D_{\a}Z_{\be}+D_{\be}Z_{\a}$ or equivalently $\zpi_{\a\be}=\mathcal{L}_{Z}g_{\a\be}$ and with respect to the conformal metric $\til{g}$ is given by
\bee\label{relationdeformationtensor}
\zgpi_{\a\be}=\mathcal{L}_{Z}\til{g}_{\a\be}=\Omega\zpi_{\a\be}+Z(\Omega)g_{\a\be}.
\ee
One can obtain the deformation tensors associated with the commutation vectorfields $\{Z_i\}$ with $i=1,2,3,4,5$ as follows\footnote{The computations are the same as section 4 in\cite{CZD1}.}.\\
\hs For $Z_1=Q=(1+t)L$, the deformation tensor $\qpi$ can be computed in the null frame $(L,\dl,X_1,X_2)$ as:
\bee\label{computationqpi}
\bes
&\qpi_{LL}=0,\hs \qpi_{\dl\dl}=4(1+t)\mu L(\eta^{-2}\mu)-4\mu\eta^{-2}\mu,\\
&\qpi_{L\dl}=-2(1+t)L\mu-2\mu,\hs \qpi_{LA}=0,\\
&\qpi_{\dl A}=2(1+t)(\eta_A+\zeta_{A}),\hs \sqpi_{AB}=2(1+t)\chi_{AB}.
\end{split}
\ee
Hence, the corresponding deformation tensor with respect to $\til{g}$ is given by:
\beeq
&&\qgpi_{LL}=0,\hs \qgpi_{\dl\dl}=4\Omega\mu((1+t)L(\eta^{-2}\mu)-\eta^{-2}\mu),\\
&&\qgpi_{L\dl}=-2(\Omega(1+t)L\mu+\Omega\mu+\mu(1+t)L\Omega),\hs
\qgpi_{LA}=0,\\
&&\qgpi_{\dl A}=2\Omega(1+t)(\eta_A+\zeta_A),\hs \pre {(Q)} {} {\hat{\til{\s{\pi}}}}_{AB}=
2\Omega(1+t)\hat{\chi}_{AB},\hs tr\sqgpi=2\Omega(1+t)\til{tr}\chi.
\eeq
\hs For $Z_2=T$, the deformation tensor $\tpi$ with respect to $g$ in the null frame $(L,\dl,X_1,X_2)$ is given by:
\bee\label{computationtpi}
\bes
&\tpi_{LL}=0,\hs \tpi_{\dl\dl}=2\mu T(\eta^{-2}\mu),\\
&\tpi_{L\dl}=-2 T\mu,\hs \tpi_{LA}=-(\eta_A+\zeta_A),\\
&\tpi_{\dl A}=-\eta^{-2}\mu(\eta_A+\zeta_A),\hs \stpi_{AB}=2\eta^{-1}\mu\ta_{AB},
\end{split}
\ee
and with respect to $\til{g}$ is given by:
\beeq
&&\tgpi_{LL}=0,\hs \tgpi_{\dl\dl}=2\Omega\mu T(\eta^{-2}\mu), \hs \tgpi_{L\dl}=-2(\Omega T\mu+\mu T\Omega),\\
&&\tgpi_{LA}=-\Omega(\eta_A+\zeta_A),\hs \tgpi_{\dl A}=-\Omega\eta^{-2}\mu(\eta_A+\zeta_A),\\
&& \pre {(T)} {} {\hat{\til{\s{\pi}}}}_{AB}=
2\Omega\eta^{-1}\mu\hat{\ta}_{AB},\hs tr\stgpi=2\Omega\eta^{-1}\mu\til{tr}\ta,
\eeq
where $\eta^{-1}\mu\til{tr}\ta=\Omega^{-1}(\Omega\eta^{-1}\mu tr\ta+T\Omega)$.\\
\hs For $Z_{i+2}=R_i$, the deformation tensor $\rpi$ with respect to $g$ in the null frame $(L,\dl,X_1,X_2)$ is given by
\bee\label{computationrpi1}
\bes
\rpi_{LL}&=0,\quad \rpi_{TT}=2\eta^{-1}\mu R_i(\eta^{-1}\mu),\quad \rpi_{LT}=-R_i\mu,\\
\rpi_{LA}&=\lam_i\eta\mu^{-1}\zeta_A+L^{\lam}\ep_{i\lam l}X_A^l-R_i^B\chi_{AB},\\
\rpi_{TA}&=\lam_i X_A(\eta^{-1}\mu)+T^{\lam}\ep_{i\lam l}X_A^l-\eta^{-1}\mu R_i^B\ta_{AB},\quad \rpi_{AB}=-2\lam_i\ta_{AB}.
\end{split}
\ee
\hs We then turn to the estimates for the deformation tensors. Since $|\ep_A|\les|\s{k}|\les\dfrac{\da^{\frac{3}{2}} M^2}{(1+t)^2}$ and
 \bee\label{estzetaetaA}
 |\eta\mu^{-1}\zeta_A|=|\eta\ep_A-\sd_A\eta|\les\dfrac
{\da^{\frac{3}{2}}M}{(1+t)^2},\ |\eta_A|=|\zeta_A+\sd_A\mu|\les \dfrac{\da^{\frac{1}{2}}M(1+\ln(1+t))}{1+t},
 \ee
due to Lemma\ref{ltfe}, the following estimates hold for the deformation tensor $\qgpi$:
\bee\label{estqgpi}
\bes
&\qgpi_{LL}=0,\hs ||\mu^{-1}\qgpi_{\dl\dl}||_{\supnormda}\les \da^{\frac{1}{2}}M(1+\ln(1+t)),\\
&||\qgpi_{L\dl}||_{\supnormda}\les \da^{\frac{1}{2}}M(1+\ln(1+t)),\\
&\qgpi_{LA}=0,\hs ||\mu^{-1}\qgpi_{\dl A}||_{\supnormda}\les \da^{\frac{1}{2}}M(1+\ln(1+t)),\\
&||\pre {(T)} {} {\hat{\til{\s{\pi}}}}_{AB}||_{\supnormda}\les\dfrac{\da^{\frac{3}{2}}M(1+\ln(1+t))}
{1+t},\hs ||tr\sqgpi||_{\supnormda}\les 1.
\end{split}
\ee
Indeed, one can obtain a more precise bound for $tr\sqgpi$ as follows. Note that
\bee
tr\sqgpi-4=2\Omega (1+t)\til{tr}\chi'+\dfrac{4\Omega u}{1-u+t}+4(\Omega-1).
\ee
It follows from Lemma\ref{hrhoeta}and \ref{2ndff} that $||tr\sqgpi-4||_{\supnormda}\les\dfrac{\da^{\frac{3}{2}}M(1+\ln(1+t))}{1+t}$.\\
\hs The following estimates for the deformation tensor $\tgpi$ hold:
\bee\label{esttgpi}
\bes
&\tgpi_{LL}=0,\hs ||\mu^{-1}\tgpi_{\dl\dl}||_{\supnormda}\les\da^{-\frac{1}{2}}M(1+\ln(1+t)),\\
&||\tgpi_{L\dl}||_{\supnormda}\les\da^{-\frac{1}{2}}M(1+\ln(1+t)),\\
&||\tgpi_{LA}||_{\supnormda}\les \dfrac{\da^{\frac{1}{2}}M(1+\ln(1+t))}
{1+t},\hs ||\mu^{-1}\tgpi_{\dl A}||_{\supnormda}\les \dfrac{\da^{\frac{1}{2}}M(1+\ln(1+t))}
{1+t},\\
&||\pre {(T)} {} {\hat{\til{\s{\pi}}}}_{AB}||_{\supnormda}\les\dfrac{\da^{\frac{3}{2}}
M(1+\ln(1+t))^2}{(1+t)^2},\hs ||tr\stgpi||_{\supnormda}\les \dfrac{\da^{\frac{1}{2}}M(1+\ln(1+t))}
{1+t}.
\end{split}
\ee
\hs Define
\[
y^i=\hat{T}^i+\dfrac{x^i}{1-u+t},\quad z^i=L^i-\dfrac{x^i}{1-u+t}.
\]
Then,
\bee\label{ziyi}
z^i=-\eta y^i+\dfrac{(\eta-1)x^i}{1-u+t}-\fai^i.
\ee
Since $x^{\lam}\ep_{i\lam l}X_A^l=\sum_l\mari{R}_i^lX_A^l=\bar{g}(\mari{R}_i,X_A)=\sg_{AB}R_i^B$, one can rewrite $\rpi_{TA}$ and $\rpi_{LA}$ as
\begin{align*}
\rpi_{TA}&=\lam_i X_A(\eta^{-1}\mu)+\eta^{-1}\mu y^{\lam}\ep_{i\lam l}X_A^l-\eta^{-1}\mu R_i^B\ta'_{AB},\\
\rpi_{LA}&=\lam_i\eta\mu^{-1}\zeta_A+z^{\lam}\ep_{i\lam l}X_A^l-R_i^B\chi'_{AB}.
\end{align*}
It suffices to estimate $\lam_i$ and $y^i$. Since on $\Si_t$, $r$ achieves its maximum on $S_{t,0}$, which equals to $1+t$. Hence,
\[
 r\leq (1+t).
 \]
 Since $|Tr|=|\eta^{-1}\mu\sum_i\dfrac{x^i\hat{T}^i}{r}|\les\da^{\frac{1}{2}}M(1+
 \ln(1+t))$, then integrating $Tr$ from $0$ to $u$ and noting that $r=1+t+\int_0^u Tr du'$ yield
 \begin{equation*}
 r\les 1+t+u+\da^{\frac{1}{2}}Mu(1+\ln(1+t)),\hs r\gtrsim 1+t-u-\da^{\frac{1}{2}}Mu(1+\ln(1+t)),
 \end{equation*}
or roughly, $C^{-1}(1+t)\leq r\leq C (1+t)$ for some constant $C$ which will be used without mention.\\
\hs Note that
\[
L\lam_i=\sum_kL(\mari{R}^k)\hat{T}^k+\sum_kL(\hat{T}^k)\mari{R}^k.
\]
and
\bee\label{Llami}
|L\lam_i|=|-\ep_{ijk}\fai^j\hat{T}^k-\eta\mu^{-1}\zeta^Ax^k\ep_{ikl}X_A^l|\les \dfrac{\da^{\frac{3}{2}}M}{1+t},
\ee
due to Lemma\ref{LTidTi} and \eqref{estzetaetaA}. Then,
 \bee\label{estlami}
 ||\lam_i||_{\supnormda}\les\da^{\frac{3}{2}}M(1+\ln(1+t)).
 \ee
\hs Let
\bee\label{defy'i}
y^{' i}=\hat{T}^i+\pa_r^i, \quad \hat{T}=\bar{g}(\pa_r,\hat{T})\pa_r+\Si\hat{T}.
 \ee
where $\Si$ is the projection from $\Si_t$ to the unit Euclidean sphere. Then,
\begin{align}
1=|\hat{T}|^2&=|\Si\hat{T}|^2+|\bar{g}(\pa_r,\hat{T})|^2=r^{-2}\sum_i|\bar{g}(
\mari{R}_i,\hat{T})|^2+|\bar{g}(\pa_r,\hat{T})|^2=r^{-2}\sum_i\lam_i^2+
|\bar{g}(\pa_r,\hat{T})|^2,\label{1-g(parT)^2}\\
|y'|^2&=|\hat{T}|^2+1+2\bar{g}(\pa_r,\hat{T})=r^{-2}\sum_i\lam_i^2+|1+\bar{g}
(\pa_r,\hat{T})|^2.\label{1-g(parT)}
\end{align}
On $S_{t,0}$, $g(\pa_r,\hat{T})=-1$, which implies the angle between $\pa_r$ and $\hat{T}$ is larger than $\dfrac{\pi}{2}$ due to the continuity of $\hat{T}$ in $u$ and the bootstrap assumptions. Then,
\bee
1+g(\pa_r,\hat{T})\les1-|g(\pa_r,\hat{T})|^2\les
\dfrac{\da^3M^2(1+\ln(1+t))^2}{(1+t)^2},
\ee
due to \eqref{1-g(parT)^2} and \eqref{estlami}. Thus, it follows from\eqref{1-g(parT)} that
\bee\label{esty'i}
|y'|\les|r|^{-1}|\lam_i|+|1-g(\pa_r,\hat{T})|\les\dfrac{\da^{\frac{3}{2}}M(1+\ln(1+t))}{1+t}.
\ee
Therefore,
\bee\label{estyi}
|y^{i}|\leq|y^{'i}|+|y^i-y^{' i}|\leq|y^{'i}|+\arrowvert\dfrac{1}{r}-\dfrac{1}{1-u+t}\arrowvert\cdot|x^i|\les
\dfrac{\da^{\frac{1}{2}}M(1+\ln(1+t))}{1+t}.
\ee
Thus, it follows from \eqref{relationdeformationtensor}, \eqref{estlami} and \eqref{estyi} that
\bee\label{estrgpi}
\bes
&\rgpi_{LL}=0,\hs ||\mu^{-1}\rgpi_{\dl\dl}||_{\supnormda}\les \da^{\frac{1}{2}}M(1+\ln(1+t)),\\
&||\rgpi_{L\dl}||_{
\supnormda}\les\da^{\frac{1}{2}}M(1+\ln(1+t)),\\
&||\rgpi_{LA}||_{\supnormda}\les\dfrac{\da^{\frac{1}{2}}M(1+\ln(1+t))}{1+t},\hs ||\rgpi_{\dl A}||_{\supnormda}\les\dfrac{\da M^2(1+\ln(1+t))^2}{1+t},\\
&||\pre {(R_i)} {} {\hat{\til{\s{\pi}}}}_{AB}||_{\supnormda}\les
\dfrac{\da^2M^2(1+\ln(1+t))^2}{(1+t)^2},\hs ||tr\srgpi||_{\supnormda}\les \dfrac{\da^{\frac{1}{2}}M(1+\ln(1+t))}{1+t}.
\end{split}
\ee
\hs We state the following Proposition which shows the relations between some elementary derivatives, whose proof can be found in\cite{christodoulou2014compressible}.
\begin{prop}\label{relationangular}
Under the bootstrap assumptions, there exists a numerical constant $C$ such that for any $S_{t,u}$ $1-$form $\xi$, the following bounds hold:
\bee\label{h0}
|\xi|^2\leq C(1+t)^{-2}\sum_i|\xi(R_i)|^2,
\ee
\bee\label{h1}
|\s{D}\xi|^2+(1+t)^{-2}|\xi|^2\leq C(1+t)^{-2}\sum_i|\lie_{R_i}\xi|^2.
\ee
Similarly, for any $S_{t,u}$ $(0,2)$ tensor $\lam$, \eqref{h0} and \eqref{h1} hold with $\xi$ replaced by $\lam$.
\end{prop}
\begin{remark}
In particular, taking $\xi=\sd f$ for any smooth function $f$ into \eqref{h0} and \eqref{h1} yields
\begin{align*}
|\sd f|^2&\leq C(1+t)^{-2}\sum_i|R_i f|^2,\\
|\s{D}^2f|^2+(1+t)^{-2}|\sd f|^2&\leq  C(1+t)^{-2}\sum_i|\lie_{R_i}\sd f|^2.
\end{align*}
\end{remark}
\section{\textbf{Multiplier and commuting vector fields methods, estimates on lower order terms}}\label{section6}
\hs We derive the energy estimates by following the framework in \cite{christodoulou2014compressible}and \cite{CZD1}. Here, the following two vector fields will be used as multipliers:
     \bee\label{defmultipliers}
      K_{0}=(1+\mu)L+\dl,\quad K_{1}=\dfrac{2(1+t)}{\frac{1}{2}\til{tr}\chi},
      \ee
      which are adapted from\cite{christodoulou2014compressible}. These two multipliers are the analog of $\pa_t$ which commutes with the linear wave equation and the Morawetz vector field $\dfrac{1}{2}[(t-r)^2(\pa_t+\pa_r)+(t+r)^2(\pa_t-\pa_r)]$ in the Minkowski space, respectively.\\
\hs Let $\mu_{\sg}$ be the area form of $S_{t,u}$. Since $|\det g|=\mu^2\det \sg$, which implies that the volume form on $W_{t}^{ u}$ equals to $\mu dt d ud\mu_{\sg}$. Therefore, for convenience, the following notations will be adopted:
\beeq
&&\int_{\Si_{t}^{ u}}f=\int_{0}^{ u}\int_{S_{t, u}}f d\mu_{\sg}d u',\\
&&\int_{C_{ u}^{t}}f=\int_{0}^{t}\int_{S_{t, u}}f d\mu_{\sg}dt',\\
&&\int_{W_{t}^{ u}}f=\int_{0}^{ u}\int_{0}^t\int_{S_{t, u}}\mu\cdot f(t', u',\ta) d\mu_{\sg}d u'dt'.\\
\eeq
\subsection{\textbf{Multiplier method, fundamental energy estimates}}
\hs Let $f$ be the right hand side of\eqref{nonlinearwave} and then \eqref{nonlinearwave} becomes: $\square_{\til{g}}\fai=f$. 
Define the energy-momentum tensor associated with $\fai$ as:
\begin{equation}\label{energymomentum}
\til{T}_{\mu\nu}:=\pa_{\mu}\fai\pa_{\nu}\fai-\dfrac{1}{2}\til{g}_{\mu\nu}\til{g}^{\a\be}\pa_{\a}\fai\pa_{\be}\fai
=\pa_{\mu}\fai\pa_{\nu}\fai-\dfrac{1}{2}g_{\mu\nu}g^{\a\be}\pa_{\a}\fai\pa_{\be}\fai:=T_{\mu\nu},
\end{equation}
which describes the density of energy and momentum in spacetime.
Define
\beeq
\til{P}^{\mu}_0:&=&-\til{T}^{\mu}_{\nu}K_0^{\nu},\\
\til{P}^{\mu}_1:&=&-\til{T}^{\mu}_{\nu}K_1^{\nu}-\til{g}^{\mu\nu}\left(2(1+t)\fai\pa_{\nu}\fai-\fai^2\pa_{\nu}(1+t)\right),
\eeq
where $K_0$ and $K_1$ are the two multipliers defined by \eqref{defmultipliers}. Then,
\bee\label{DmuPmu}
\bes
\til{D}_{\mu}\til{P}^{\mu}_0&=-\left(\p K_0\fai+\dfrac{1}{2}\til{T}^{\mu\nu}\til{\pi}_{0,\mu\nu}\right):=\til{Q}_0,\\
\til{D}_{\mu}\til{P}^{\mu}_1&=-\p(K_1\fai+2(1+t)\fai)-\dfrac{1}{2}\til{T}^{\mu\nu}\til{\pi}'_{1,\mu\nu}-\fai^2\square_{\til{g}}(1+t):=\til{Q}_1,
\end{split}
\ee
where $\til{\pi}_{\a,\mu\nu}$ are the deformation tensors of $K_{\a} (\a=0,1)$ with respect to $\til{g}$ and $\til{\pi}'_{1,\mu\nu}=\Omega\pi_{1,\mu\nu}+(K_{1}\Omega-4\Omega(1+t))g_{\mu\nu}$. Taking $f=(1+t)$ in \eqref{decomposewaveeq} and noticing that $L(1+t)=1$, $\dl(1+t)=\eta^{-2}\mu$ and $\sd(1+t)=0$ yield
\bee\label{eq1+t}
\Omega\square_{\til{g}}(1+t)=-\mu^{-1}L(\eta^{-2}\mu)-\dfrac{1}{2\mu}\left(\eta^{-2}\mu\til{tr}\chi+\til{tr}\ud{\chi}\right).
\ee
Note that
\bee\label{tilDtilPtilQ}
\til{D}_{\mu}\til{P}_{\a}^{\mu}=\dfrac{1}{\sqrt{|\det\til{g}|}}\pa_{\mu}\left(\sqrt{|\det\til{g}|}\til{P}_{\a}^{\mu}\right)
=\dfrac{1}{\Omega\sqrt{|\det g|}}\pa_{\mu}\left(\Omega^2\sqrt{|\det g|}\til{P}_{\a}^{\mu}\right)=\Omega^{-2}D_{\mu}P_{\a}^{\mu},
\ee
where $P_{\a}^{\mu}=\Omega^{2}\til{P}_{\a}^{\mu}$ and $Q_{\a}=\Omega^2\til{Q}_{\a}$ for $\a=0,1$. Then, it follows from\eqref{DmuPmu} that
\bee\label{DPQ}
D_{\mu}P_{\a}^{\mu}=Q_{\a},\ \a=0,1.
\ee
The following divergence theorem holds:
\begin{lem}\label{divergencethm}
For any spacetime vector field $J=J^{t}\dfrac{\pa}{\pa t}+J^{u}\dfrac{\pa}{\pa u}+\s{J}^A\dfrac{\pa}{\pa \ta^A}$, it holds that
\beeq
\int_{W_t^u}\mu D_{\a}J^{\a}dt'du'd\mu_{\sg}&=&\int_{\Si_t^u}\mu J^t du'd\mu_{\sg}-\int_{\Si_0^u}\mu J^t du'd\mu_{\sg}\\
&+&\int_{C_u^t}\mu J^u dt'd\mu_{\sg}-\int_{C_0^t}
\mu J^u dt'd\mu_{\sg}.
\eeq
\end{lem}
\hs Note that
\bee\label{P0tP0u}
\bes
\mu P_0^t&=-\dfrac{1}{2}\left(P_{0,\dl}+\eta^{-2}\mu P_{0,L}\right)=\dfrac{\Omega}{2}\left((1+\mu)T_{\dl L}+T_{\dl\dl}+\eta^{-2}\mu T_{L\dl}+\eta^{-2}\mu(1+\mu)T_{LL}\right)\\
&=\dfrac{\Omega}{2}\left(\eta^{-2}\mu(1+\mu)(L\fai)^2+(1+\mu+\eta^{-2}\mu)\mu\sdfai^2+(\dl\fai)^2\right),\\
\mu P_0^u&=-P_{0,L}=\Omega\left((1+\mu)T_{LL}+T_{L\dl}\right)=\Omega\left((1+\mu)(L\fai)^2+\mu\sdfai^2\right),
\end{split}
\ee
and
\bee\label{P1tP1u}
\bes
\mu P_1^t&=-\dfrac{1}{2}(P_{1,\dl}+\eta^{-2}\mu P_{1,L})\\
&=\dfrac{\Omega}{2}(\dfrac{4(1+t)}{\til{tr}\chi}\mu|\sd\fai|^2+\eta^{-2}\mu\dfrac{4(1+t)}{\til{tr}\chi}(L\fai)^2\\
&+2(1+t)\fai[\dl\fai+\eta^{-2}\mu L\fai]-\fai^2[\dl(1+t)+
\eta^{-2}\mu L(1+t)]),\\
\mu P_1^u&=-P_{1,L}=\Omega\left(\dfrac{4(1+t)}{\til{tr}\chi}(L\fai)^2+2(1+t)\fai L\fai-\fai^2 L(1+t)\right).
\end{split}
\ee
Denote
\beeq
E'_1(t)&=&\int_{\Si_t^u}2\Omega\dfrac{1+t}{\til{tr}\chi}\left(\eta^{-2}\mu\left(L\fai+\dfrac{1}{2}\til{tr}\chi\fai\right)^2+\mu\sdfai^2\right) d u' d\mu_{\sg},\\
F_1'(t)&=&\int_{C_u^t}4\Omega\dfrac{1+t}{\til{tr}\chi}\left(L\fai+\dfrac{1}{2}\til{tr}\chi\fai\right)^2 dt'd\mu_{\sg}.
\eeq
Then,
\beeq
\int_{C_u^t}\mu P_1^{\mu}\ dt'd\mu_{\sg}-F_1'(t)&=&\int_{C_u^t}-\Omega\left(2(1+t)\fai L\fai+(1+t)\til{tr}\chi\fai^2+\fai^2L(1+t)\right)\\
&=&-\int_{C_u^t}\Omega\left(L((1+t)\fai^2)+(1+t)\til{tr}\chi\fai^2\right).
\eeq
Since $\lie_Ld\mu_{\sg}=\lie_L\sqrt{\det\sg}d\ta^1\wedge d\ta^2=tr\chi d\mu_{\sg}$, then for any function $f$
\bee\label{patCut}
\bes
\dfrac{\pa}{\pa t}\left(\int_{C_u^t}\Omega f dt'd\mu_{\sg}\right)&=\int_{C_u^t}\Omega\left(Lf+(tr\chi+L(\log\Omega))f\right) dt'd\mu_{\sg}\\
&=\int_{C_u^t}\Omega\left(Lf+\til{tr}\chi f\right) dt'd\mu_{\sg}.
\end{split}
\ee
Taking $f=(1+t)\fai^2$ in \eqref{patCut} yields
\bee\label{differenceflux}
\int_{C_u^t}\mu P_1^{\mu}\ dt'd\mu_{\sg}-F_1'(t)=-\int_{S_{t,u}}\Omega(1+t)\fai^2 d\mu_{\sg}+\int_{S_{0,u}}\Omega\fai^2 d\mu_{\sg}.
\ee
Similarly,
\bee
\bes
&\int_{\Si_t^u}\mu P_1^t\ du'd\mu_{\sg}-E'_1(t)\\
&=
\int_{\Si_t^u}\dfrac{\Omega}{2}\left(2(1+t)\fai(2T\fai)-\fai^2\left(\eta^{-2}\mu L(1+t)+\dl(1+t)\right)-\eta^{-2}\mu(1+t)\til{tr}\chi\fai^2\right) du'd\mu_{\sg}.
\end{split}
\ee
Since $\lie_T d\mu_{\sg}=\eta^{-1}\mu tr\ta d\mu_{\sg}$, then for any function $f$
\bee\label{pauStu}
\dfrac{\pa}{\pa u}\left(\int_{S_{t,u}}\Omega f d\mu_{\sg}\right)=\int_{S_{t,u}}\Omega\left(Tf+(\eta^{-1}\mu tr\ta+T\log\Omega)f\right)d\mu_{\sg}.
\ee
Taking $f=(1+t)\fai^2$ in \eqref{pauStu} yields
\beeq
&&\int_{\Si_t^u}\Omega\left((1+t)\fai(2T\fai)+\fai^2T(1+t)+(\eta^{-1}\mu tr\ta+T\log\Omega)(1+t)^2\fai^2\right) du'd\mu_{\sg}\\
=&&\int_{S_{t,u}}\Omega(1+t)\fai^2 d\mu_{\sg}-\int_{S_{t,0}}\Omega(1+t)\fai^2 d\mu_{\sg},
\eeq
where the second integral on the last line vanishes. Thus,
\bee\label{differenceenergy}
\int_{\Si_t^u}\mu P_1^t\ du'd\mu_{\sg}-E_1'(t)=\int_{S_{t,u}}\Omega(1+t)\fai^2 d\mu_{\sg}-I^u(t),
\ee
where
\bee
I^u(t)=\int_{\Si_t^u}\Omega\fai^2\left(\dl(1+t)+(1+t)\til{tr}\dc\right) du'd\mu_{\sg}.
\ee
Applying Lemma\ref{divergencethm} to \eqref{DPQ}, collecting \eqref{P0tP0u}, \eqref{P1tP1u}, \eqref{differenceflux} and \eqref{differenceenergy} and noticing that $\fai$ vanishes on $C_0^u$ yield
\bee\label{energy1'}
\bes
\int_{\Si_t^u}\mu P_{0}^t+\int_{C_u^t}\mu P_{0}^u&=\int_{\Si_{0}^u}\mu P_{0}^t+\int_{W_t^u}Q_0d\mu_g,\\
E'_1(t)+F'_1(t)&=\int_{\Si_t^u}\mu P_1^t\ du'd\mu_{\sg}+\int_{C_u^t}\mu P_1^{\mu}\ dt'd\mu_{\sg}+I^u(t)-\int_{S_{0,u}}(1+t)\Omega\fai^2 d\mu_{\sg}\\
&=\int_{\Si_0^u}\mu P_1^t\ du'd\mu_{\sg}+I^u(t)-\int_{S_{0,u}}(1+t)\Omega\fai^2 d\mu_{\sg}+\int_{W_t^u}Q_1 d\mu_{\sg}\\
&=I^u(t)-I^u(0)+E_1'(0)+\int_{W_t^u}Q_1d \mu_{\sg}.
\end{split}
\ee
\hs Since $\Omega\sim 1$ and $\eta\sim 1$, one may define the energies and the fluxes as follows:
\beeq
E_0(t,u)&=&\int_{\Si_t^u}\mu(1+\mu)[(L\fai)^2+\sdfai^2]+(\dl\fai)^2 du'd\mu_{\sg},\\
F_0(t,u)&=&\int_{C_u^t}(1+\mu)(L\fai)^2+\mu\sdfai^2 dt'd\mu_{\sg},\\
E_1(t,u)&=&\int_{\Si_t^u}\dfrac{1+t}{\til{tr}\chi}\left(\mu(L\fai+\dfrac{1}{2}\til{tr}
\chi\fai)^2+\mu\sdfai^2\right)du'd\mu_{\sg},\\
F_1(t,u)&=&\int_{C_u^t}\dfrac{1+t}{\til{tr}\chi}(L\fai+\dfrac{1}{2}\til{tr}\chi\fai
)^2dt'd\mu_{\sg}.
\eeq
Then, \eqref{energy1'} implies
\bee\label{energy1}
\bes
E_{0}(t)+F_{0}(t)&\les E_{0}(0)+\int_{W_t^u}Q_{0},\\
E_{1}(t)+F_{1}(t)&\les I^u(t)-I^u(0)+E_1(0)+\int_{W_t^u}Q_1,
\end{split}
\ee
where $Q_0,Q_1$ are given in\eqref{DmuPmu}. Denote
\bee\label{modifiedef}
\bes
\til{E}_0(t,u)&=\sup_{t'\in[0,t]}E_0(t',u),\\
\til{F}_0(t,u)&=F_0(t,u),\\
\til{E}_1(t,u)&=\sup_{t'\in[0,t]}(1+\ln(1+t'))^{-4}E_1(t',u),\\
\til{F}_1(t,u)&=\sup_{t'\in[0,t]}(1+\ln(1+t'))^{-4}F_1(t',u),
\end{split}
\ee
which will be used later. Note that $\til{\pi}_{0,\mu\nu}=\Omega\pi_{0,\mu\nu}+(K_0\Omega)g_{\mu\nu}$ and $\til{\pi}'_{1,\mu\nu}=\Omega\pi_{1,\mu\nu}+(K_{1}\Omega-4\Omega(1+t))g_{\mu\nu}$ can be computed as follows.
\beeq
&&\til{\pi}_{0,LL}=0,\hs \til{\pi}_{0,L\dl}=-2\Omega\left((1+2\mu)L\mu+\dl\mu+\mu L(\eta^{-2}\mu)\right)-2\mu((1+\mu)L\Omega+\dl\Omega),\\
&&\til{\pi}_{0,\dl\dl}=4\Omega\mu\left((1+\mu)L(\eta^{-2}\mu)-\dl\mu\right),\hs \til{\pi}_{0,LA}=-2\Omega(\zeta_A+\eta_A),\\
&&\til{\pi}_{0,\dl A}=2\Omega\left((1+\mu)(\eta_A+\zeta_A)-\mu(\sd_A\mu+\sd_A(\eta^{-2}\mu))\right),\\
&&\hat{\til{\s{\pi}}}_{0,AB}=2\Omega[(1+\mu)\hat{\chi}_{AB}+\hat{\dc}_{AB}],\hs tr\til{\s{\pi}}_0=2\Omega[(1+\mu)\til{tr}\chi+\til{tr}\dc].
\eeq
\beeq
&&\til{\pi}'_{1,LL}=0,\hs \til{\pi}'_{1,\dl\dl}=4\Omega\mu\left(\dfrac{4(1+t)}{\til{tr}\chi}L(\eta^{-2}\mu)
-\dl\left(\dfrac{4(1+t)}{\til{tr}\chi}\right)\right),\\
&&\til{\pi}'_{1,L\dl}=8\Omega\left(\dfrac{1+t}{\til{tr}\chi}(\mu L(\eta^{-2}\mu)-L\mu)-\mu L\left(\dfrac{1+t}{\til{tr}\chi}\right)\right)-8\mu\left[
\dfrac{1+t}{\til{tr}\chi}L\Omega-\Omega(1+t)\right],\\
&&\til{\pi}'_{1,LA}=0,\hs \til{\pi}'_{1,\dl A}=8\Omega\left(\dfrac{1+t}{\til{tr}\chi}(\eta_A+\zeta_A)-\mu X_A\left(\dfrac{1+t}{\til{tr}\chi}\right)\right),\\
&&\hat{\til{\s{\pi}}}'_{1,AB}=8\Omega\dfrac{1+t}{\til{tr}\chi}\hat{\chi}_{AB},\hs tr\hat{\til{\s{\pi}}}'_{1,AB}=0.
\eeq
\hs Then, one can rewrite $Q_0=-\Omega^2 f K_0\fai-\dfrac{1}{2}T^{uv}\til{\pi}_{0,uv}=\sum_{\a=0}^7Q_{0,\a},$ where
\begin{align*}
Q_{0,0}&=-\Omega^2 f K_0\fai\\
&=\Omega^2\mu^{-1}K_0\fai\left(-\dfrac
{2\mu\eta'}{\eta^2}A(\dfrac{a}{(1+t)^{\lam}}\fe)\de\fe+\dfrac{\mu}{\p\eta}A\left[\dfrac{a}{(1+t)^{\lam}}\left(\dfrac{\pa\fe}{\pa t}-\frac{1}{2}|\nabla\fe|^2-\frac{\lam}{(1+t)}\fe\right)\right]\right),\\
Q_{0,1}&=-\dfrac{1}{2}T^{LL}\til{\pi}_{0,LL}=-\dfrac{1}{8}\mu^{-2}(\dl\fai)^2\til{\pi}_{0,LL}=0,\\
Q_{0,2}&=-\dfrac{1}{2}T^{\dl\dl}\til{\pi}_{0,\dl\dl}=
-\dfrac{1}{8}\mu^{-2}(L\fai)^2\til{\pi}_{0,\dl\dl}=-\dfrac{1}{2}\mu^{-1}\Omega((1+\mu)L(\eta^{-2}\mu)-\dl\mu)(L\fai)^2,\\
Q_{0,3}&=-T^{L\dl}\til{\pi}_{0,L\dl}=
-\dfrac{1}{4}\mu^{-1}\sdfai^2\til{\pi}_{0,L\dl}\\
&=\dfrac{\Omega}{2}\sdfai^2\left(\boxed{\mu^{-1}L\mu+\mu^{-1}\dl\mu}+2L\mu
+L(\eta^{-2}\mu)+\Omega^{-1}K_0\Omega\right),\\
Q_{0,4}&=-T^{LA}\til{\pi}_{0,LA}=
\dfrac{1}{2}\mu^{-1}(\dl\fai)(\sd^A\fai)\til{\pi}_{0,LA}=-\Omega\mu^{-1}(\dl\fai)(\sd^A\fai)(\zeta_A+\eta_A),\\
Q_{0,5}&=-T^{\dl A}\til{\pi}_{0,\dl A}=\dfrac{1}{2}\mu^{-1}(L\fai)(\sd^A\fai)\til{\pi}_{0,\dl A}\\
&=\Omega\mu^{-1}(L\fai)(\sd^A\fai)((1+\mu)(\eta_A+\zeta_A)-\mu X_A\mu+
\mu X_A(\eta^{-2}\mu)).
\end{align*}
 And, $-\dfrac{1}{2}T^{AB}\til{\pi}_{0,AB}=Q_{0,6}+Q_{0,7}$, where
\beeq
Q_{0,6}&=&-\Omega[\hat{\dc}_{AB}+(1+\mu)\hat{\chi}_{AB}](\sd^A\fai)(\sd^B\fai),\\
Q_{0,7}&=&-\dfrac{1}{4\mu}tr\til{\s{\pi}}_0(L\fai)(\dl\fai)
=-\dfrac{1}{2\mu}\Omega(\til{tr}\dc+(1+\mu)\til{tr}\chi)(L\fai)(\dl\fai).
\eeq
Similar computations lead to the decomposition of $Q_1=\sum_{\a=0}^8Q_{1,\a}$, where
\begin{equation*}
\bes
Q_{1,0}&=-\Omega^2\p (K_1\fai+2(1+t)\fai)\\
&=\Omega^2\mu^{-1}(K_1\fai+2(1+t)\fai)\cdot\\
&\left(-\dfrac{2\mu\eta'}{\eta^2}A(\dfrac{a}{(1+t)^{\lam}}\fe)\de\fe+\dfrac{\mu}{\p\eta}A\left[\dfrac{a}{(1+t)^{\lam}}\left(\dfrac{\pa\fe}{\pa t}-\frac{1}{2}|\nabla\fe|^2-\frac{\lam}{(1+t)}\fe\right)\right]\right),\\
Q_{1,1}&=-\dfrac{1}{2}T^{LL}\til{\pi}_{1,LL}=-\dfrac{1}{8}\mu^{-2}(\dl\fai)^2\til{\pi}_{1,LL}=0,\\
Q_{1,2}&=-\dfrac{1}{2}T^{\dl\dl}\til{\pi}_{1,\dl\dl}=
-\dfrac{1}{8}\mu^{-2}(L\fai)^2\til{\pi}_{1,\dl\dl}=-\dfrac{1}{2}\mu^{-1}\Omega\left(\dfrac{4(1+t)}{\til{tr}\chi}L(\eta^{-2}\mu)
-\dl\dfrac{4(1+t)}{\til{tr}\chi}\right)(L\fai)^2,\\
Q_{1,3}&=-T^{L\dl}\til{\pi}_{1,L\dl}=
-\dfrac{1}{4}\mu^{-1}\sdfai^2\til{\pi}_{1,L\dl}\\
&=-2\Omega\sdfai^2\left(\dfrac{1+t}{\til{tr}\chi}[L(\eta^{-2}\mu)-
\boxed{\mu^{-1}L\mu}]-L\dfrac{1+t}{\til{tr}\chi}-\dfrac{1+t}{\til{tr}\chi}
\Omega^{-1}L\Omega+(1+t)\right),\\
Q_{1,4}&=-T^{LA}\til{\pi}_{1,LA}=
\dfrac{1}{2}\mu^{-1}(\dl\fai)(\sd^A\fai)\til{\pi}_{1,LA}=0,\\
Q_{1,5}&=-T^{\dl A}\til{\pi}_{1,\dl A}=\dfrac{1}{2}\mu^{-1}(L\fai)(\sd^A\fai)\til{\pi}_{1,\dl A}=4\Omega\mu^{-1}(L\fai)(\sd^A\fai)(\dfrac{1+t}{\til{tr}\chi}(\eta_A+\zeta_A)
-\mu X_A\dfrac{1+t}{\til{tr}\chi}),\\
Q_{1,6}&=-\dfrac{1}{2}\left[(\sd^A\fai)(\sd^B\fai)-\dfrac{1}{2}\sg^{AB}\sdfai^2
\right]\cdot\hat{\til{\s{\pi}}}_{1,AB}=-4\Omega\dfrac{1+t}{\til{tr}\chi}\hat{\chi}_{AB}\sd^A\fai\sd^B\fai,\\
Q_{1,7}&=-\dfrac{1}{4\mu}tr\til{\s{\pi}}_1(L\fai)(\dl\fai)
=0\footnotemark,\\
Q_{1,8}&=\Omega^2\fai^2\square_{\til{g}}(1+t).
\end{split}
\end{equation*}
\footnotetext{This is an important fact and is also the reason of choosing $K_1$ as a multiplier. Otherwise, the growth of $t$ will prevent us to bound $\dl\fai$.}
Due to the definitions of $E$ and $F$, it turns out that the boxed terms are most difficult to handle since they involve $\mu^{-1}\sdfai^2$ and possibly shock formation ($\mu\to0$).
\subsubsection{\textbf{Estimate for the error integrals $Q_i$}}
\begin{lem}\footnote{The proof is the same as Lemma5.5 in\cite{CZD1}.}\label{faibound}
There is numerical constant $C$ such that
\bee
\int_{S_{t,u}}\fai^2 d\mu_{\sg}\leq C\tilde{\da}\int_0^u\int_{S_{t,u}}\eta^{-4}\mu^2(L\fai)^2+(\dl\fai)^2\leq C\tilde{\da} E_0(t,u).
\ee
\end{lem}
\begin{remark}
As a corollary,
\[
\int_{W_t^u}\mu^{-1}(1+t)(1+\ln(1+t))^3(\til{tr}\chi\fai)^2\les\da(1+\ln(1+t))^4\til{E}_0(t,u),
\]
which will be used without mention.
\end{remark}
\begin{remark}\label{feL2}
Similar argument leads to
\bee
\int_{W_t^u}\mu^{-1}\dfrac{(1+\ln(1+t))\fe^2}{1+t}\les\til{\da}(1+\ln(1+t))^4\til{E}_0(t,u).
\ee
\end{remark}
\hs We start with the most difficult terms. First,
\beeq
\int_{W_t^u}Q_{0,3}&=&\boxed{\frac{\Omega}{2}\int_{W_t^u}
(\mu^{-1}L\mu+\mu^{-1}T\mu)\sdfai^2}+\frac{\Omega}{2}\int_{W_t^u}(2L\mu+L(\eta^{-2}\mu)+
\Omega^{-1}K_0\Omega)\sdfai^2\\
&\les&\int_{W_t^u}
(\mu^{-1}L\mu+\mu^{-1}T\mu)\sdfai^2+\da^{\frac{1}{2}}\int_0^uF_0(u')du'.
\eeq
\textbf{To deal with any integral involving $\mu^{-1}L\mu$ or $\mu^{-1}\dl\mu$, one can split it into shock part (in $W_{shock}$) and non-shock part. The difficult shock part can be bounded by using two key propositions of $\mu$: Proposition\ref{keymu1} and \ref{keymu2} as follows.}
\[
\int_{W_t^u}
(\mu^{-1}L\mu+\mu^{-1}T\mu)\sdfai^2=\left(\int_{W_t^u\cap W_{shock}}
+\int_{W_t^u\cap W_{n-s}}\right)(\mu^{-1}L\mu+\mu^{-1}T\mu)\sdfai^2,
\]
where $W_{n-s}=W_t^u\setminus W_{shock}$.\\
\hs Since $\mu$ has a positive lower bound on $W_t^u\cap W_{n-s}$, then
\begin{equation*}
\int_{W_t^u\cap W_{n-s}}(\mu^{-1}L\mu+\mu^{-1}T\mu)\sdfai^2\les
\int_{W_t^u\cap W_{n-s}}(L\mu+T\mu)\sdfai^2\les\da^{-\frac{1}{2}}\int_0^t\dfrac{1+\ln(1+t')}{(1+t')^2}E_1(t')dt'.
\end{equation*}
In the shock region, it follows from Proposition \ref{keymu1}, \ref{keymu2} and $\sdfai^2\ge 0$ that
\bee\label{q031}
\int_{W_t^u\cap W_{shock}}(\mu T\mu)\sdfai^2\les
\int_0^t\dfrac{\da^{-\frac{3}{4}}(1+\ln(1+t'))}{\sqrt{\ln(1+s)-\ln(1+t')}}
\dfrac{E_1(t)}{(1+t')^2}dt',
\ee
\bee\label{q032}
\int_{W_t^u\cap W_{shock}}(\mu L\mu)\sdfai^2\les
\int_{W_t^u\cap W_{shock}}-\mu^{-1}\sdfai^2\dfrac{1}{(1+t')\ln(1+t')}
\les-\dfrac{1}{(1+t)^2}K(t,u),
\ee
where $K(t,u)$ is the non-negative integral
\[
K(t,u)=\int_{W_t^u\cap W_{shock}}\mu^{-1}\sdfai^2\dfrac{1+t'}{1+\ln(1+t')},
\]
which is the key role to control $\mu^{-1}\sdfai^2$ in $W_{shock}$. For $Q_{1,3}$, it holds that
\beeq
\int_{W_t^u}Q_{1,3}&=&\boxed{\int_{W_t^u}(\mu^{-1} L\mu)(1+t')^2\sdfai^2}\\
&+&\int_{W_t^u}\dfrac{1+t'}{\til{tr}\chi}\left[
L(\eta^{-2}\mu-\Omega^{-1}L\Omega)+(\dfrac{1}{2}\til{tr}\chi-\dfrac{1}{1+t'})
+(\dfrac{1}{2}\til{tr}\chi+\dfrac{L\til{tr}\chi}{\til{tr}\chi})\right]\sdfai^2\\
&=&I_1+I_2.
\eeq
Since
\beeq
\left|\dfrac{1}{2}\til{tr}\chi-\dfrac{1}{1+t}\right|&=&
\left|\dfrac{u}{(1+t)(1-u+t)}+\dfrac{1}{2}\til{tr}\chi'+\dfrac{1}{2}L(\ln\Omega)
\right|\les\dfrac{\da^{\frac{1}{2}}(1+\ln(1+t))}{(1+t)^2},\\
\left|\dfrac{1}{\til{tr}\chi}(L\til{tr}\chi+\dfrac{1}{2}(\til{tr}\chi)^2)\right|
&=&\left|\dfrac{1}{\til{tr}\chi}(L(tr\chi'+\dfrac{2}{1-u+t}+L(\ln\Omega))+\dfrac{1}
{2}(tr\chi'+\dfrac{2}{1-u+t}+L(\ln\Omega))^2)\right|\\
&\les&\dfrac{1}{(1+t)(1+\ln(1+t))},
\eeq
due to Lemma\ref{2ndff}, then
\bee\label{q131}
I_2\les\int_0^t\dfrac{1}{(1+t')(1+\ln(1+t'))}E_1(t')dt'\leq\dfrac{1}{4}(1+\ln(1+t))^4\til{E}_{1}(t,u).
\ee
\hs As before,
\[
\int_{W_t^u}(\mu^{-1} L\mu)(1+t')^2\sdfai^2=\left(
\int_{W_t^u\cap W_{shock}}+\int_{W_t^u\cap W_{n-s}}\right)(\mu^{-1} L\mu)(1+t')^2\sdfai^2.
\]
\hs In the non-shock region,
\bee
\int_{W_t^u\cap W_{n-s}}(\mu^{-1}L\mu)(1+t')^2\sdfai^2\les\int_0^t\da^{\frac{1}{2}}\dfrac{1}{1+t'}
E_1(t')dt'\leq\dfrac{1}{5}\da^{\frac{1}{2}}(1+\ln(1+t))^5\til{E}_1(t,u).
\ee
Note that $\mu\geq\dfrac{1}{10}$ in $W_{n-s}$, which means
$\mu\ge 1-\da^{\frac{1}{2}}M(1+\ln(1+t))\sim\dfrac{1}{10}$ and $\da^
{\frac{1}{2}}M(1+\ln(1+t))\leq 1$. Hence,
\bee\label{q132}
\int_{W_t^u\cap W_{n-s}}(\mu^{-1}L\mu)(1+t')^2\sdfai^2\leq
\dfrac{1}{5}(1+\ln(1+t))^4\til{E}_1(t,u).
\ee
\hs In the shock region, it follows from Proposition\ref{keymu1} that
\bee\label{q133}
\int_{W_t^u\cap W_{shock}}(\mu^{-1} L\mu)(1+t')^2\sdfai^2\les\int_{W_t^u\cap
W_{shock}}
-\dfrac{1+t'}{1+\ln(1+t')}\mu^{-1}\sdfai^2=-K(t,u).
\ee
Hence,
\bee
\bes
\int_{W_t^u}Q_{1,3}\leq\dfrac{1}{5}(1+\ln(1+t))^4\til{E}_1(t,u)-K(t,u).
\end{split}
\ee
This finishes the estimates for the two boxed terms. Next, $Q_{1,2},Q_{1,6}$, $Q_{0,2},Q_{0,6},Q_{0,7}$ are estimated directly by the definition of $E_{\a},F_{\a}$ ($\a=0,1$) as follows.
\begin{equation*}\label{qeasy}
\bes
\int_{W_t^u}Q_{1,2}&\leq\int_{W_t^u}\mu^{-1}\dfrac{1+\ln(1+t')}{1+t'}(1+t')^2\left[
(L\fai+\dfrac{1}{2}\til{tr}\chi\fai)^2+(\til{tr}\chi\fai)^2\right]\\
&\leq\int_0^uF_1(u')du'+\da(1+\ln(1+t))^4\til{E}_0(t,u),\\
\int_{W_t^u}Q_{1,6}&\les\int_0^t\dfrac{\da^{\frac{3}{2}}(1+\ln(1+t'))}
{(1+t')^2}E_1(t')dt',\quad \int_{W_t^u}Q_{0,2}\les\da^{-\frac{1}{2}}\dfrac{1+\ln(1+t)}{(1+t)^2}\int_0^uF_1(u')du',\\
\int_{W_t^u}Q_{0,6}&\les\da^{\frac{3}{2}}\int_0^uF_0(u')du',\\
\int_{W_t^u}Q_{0,7}&\les\int_{W_t^u}\dfrac{\da^{\frac{1}{2}}(1+\ln(1+t'))}
{1+t'}\mu^{-1}(L\fai)(\dl\fai)\\
&\les\int_{W_t^u}\dfrac{\da^{\frac{1}{2}}(1+\ln(1+t'))}
{(1+t')^2}\mu^{-1}(\dl\fai)^2+\da^{\frac{1}{2}}
(1+\ln(1+t'))\mu^{-1}\left[(L\fai+\dfrac{1}{2}\til{tr}\chi\fai)^2+(\til{tr}\chi\fai)^2
\right]\\
&\les\da^{\frac{1}{2}}\int_0^t\dfrac{1+\ln(1+t')}{(1+t')^2}E_0(t')dt'+\da^{\frac{1}{2}
}\dfrac{1+\ln(1+t)}{(1+t)^2}\int
_0^uF_1(u')du'+\da^{\frac{1}{2}}\til{E}_0(t,u).
\end{split}
\end{equation*}
Next, we turn to estimate $Q_{1,0},Q_{0,0}$ and $Q_{1,8}$. It follows from
\bee\label{mudefe}
\mu\de\fe=\mu\da^{ij}\dfrac{\pa\fai_j}{\pa x^i}=\eta\hat{T}^iT\fai_j+\mu\sg^{AB}X_A^jX_B\fai_j,
\end{equation}
\eqref{inhomogeneous1} and \eqref{inhomogeneous2} that\footnote{This estimate is the only reason why we require $\lam>\frac{3}{2}$.}
\bee\label{q10}
\bes
\int_{W_t^u}Q_{1,0}&\les\int_{W_t^u}\mu^{-1}\dfrac{1+t'}{\til{tr}\chi}(L\fai+\frac{1}{2}\til{tr}\chi\fai)\\
&\cdot\left[
\dfrac{\da^{\frac{3}{2}}M\ln(1+t)}{(1+t)^{\lam+1}}(|T\fai|+\mu|\sd\fai|)+\dfrac{a}{(1+t')^{\lam}}\left(\dl\fai+\mu L\fai+\frac{1}{1+t}\fai\right)+\dfrac{1}{(1+t)^{\lam+1}}\left(\fai+\frac{\fe}{1+t}\right)\right]\\
&\les\int_{W_t^u}\mu^{-1}\left[\da^{-\frac{1}{2}}(1+t')^{\frac{7}{2}-\lam}(L\fai+\frac{1}{2}\til{tr}\chi\fai)^2
+\da^{\frac{1}{2}}(1+t')^{\frac{1}{2}-\lam}((\dl\fai)^2+\mu^2(L\fai)^2)\right]\\
&+\int_{W_t^u}\mu^{-1}\left[\dfrac{\fai^2+\fe^2}{1+t}+\dfrac{1}{(1+t')^{\lam-1}}\dfrac{1+t'}{\til{tr}\chi}(L\fai+\frac{1}{2}\til{tr}\chi\fai)^2\right]\\
&\les\da^{\frac{1}{2}}\int_0^t\dfrac{1}{(1+t')^{\lam-\frac{1}{2}}}E_0(t')t'
+\da^{-\frac{1}{2}}\int_0^uF_1(u')du'+\da\til{E}_0.
\end{split}
\ee
\hs Similarly, it holds that
\bee\label{q00}
\int_{W_t^u}Q_{0,0}\les\int_0^t\dfrac{1}{(1+t')^{\lam}}E_0(t')dt'+\da\til{E}_0(t).
\ee
For $Q_{1,8}$, it follows from \eqref{eq1+t} that
\bee\label{q18}
\int_{W_t^u}Q_{1,8}\les\int_{W_t^u}\mu^{-1}\fai^2(L(\eta^{-2}\mu+\mu\til{tr}\chi
+\til{tr}\dc))\les\da E_0(t,u).
\ee
Finally, we treat $Q_{1,5},Q_{0,4}$ and $Q_{0,5}$. It follows from \eqref{estzetaetaA} that
\beeq
\int_{W_t^u}Q_{1,5}&\les&\int_{W_t^u}\da^{\frac{1}{2}}\mu^{-1}(L\fai)(\sd^A\fai)
(1+t')(1+\ln(1+t'))\\
&+&\int_{W_t^u}\mu^{-1}(L\fai)(\sd^A\fai)\mu(1+t')(1+\ln(1+t'))^{\frac{1}
{2}}\\
&=&I_1+I_2.
\eeq
Then,
\bee\label{q151}
\bes
I_2&\les\int_{W_t^u}\da^{-\frac{1}{2}}(1+t')(1+\ln(1+t'))^2(L\fai)^2+
\int_{W_t^u}\da^{\frac{1}{2}}(1+t')(1+\ln(1+t'))^{-1}\sdfai^2\\
&\les\da^{-\frac{1}{2}}\int_0^uF_1(u')du'+\da^{\frac{1}{2}}\til{E}_0(t)+\da^{\frac{1}{2}}\til{E}_1(t),
\end{split}
\ee
due to Lemma\ref{estmu}. Rewrite $I_1$ as
\[
I_1=\left(\int_{W_t^u\cap W_{shock}}+\int_{W_t^u\cap W_{n-s}}\right)\da^{\frac{1}{2}}\mu^{-1}(L\fai)(\sd^A\fai)
(1+t')(1+\ln(1+t')).
\]
In the non-shock region, it holds that
\bee\label{q152}
\int_{W_t^u\cap W_{n-s}}\da^{\frac{1}{2}}\mu^{-1}(L\fai)(\sd^A\fai)
(1+t')(1+\ln(1+t'))\les\da^{-\frac{1}{2}}\int_0^uF_1(u')du'+\da^{\frac{1}{2}}\til{E}_0(t)+\da^{\frac{1}{2}}\til{E}_1(t);
\ee
while in the shock region, it holds that
\bee\label{q153}
\bes
&\int_{W_t^u\cap W_{shock}}\da^{\frac{1}{2}}\mu^{-1}(L\fai)(\sd^A\fai)
(1+t')(1+\ln(1+t'))\\
&\les\da^{\frac{1}{2}}\int_{W_t^u}\mu^{-1}(1+t')(1+\ln(1+t'))^2(L\fai)^2+\da^{\frac{1}{2}}\int_{W_t^u\cap W_{shock}}\mu^{-1}\sdfai^2(1+t')(1+\ln(1+t'))^{-1}\\
&\les\da^{\frac{1}{2}}\int_0^uF_1(u')du'+\da^{\frac{1}{2}}
\til{E}_0(t,u)+\da^{\frac{1}{2}}K(t,u).
\end{split}
\ee
For $Q_{0,5}$, it holds that
\bee
\int_{W_t^u}Q_{0,5}=\int_{W_t^u}\Omega\mu^{-1}(L\fai)(\sd^A\fai)(\eta_A+\zeta_A)
+\int_{W_t^u}\mu^{-1}(L\fai)(\sd^A\fai)(\mu(\eta_A+\zeta_A)+\mu\sd_A\mu):=I_1+I_2.
\ee
As before,
\bee\label{q051}
I_2\les\da^{\frac{1}{2}}\int_0^t\dfrac{1+\ln(1+t')}{(1+t')^{\frac{3}{2}}}E_1(t')dt'+
\da^{\frac{1}{2}}\int_0^uF_1(u')du'+\da^{\frac{1}{2}}\til{E}_0(t,u).
\ee
Rewrite $I_1$ as
\bee
I_1=\left(\int_{W_t^u\cap W_{shock}}+
\int_{W_t^u\cap W_{n-s}}\right)\Omega\mu^{-1}(L\fai)(\sd^A\fai)(\eta_A+\zeta_A).
\ee
\hs Similarly, in the non-shock region,
\bee\label{q052}
\int_{W_t^u\cap W_{n-s}}\Omega\mu^{-1}(L\fai)(\sd^A\fai)(\eta_A+\zeta_A)\les\da^{\frac{1}{2}}\int_0^t\dfrac{1+\ln(1+t')}{(1+t')^{\frac{3}{2}}}E_1(t')dt'+
\da^{\frac{1}{2}}\int_0^uF_1(u')du'+\da^{\frac{1}{2}}\til{E}_0(t,u);
\ee
while in the shock region, it holds that
\begin{equation*}\label{q053}
\int_{W_t^u\cap W_{shock}}\Omega\mu^{-1}(L\fai)(\sd^A\fai)(\eta_A+\zeta_A)
\les\da^{\frac{1}{2}}\dfrac{1+\ln(1+t)}{1+t}\int_0^uF_1(u')du'+\da^{\frac{1}{2}}\til{E}_0(t,u)
+\da^{\frac{1}{2}}\dfrac{(1+\ln(1+t))^2}{(1+t)^2}K(t,u).
\end{equation*}
Similarly, for $Q_{0,4}$, one can obtain
\begin{equation*}\label{q04}
\int_{W_t^u}Q_{0,4}\les\da^{\frac{1}{2}}\int_0^t\dfrac{(1+\ln(1+t'))^2}{(1+t')^2}
E_0(t')dt'+\da^{\frac{1}{2}}\dfrac{(1+\ln(1+t))^2}{(1+t)^2}K(t,u)+\da^{\frac{1}{2}}\int_0^uF_0(u')du'.
\end{equation*}
\hs Denote $M(t)=(1+\ln(1+t))^{-4}$. Then, it follows from \eqref{energy1} that
\bee\label{energy2}
\bes
&E_0(t,u)+F_0(t,u)+M(t)(E_1(t,u)+F_1(t,u))\leq E_0(0,u)+E_1(0,u)+I(0)\\
&+\left|\int_{W_t^u}\sum_{i=1}^7Q_{0,i}\right|+M(t)\left|\int_{W_t^u}\sum_{i=1}^8
Q_{1,i}\right|+\left|\int_{W_t^u}\underbrace{f\cdot K_0\fai}_{Q_{0,0}}\right|
+M(t)\left|\int_{W_t^u}\underbrace{f\cdot K_1\fai}_{Q_{1,0}}\right|+M(t)I(t).
\end{split}
\ee
It follows from the estimates for $Q_{\a}(\a=0,1)$ above that for sufficiently small $\da$, the right hand side of\eqref{energy2} is bounded by
\bee
\bes
&\int_0^t\dfrac{\da^{-\frac{3}{4}}(1+\ln(1+t'))}{\sqrt{\ln(1+s)-
\ln(1+t')}}
\dfrac{M(t')E_1(t')(1+\ln(1+t'))^4}{(1+t')^2}dt'\\
&+\int_0^t\left(\da^{-\frac{1}{2}}\dfrac{1+\ln(1+t')}{(1+t')^2}+\dfrac{1}
{(1+t')^{\lam}}+\da^{\frac{1}{2}}\dfrac{(1+\ln(1+t'))^5}{(1+t')^{\frac{3}{2}}}
\right)M(t')E_1(t')dt'\\
&+\int_0^t\left(\da^{\frac{1}{2}}\dfrac{1+\ln(1+t')}{(1+t')^2}+\dfrac{1}{
(1+t')^{\lam}}+\da^{\frac{1}{2}}\dfrac{1}{(1+t')^{\lam-\frac{1}{2}}}\right)E_0(t')dt'-M(t)K(t,u)\\
&+\da^{\frac{1}{2}}\til{E}_0(t,u)+\dfrac{1}{5}\til{E}_1(t,u)
+\da^{-\frac{1}{2}}\int_0^uM(t)F_1(u')du'+\int_0^uF_0(u')du'.
\end{split}
\ee
\begin{remark}
Note that all the coefficients of $M(t')E_1(t')$ and $E_0(t')$ in the integral is integrable in time.
\end{remark}
Keeping only $M(t)E_1
(t,u)$ on the left hand side of \eqref{energy2} yields
\bee\label{me1}
M(t)E_1(t,u)\leq G(t)+D_0+\int_0^tC(t')\til{E}_1(t')dt'+\dfrac{1}{5}\til{E}_1(t,u),
\ee
where $D_0=\int_{\Si_0^u}(L\fai)^2+(\dl\fai)^2+|\sd\fai|^2$, $C(t')$ representing the coefficients of $M(t')E_1(t')$ and $G(t)$ representing the remaining terms. Since right hand side of \eqref{me1} is increasing in $t$, then
\bee\label{me2}
\til{E}_1(t,u)\les G(t)+D_0+\int_0^tC(t')\til{E}_1(t')dt'.
\ee
Applying Gronwall inequality to \eqref{me2} yields
\bee\label{tilE1}
\til{E}_1(t,u)\les G(t)+D_0,
\ee
where the right hand side of \eqref{tilE1} doesn't contain $\til{E}_1(t)$. Define
\[
\til{K}(t,u)=\sup_{t'\in[0,t]}M(t')K(t',u).
\]
Then, it follows from \eqref{tilE1} and \eqref{energy2} that
\bee\label{energy3}
\bes
&\til{E}_0(t,u)+\til{F}_0(t,u)+\til{E}_1(t,u)+\til{F}_1(t,u)+\til{K}(t,u)\les
D_0\\
&+\int_0^t\left(\da^{\frac{1}{2}}\dfrac{1+\ln(1+t')}{(1+t')^2}+\dfrac{1}{
(1+t')^{\lam}}+\da^{\frac{1}{2}}\dfrac{1}{(1+t')^{\lam-\frac{1}{2}}}\right)\til{E}_0(t')dt'\\
&+\da^{\frac{1}{2}}\til{E}_0(t,u)+\da^{-\frac{1}{2}}\int_0^u\til{F}_1(u')du'
+\int_0^u\til{F}_0(u')du'.
\end{split}
\ee
Keeping only $\til{E}_0(t,u)$ on the left hand side of \eqref{energy3} and applying the same approach as before yield
\bee
\bes
&\til{E}_0(t,u)+\til{F}_0(t,u)+\til{E}_1(t,u)+\til{F}_1(t,u)+\til{K}(t,u)\les
D_0\\
&+\da^{-\frac{1}{2}}\int_0^u\til{F}_1(u')du'
+\int_0^u\til{F}_0(u')du'.
\end{split}
\ee
Similarly argument leads to
\bee\label{energyfinal}
\til{E}_0(t,u)+\til{F}_0(t,u)+\til{E}_1(t,u)+\til{F}_1(t,u)+\til{K}(t,u)\les
D_0.
\ee
\begin{remark}\label{energyremark}
\textbf{To derive the energy estimates for  high order derivatives of $\fai$, the corresponding error integrals $Q_{0,i},Q_{1,i}$ (for $i\ge 1$) can be estimated as above exactly by high order energies and fluxes directly. However, the treatment for $Q_{0,0},Q_{1,0}$ will involve complicated high order acoustical terms and will be given details later.}
\end{remark}
\subsection{\textbf{Some top order acoustical terms}}
\hs We compute the top order acoustical terms by following the framework in \cite{CZD1} section5.2. For the wave equation \eqref{nonlinearwave}: $\Box_{\til{g}}\fai=f$, let $\fai_1=\fai,\p_1=f, J_1=J$ and $\fai_n=Z\fai_{n-1}$. Then, the following relations hold:
\bee\label{recursion1'}
\bes
\til{\p}_n&=\Omega^2\mu\p_n=Z\til{\p}_{n-1}+\left(\dfrac{1}{2}tr_{\til{g}}\zgpi
-\mu^{-1}Z\mu-2Z\log \Omega\right)\cdot\til{\p}_{n-1}+\mu div_g\pre Z {} J_{n-1}\\
&=Z\til{\p}_{n-1}+\pre {Z} {} {\da}\cdot\til{\p}_{n-1}+\pre Z {} {\sigma_{n-1}},\\
\til{\p}_1&=\Omega^2\mu\left(-\dfrac{2\eta'}{\eta^2}A(\dfrac{a}{(1+t)^{\lam}}\fe)\de\fe+\dfrac{1}{\p\eta}A\left[\dfrac{a}{(1+t)^{\lam}}\left(\dfrac{\pa\fe}{\pa t}-\frac{1}{2}|\nabla\fe|^2-\frac{\lam}{(1+t)}\fe\right)\right]\right),
\end{split}
\ee
where
\bee
\pre {Z} {} {\da}=\dfrac{1}{2}tr_{\til{g}}\zgpi
-\mu^{-1}Z\mu-2Z(\log \Omega),\quad \pre Z {} {\sigma_{n-1}}=\mu div_g\pre Z {} J_{n-1}.
\ee
It follows from the estimates for the deformation tensors \eqref{estqgpi}, \eqref{esttgpi} and \eqref{estrgpi} and the relation $tr_{\til{g}}\zgpi=-\mu^{-1}\Omega^{-1}\zgpi_{L\dl}+tr_{\til{g}}\szgpi$ that $|\pre {Z} {} {\da}|\les\da$ for $Z=R_i, Q$, $|\pre {Z} {} {\da}|\les 1$ for $Z=T$.
\begin{lem}\label{recur}
Let $y_n (n=1,2,...)$ be a sequence 
and $A_n$ be a given sequence of operators. Suppose that $x_n$ is a sequence satisfying: $x_n=A_{n-1}x_{n-1}+y_{n-1}.$ Then,
\[
x_n=A_{n-1}A_{n-2}\cdots A_1 x_{1}+\sum_{m=0}^{n-2}A_{n-1}\cdots A_{n-m}y_{n-m-1}.\]
\end{lem}
Applying this lemma to \eqref{recursion1'} yields
\bee\label{recursion2}
\til{\p}_n=(Z_{n-1}+\pre {Z_{n-1}} {} {\da})\cdots(Z_{1}+\pre {Z_{1}} {} {\da})\til{\p}_{1}+\sum_{k=0}^{n-2}(Z_{n-1}+\pre {Z_{n-1}} {} {\da})\cdots(Z_{n-k}+\pre {Z_{n-k}} {} {\da})\pre {Z_{n-k-1}} {} {\sigma_{n-k-1}}.
\ee
One has the following decomposition\footnote{For the detailed expressions, see\cite{CZD1} section5.2.}:
\bee\label{decomsi}
\pre Z {} {\sigma_{n-1}}=\pre Z {} {\sigma_{1,n-1}}+\pre Z {} {\sigma_{2,n-1}}+\pre Z {} {\sigma_{3,n-1}},
\ee
\textbf{where $\pre Z {} {\sigma_{2,n-1}}$ is the major term which contains the product of $1^{st}$ order derivatives of both deformation tensors of $Z$ and $\fai_{n-1}$, $\pre Z {} {\sigma_{1,n-1}}$ contains the product of deformation tensors of $Z$ and $2^{nd}$ order derivatives of $\fai_{n-1}$ and $\pre Z {} {\sigma_{3,n-1}}$ contains the rest which is easy to handle.}\\
\hs In the paper, the highest order of an object will be $N_{top}+1$ (which is called the top order) and we consider the case of $n=N_{top}+1$ in \eqref{recursion1'} with $\fai_n=Z_{n-1}\cdots Z_1\fai$.\\
\hs Then we will consider the top order acoustical terms in $\til{\p}_n$ associated with the different choices of $\fai_n$ and there are following 3 possibilities of $\fai_n$ :
\[
\fai_n=R_i^{\a+1}\fai,\quad \fai_n=R_i^{\a'}T^{l+1}\fai,\quad \fai_n=R_i^{\a'}T^{l}Q\fai,
\]
where $N_{top}=|\a|+1=|\a'|+l+1$. It follows from the discussion in \cite{CZD1} section 5.2 that the possible top order acoustical terms will be contained in $Z_{n-1}\cdots Z_2\pre {Z_1} {} {\sigma_{2,1}}$.\\
 \begin{itemize}
 \item[Case1:] For $Z=Q$, it holds that\footnote{The notation $[\quad]_{P.A.}$ means the principle acoustical terms of an entity, which only keeps the spatial derivatives of $\mu$ and $\chi$}
     \beeq
     &&[tr\sqgpi]_{P.A.}=2\Omega (1+t)\til{tr}\chi,\hs [\qgpi_{L}]_{P.A.}=0,\hs
     [\qgpi_{L\dl}]_{P.A.}=0,\\
     &&[\qgpi_{\dl}]_{P.A.}=2\Omega (1+t)\sd\mu,\hs [\mu\pre {(Z)} {} {\hat{\til{\s{\pi}}}}]_{P.A.}=2\mu\Omega (1+t)\hat{\chi},
     \eeq
     and then\footnote{The structure equation\eqref{Ttrchi} is used.}
     \bee
     \bes
     [\pre {(Q)} {} {\sigma_{2,1}}]_{P.A.}&=[\Omega (1+t)L\fai(T\til{tr}\chi-\s{\de}
     \mu)+2\mu\Omega (1+t)\s{div}\hat{\chi}\cdot\sd\fai]_{P.A.}\\
     &=2\mu \Omega (1+t)\sd tr\chi\cdot\sd\fai.
     \end{split}
     \ee
     Hence, for the variation $R_i^{\a'}T^lQ\fai$, i.e. $Z_1=Q,$ $Z_{n-1}\cdots Z_2=R_i^{\a'}T^{l}$, the corresponding top order acoustical terms are
     \bee\label{paq}
     \left\{\bes
     &\mu\Omega(1+t)\sd R_i^{\a}tr\chi\cdot\sd\fai,\hs if\ l=0,\\
     &\mu\Omega(1+t)R_i^{\a'+1}T^{l-1}\s{\de}\mu\cdot\sd\fai,\hs if\ l\ge 1.
     \end{split}\right.
     \ee
 \item[Case2:] For $Z=T$, it holds that
       \bee
       \bes
       [\pre {(T)} {} {\sigma_{2,1}}]_{P.A.}&=[-\Omega\eta^{-2}\mu L\fai(T\til{tr}\chi-
       \s{\de}\mu)+\Omega T\fai\cdot\s{\de}\mu-2\Omega\eta^{-2}\mu^2\s{div}\hat{\chi}\cdot\sd\fai]_{P.A.}\\
       &=\Omega T\fai\cdot\s{\de}\mu-2\Omega\eta^{-2}\mu^2\s{div}\hat{\chi}\cdot\sd\fai.
       \end{split}
       \ee
       Hence, for the variation $R_i^{\a'}T^{l+1}\fai$, the corresponding top order acoustical terms are
       \bee\label{pqt}
       \left\{\bes
       &\Omega T\fai\cdot R_i^{\a}\s{\de}\mu-2\Omega\eta^{-2}\mu^2\sd\fai\cdot\sd R_i^{\a}tr\chi,\hs if\ l=0,\\
       &\Omega T\fai\cdot R_i^{\a'}T^l\s{\de}\mu-2\Omega\eta^{-2}\mu^2\sd\fai\cdot R_i^{\a'+1} T^
       {l-1}\s{\de}\mu,\hs if\ l\ge 1.
       \end{split}\right.
       \ee
       \item[Case3:] For $Z=R_i$, it holds that
       \bee
       \bes
       [\pre {(R_i)} {} {\sigma_{2,1}}]_{P.A.}&=[\Omega\mu^{-1}\lam_i L\fai(T\til{tr}\chi-\s{\de}\mu)+\Omega R_i\fai(T tr\chi-\s{\de}\mu)\\
       &+\Omega T\fai\cdot R_i tr\chi+2\eta^{-1}\mu\Omega\lam_i\sd\fai\cdot\sd tr\chi]_{P.A.}\\
       &=\Omega T\fai\cdot R_i tr\chi+2\eta^{-1}\mu\Omega\lam_i\sd\fai\cdot\sd tr\chi.
       \end{split}
       \ee
       Hence, for the variation $R_i^{\a+1}\fai$, the corresponding top order acoustical term is
       \bee\label{par}
       \Omega T\fai\cdot\sd R_i^{\a}tr\chi+2\eta^{-1}\mu\Omega\lam_i\sd\fai\cdot
       \sd R_i^{\a}tr\chi.
       \ee
 \end{itemize}
 \hs Note that $|\sd\fai|\les\dfrac{\da^{\frac{3}{2}} M}{(1+t)^{2}}$ and $|T\fai|\les \dfrac{\da^{\frac{1}{2}}M}{1+t}$. Then, it follows from \eqref{paq}, \eqref{pqt} and \eqref{par} that to estimate the top order acoustical terms, one needs only to estimates
 \[
 T\fai\cdot \s{d}R_i^{a}tr\chi,\hs T\fai\cdot R_i^{\a'}T^l\s{\de}\mu,
 \]
 where $|\a|+1=|\a'|+l+1=N_{top}$.
 \subsection{\textbf{Estimates for lower order terms}}
 \hs In this subsection, we will give the estimates for lower order terms (the terms of order $\leq N_{top}$), which will be necessary for the estimates of top order acoustical terms. To this end,  
 it will be shown later that the key role in this subsection is to estimate the difference between the acoustical coordinates and the rectangular coordinates.\\
 \hs We define the modified energies and fluxes as follows. First, define:
 \[
 \wi{E}_{0,k+1}=\sum_{\fai}\sum_{|\a|=k}\da^{2l}\til{E}_0(Z^{\a}\fai),
 \]
 which means the sum of energies of all the $n$-th order variations and similar for $\wi{F}_{0,k+1},\wi{E}_{1,k+1},\wi{F}_{1,k+1}$. For $k=0$, we set $\wi{E}_{\a,1}=\til{E}_{\a}$ and $ \wi{F}_{\a,1}=\til{F}_{\a}$ ($\a=0,1$).\\
\hs Let $0\leq k\leq N_{top}$ and $(b_k)$ be a sequence of nonnegative integers to be chosen later. Recall $\mu_m^u$ defined by \eqref{defmumu} and define:
\bee\label{modifiedenergy}
\bes
\overline{E}_{0,k+1}(t, u )&=\sup_{\tau\in [0,t]}\{(1+\ln(1+t))^{-2p}\mu^{ u }_m(\tau)^{2b_{k+1}}\wi{E}_{0,k+1}(\tau, u )\},\\
\overline{F}_{0,k+1}(t, u )&=\sup_{\tau\in [0,t]}\{(1+\ln(1+t))^{-2p}\mu^{ u }_m(\tau)^{2b_{k+1}}\wi{F}_{0,k+1}(\tau, u )\},\\
\overline{E}_{1,k+1}(t, u )&=\sup_{\tau\in [0,t]}\{(1+\ln(1+t))^{-2q}\mu^{ u }_m(\tau)^{2b_{k+1}}\wi{E}_{1,k+1}(\tau, u )\},\\
\overline{F}_{1,k+1}(t, u )&=\sup_{\tau\in [0,t]}\{(1+\ln(1+t))^{-2q}\mu^{ u }_m(\tau)^{2b_{k+1}}\wi{F}_{1,k+1}(\tau, u )\},
\end{split}
\ee
where $l$ is the number of $T's$ in $Z^{\a+1}\in\{Q,T,R_i\}$ and $q\ge p$ are two constants.
 \begin{lem}
Let $A(t,u)$ be the area of $S_{t,u}$. Then there is a universal constant $C$ such that
\[
C^{-1}(1+t)\leq A(t,u)\leq C(1+t).
\]
As a corollary, for any function $f$, it holds that
\bee
C^{-1}(1+t)||f||_{L^2([0,u]\times S^2)}\leq ||f||_{L^2(\Si_{t,u})}\leq C(1+t)||f||_{L^2([0,u]\times S^2)}.
\ee
\end{lem}
\begin{pf}
It follows from the definition of $\chi$ and Lemma\ref{2ndff} that
\[
|tr\chi|=|\frac{1}{2}\lie_L\sg|=|\dfrac{1}{\sqrt{\det\sg}}(L\sqrt{\det\sg})|\les\dfrac{1}{1+t}.
\]
Denote $\mu_{\sg}(0,0)$ be the area form of the sphere $S_{0,0}$.
Then,
\bee\label{dmusgdmu0}
d\mu_{\sg}=\dfrac{\sqrt{\det\sg(t,u)}}{\sqrt{\det\sg(0,0)}}d\mu_{\sg}(0,0)=
\dfrac{\sqrt{\det\sg(t,u)}}{\sqrt{\det\sg(0,u)}}\dfrac{\sqrt{\det\sg(0,u)}}{\sqrt{\det\sg(0,0)}}d\mu_{\sg}(0,0).
\ee
Note that
\bee
|\ln\dfrac{\sqrt{\det\sg(t,u)}}{\sqrt{\det\sg(0,u)}}|=|\int_{0}^t\dfrac{1}{\sqrt{\det\sg}}(L\sqrt{\det\sg})dt'|\les \ln(1+t).
\ee
Then, integrating \eqref{dmusgdmu0} over $S^2$ yields
\bee
C^{-1}(1+t)\int_{S^2}d\mu_{\sg(0,0)}\leq \int_{S^2}d\mu_{\sg(t,u)}\leq C(1+t)\int_{S^2}d\mu_{\sg(0,0)}.
\ee
That is, $C^{-1}(1+t)\leq A(t,u)\leq C(1+t)$.
\end{pf}
\hs Given any $k$ and $j$, define $x^j_{i_1\cdots\i_k}$ and $\pre {(k)} {} {\da}^j_{i_i\cdots i_k}$ as
\[
 R_{i_k}\cdots R_{i_1}x^j=\mari{R}_{i_k}\cdots\mari{R}_{i_1}x^j-\pre {(k)} {} {\da}^j_{i_i\cdots i_k}=x^j_{i_1\cdots\i_k}-\pre {(k)} {} {\da}^j_{i_i\cdots i_k},
\]
where for $k=0$, we set $\pre {(0)} {} {\da}^j=0$.
Obviously, $x^j_{i_1\cdots\i_k}$ are linear functions of the rectangular coordinates and in particular, $\dfrac{\pa(x^j_{i_1\cdots\i_k})}{\pa x^l}:=\pre {(k)} {} {c}^j_{l,i_i\cdots i_k}$ are constants. $\pre {(k)} {} {\da}^j_{i_i\cdots i_k}$ can be computed as 
\bee\label{recursionda}
\pre {(k)} {} {\da}^j_{i_i\cdots i_k}=R_{i_k}\pre {(k-1)} {} {\da}^j_{i_i\cdots i_{k-1}}+\pre {(k-1)} {} {c}^j_{l,i_i\cdots i_{k-1}}\lam_{i_k}\hat{T}^{l}.
\ee
Applying Lemma\ref{recur} to \eqref{recursionda} yields
\bee\label{recursionda2}
\pre {(k)} {} {\da}^j_{i_i\cdots i_k}=\sum^{k}_{m=1}\pre {(m-1)} {} {c}^j_{l,i_1\cdots i_{m-1}}R_{i_k}\cdots R_{i_{k-m+2}}(\lam_{i_{m}}\hat{T}^l).
\ee
It follows from the definition that $x^j_{i_1\cdots\i_k}$ are linear functions of $x^j$ with uniform coefficients, so can be bounded by $r$, then by $C(1+t)$. One has the following lemma.
\begin{lem}\label{RikRi1yjLinfty}
For $k=0,1,\cdots N_{\infty}-1$, if $||R_{i_k}\cdots R_{i_i}y^j||_{\supnormda}\les\dfrac{\da^{\frac{1}{2}}M(1+\ln(1+t))}{1+t}$, then the following estimates hold for sufficiently small $\da$
\begin{equation*}
||R_{i_k}\cdots R_{i_i}\lam_i||_{\supnormda}\les\da^{\frac{3}{2}}M(1+\ln(1+t)),\quad ||\pre {(k+1)} {} {\da}^j_{ii_i\cdots i_k}||_{\supnormda}\les\da^{\frac{3}{2}}M(1+\ln(1+t)).
\end{equation*}
\end{lem}
\begin{remark}
Since $|x^j_{i_1\cdots i_k}|\les (1+t)$, then one can bound $R_{i_k}\cdots R_{i_i}x^j$ by $C(1+t)$ provided that the assumption of Lemma\ref{RikRi1yjLinfty} holds.
\end{remark}
\hs It remains to estimate $||R_i^{\a}y^j||_{\supnormda} $ for $|\a|=k+1\leq N_{\infty}-1$, and one notes that\footnote{The detailed analysis is given in \cite{CZD1} (5.48)-(5.51).}
\bee\label{RiayjlieRibechi'}
R_i^{\a}y^j=O(1)\cdot \lie_{R_i}^{\be}\chi'+O_{\frac{1}{1},\frac{1}{2}}^{\leq |\a|},
\ee
due to Lemma\ref{RikRi1yjLinfty} and an induction argument. Then, it suffices to estimate $\lie_{R_i}^{\be}\chi'$ with $|\be|\leq N_{\infty}-2$.
\begin{lem}\label{Riachi'Linfty}
For sufficiently small $\da$, and all $|\a|\leq  N_{\infty}-2=[\dfrac{N_{top}}{2}]+1$, it holds that
\bee
Z_i^{\a+1}y^j\in O^{|\a|+1}_{\frac{1}{1},1-2l},\hs\dfrac{1}{1+t}Z_i^{\a+1}x^j\in O^{|\a|+1}_{\frac{0}{0},-2l},\hs \lie_{Z_i}^{\a}\chi'\in O^{|\a|+1}_{\frac{1}{2},3-2l},
\ee
where $l$ is the number of $T$'s in $Z_i^{\a}$.
\end{lem}
\begin{pf}
We will give the proof for $\lie_{R_i}^{\a}\chi'\in O^{|\a|+1}_{\frac{1}{2},3}$, and the proofs for the others are similar.\\
The case for $|\a|=0$ follows from Lemma\ref{2ndff}. For the case $|\a|\geq1$, we will prove by an induction argument as follows. Commuting \eqref{Lchi'AB} with $\lie_{R_i}^{\a}$ yields
\bee\label{zax}
\bes
\lie_{L}\lie_{R_i}^{\a}\chi'&=[\lie_{L},\lie_{R_i}^{\a}]\chi'+(e+a\frac{\eta^{-1}\hat{T}^j\fai_j}{(1+t)^{\lam}})\cdot\lie_{R_i}^{\a}\chi'+2\chi'\cdot\lie_{R_i}^{\a}\chi'\\
&+\sum_{|\be_1|+|\be_2|=|\a|,|\be_1|>0}R_i^{\be_1}e\cdot\lie_{R_i}^{\be_2}\chi'\\
&+\sum_{|\be_1|+|\be_2|+|\be_3|=|\a|,|\be_2|>0}
\lie_{R_i}^{\be_1}\chi'\cdot\lie_{R_i}^{\be_2}\sg\cdot\lie_{R_i}^{\be_3}\chi'\\
&+\sum_{|\be_1|+|\be_2|=|\a|,|\be_1|>0}R_i^{\be_1}(a\frac{\eta^{-1}\hat{T}^j\fai_j}{(1+t)^{\lam}})\cdot\lie_{R_i}^{\be_2}\chi'\\
&+\lie_{R_i}^{\a}\left(\dfrac{\sg_{AB}}{1- u +t}\left(e+a\frac{\eta^{-1}\hat{T}\fe}{(1+t)^{\lam}}\right)-\a'_{AB}\right).
\end{split}
\ee
The right hand side of \eqref{zax} can be estimated with the same argument in proof of Lemma5.8 \cite{CZD1},
and then
\bee
\lie_{L}\lie_{R_i}^{\a}\chi'=(e+a\frac{\eta^{-1}\hat{T}^j\fai_j}{(1+t)^{\lam}}+\dfrac{\da^{\frac{1}{2}}M(1+\ln(1+t))}{(1+t)^2})\cdot\lie_{R_i}^{\a}\chi'
+\chi'\cdot\lie_{R_i}^{\a}\chi'+O^{|\a|+2}_{\frac{0}{3},3},
\ee
which implies that
\bee\label{lieLlieRiachi'}
\bes
|\lie_L\left((1-u+t)^2\lie_{R_i}^{a}\chi'\right)|&\les
(1-u+t)^2|\lie_L\lie_{R_i}^{\a}\chi'|
+2(1-u+t)|\chi'|\cdot|\lie_{R_i}^{\a}\chi'|\\
&\les\dfrac{\da^{\frac{1}{2}}M(1+\ln(1+t))}{(1+t)^2}|(1-u+t)^2\lie_{R_i}^{\a}\chi'|
+O^{|\a|+2}_{\frac{0}{1},3}.
\end{split}
\ee
Let $f=|(1-u+t)^2\lie_{R_i}^{a}\chi'|$. Then, \eqref{lieLlieRiachi'} can be rewritten as
\bee
\dfrac{df}{dt}\les\dfrac{\da^{\frac{1}{2}}M(1+\ln(1+t))}{(1+t)^2}f+\dfrac{\da^{
\frac{3}{2}}M}{1+t},\\
\ee
which implies that
\bee
|f|\les\da^{\frac{3}{2}}M(1+\ln(1+t)).
\ee
That is,
\bee
||\lie_{R_i}^{\a}\chi'||_{\supnormda}
\les\dfrac{\da^{\frac{3}{2}}M(1+\ln(1+t))}{(1+t)^2}.
\ee
\end{pf}
\begin{remark}
Since $\chi=\chi'+\dfrac{\sg}{1-u+t}$, then for $|\a|\ge 1$, $||\lie_{R_i}^{\a}\chi||_{\supnormda}\les\dfrac{\da^{\frac{3}{2}}M
(1+\ln(1+t))}{(1+t)^2}.$
\end{remark}

We now turn the $L^{\infty}$ estimate for $\mu.$
\begin{prop}\label{ZiamuLinfty}
For sufficiently small $\da$ and all $|\a|\leq N_{\infty}-2$, it holds that
\[||Z_i^{\a+1}(\mu-1)||_{\supnormda}\les\da^{\frac{1}{2}-l}M(1+\ln(1+t)),
\]
where $l$ is the number of $T$'s in $Z_i^{\a+1}$.
\end{prop}
\begin{pf}
For the case $Z_i^{\a+1}=R_i^{\a+1}$ or $Z_i^{\a+1}=R_i^{\be+1}Q^{\a-\be}$, then the proof will be same as in \cite{CZD1} Proposiiton5.2. However, when $Z_i^{\a+1}=R_i^{\a'}T_i^{m+1}$ with $|\a'|+m=|\a|$, things will be more complicated since the estimates for $T^m\fe$ for $m>1$ are required. In fact, \textbf{it can be shown by a bootstrap argument and an induction argument that $||R_i^{\a'}T^{m+1}\fe||_{\supnormda}\les\da^{-m+\frac{1}{2}}M$ for all $|\a'|+m\leq N_{\infty}-1$ and as a corollary, $||R_i^{\a'}T^{m}\mu||_{\supnormda}\les\da^{\frac{1}{2}-m}M(1+\ln(1+t))$}, which will be sketched as follows (see also Lemma\ref{ltfe}).\\
\hs Let $|\a'|+m=N_{\infty}-1$ and we will proceed with induction on $m$. For $m=0$, it follows from Lemma\ref{LTidTi}, the bootstrap assumptions and the similar argument in Lemma\ref{ltfe} that $|R_i^{\a'}T\fe|\les\dfrac{\da^{\frac{3}{2}}M^2(1+\ln(1+t))}{1+t}\da^{\frac{1}{2}}M<\da^{\frac{1}{2}} M$ and $|R_i^{\a'}\mu|\les \da^{\frac{1}{2}}M(1+\ln(1+t))$.\\
\hs Assume that $|R_i^{\be}T^{k+1}\fe|\les\da^{-k+\frac{1}{2}}M$ and $|R_i^{\be}T^{k}\mu|\les\da^{-k+\frac{1}{2}}M(1+\ln(1+t))$ for all $|\be|+k=N_{\infty}-1$ and $0\leq k\leq m-1$. Suppose that $|R_i^{\a'}T^{m+1}\fe|\les\da^{-m+\frac{1}{2}}M$. Then, it follows that
\bee
\bes
R_i^{\a'}T^{m+1}\fe&=R_i^{\a'}T^{m}(\eta^{-1}\mu\hat{T}^i\fai_i)\\
&=\eta^{-1}R_i^{\a'}T^{m}(\mu)\hat{T}^i\fai_i+\eta^{-1}\mu\fai_iR_i^{\a'}T^{m}(\hat{T}^i)+O^{\leq N_{\infty}-1}_{\frac{1}{1},2-2m}.
\end{split}
\ee
\been
\item For $R_i^{\a'}T^m\mu$, commuting $R_i^{\a'}T^m$ with \eqref{transportmu} yields
\bee\label{tkmu}
\bes
LR_i^{\a'}T^m\mu&=R_i^{\a'}T^mm+eR_i^{\a'}T^m\mu+\sum_{|\a_1|+|\a_2|=|\a'|,\be_1+\be_2=m,|\a_2|+\be_1>0}R_i^{\a_1}T^{\be_1}e\cdot R_i^{\a_2}T^{\be_2}\mu\\
&+\sum_{\be_1+\be_2=|\a'|-1}R_i^{\be_1}\srpi_L\cdot\sd R_i^{\be_2-1}T^m\mu+\sum_{\be_1+\be_2=m-1}R_i^{\a'}T^{\be_1}\stpi_{L}\cdot\sd T^{\be_2-1}\mu.
\end{split}
\ee
It holds that $|R_i^{\a'}T^mm|\les|R_i^{\a'}T^{m+1}\fai|+|R_i^{\a'}T^{m+1}\fe|\les\da^{-m+\frac{1}{2}}M$. By induction, \eqref{tkmu} can be rewritten as
$LR_i^{\a'}T^{m}\mu=(e+\dfrac{\da^{\frac{1}{2}}M(1+\ln(1+t))}{(1+t)^2})R_i^{\a'}T^{m}\mu+O^{\leq N_{\infty}}_{\frac{0}{1},2-2k}$. Then, integrating along integral curves of $L$ yields $|R_i^{\a'}T^m\mu|\les\da^{\frac{1}{2}-m}M(1+\ln(1+t))$.
\item For $R_i^{\a'}T^m(\hat{T}^i)$, it follows from Lemma\ref{LTidTi} that
\bee
R_i^{\a'}T^m(\hat{T}^i)=R_i^{\a'}T^{m-1}\left(-X_B(\eta^{-1}\mu)\cdot\s{d}_Bx^i\right).
\ee
Then, $|R_i^{\a'}T^{m}(\hat{T}^i)|=|O(1)R_i^{\a'+1}T^{m-1}\mu|+\text{l.o.ts}\les\dfrac{1+\ln(1+t)}{1+t}\da^{\frac{3}{2}-m}M$.
\een
Hence, $|R_i^{\a'}T^{m+1}\fe|\les\dfrac{\da^{-m+\frac{1}{2}}M^2(1+\ln(1+t))}{1+t}<\da^{-m+\frac{1}{2}}M$ provided that $\da$ is sufficiently small, which recovers the assumption $|R_i^{\a'}T^{m+1}\fe|\les\da^{-m+\frac{1}{2}}M$.
\end{pf}
\begin{prop}
It follows from \eqref{computationqpi}, \eqref{computationtpi} and \eqref{computationrpi1} that for all $|\a|\leq N_{\infty}-2$
\bee
\bes
&\lie_{Z_i}^{\a}\stpi_{L}\in O^{|\a|+1}_{\frac{1}{1},1-2l},\hs
\lie_{Z_i}^{\a}\stpi_{\dl}\in O^{|\a|+1}_{\frac{2}{1},2-2l},\\
&\lie_{Z_i}^{\a}\stpi\in O^{|\a|+1}_{\frac{2}{2},4-2l},\hs
\lie_{Z_i}^{\a}\srpi_L\in O^{|\a|+1}_{\frac{1}{1},1-2l},\\
&\lie_{Z_i}^{\a}\srpi_{\dl}\in O^{|\a|+1}_{\frac{2}{1},2-2l},\hs
\lie_{Z_i}^{\a}\srpi\in O^{|\a|+1}_{\frac{1}{1},3-2l}.
\end{split}
\ee
\end{prop}
\hs Next we turn to $L^2$ estimate for the above entities, which follows from the same framework in deriving $L^{\infty}$ estimates.
\begin{lem}
For $k=0,1,2...[\frac{l}{2}]:=[\frac{N_{top}}{2}]$, if $\max_{j,i_1\cdots i_k}||R_{i_k}\cdots R_{i_1}y^j||_{\supnormu}\les
\dfrac{\da^{\frac{1}{2}}M(1+\ln(1+t))}{1+t}$, then for sufficiently small $\da$, it holds that
\bee
\bes
&\max_{j,i_1\cdots i_k}||R_{i_k}\cdots R_{i_i}\lam_j||_{\normu2}\les(1+t)
\sum_{k=0}^l||R_{i_k}\cdots R_{i_1}y^{j}||_{\normu2},\\
&\max_{j,i_1\cdots i_k}||R_{i_k}\cdots R_{i_i}\pre {(k+1)} {} {\da}^j_{i,i_1\cdots i_k}||_{\normu2}\les(1+t)
\sum_{k=0}^l||R_{i_k}\cdots R_{i_1}y^{j}||_{\normu2}.
\end{split}
\ee
\end{lem}
\hs Next, we estimate $||\lie_{R_i}^{\a}y^j||_{\normu2}$ with $|\a|\leq N_{top}$. It suffices to estimate $||\lie_{R_i}^{\be}\chi'||_{\normu2}$ with $|\be|\leq N_{top}-1$ due to \eqref{RiayjlieRibechi'} and an induction argument. In fact, one has\footnote{The proof is the same as Proposition5.3 in\cite{CZD1}.}:
\begin{prop}\label{Riachi'L2}
For sufficiently small $\da$ and all $|\a|\leq N_{top}-1$, it holds that
\begin{equation*}
||\lie_{R_i}^{\a}\chi'||_{\normu2}\les\dfrac{1}{1+t}\left(\sum_{|\a'|\leq|\a|}
||\lie_{R_i}^{\a'}\chi'||_{L^2(\Si_0^{\tilde{\da}})}+\int_0^t
\dfrac{\da^{\frac{3}{2}}M(1+\ln(1+t'))}{1+t'}\mu_m^{-\frac{1}{2}}\sqrt{
\wi{E}_{1,|\a|+2}(t',u)} dt'\right).
\end{equation*}
\hs As a corollary, it holds that for $|\a|\leq N_{top}-1$
    \bee
    \bes
    ||R_i^{\a+1}y^j||_{\normu2}&\les ||\lie_{R_i}^{\a}\chi'||_{\normu2},\\
    ||R_i^{\a+1}\lam_j||_{\normu2}&\les (1+t)||\lie_{R_i}^{\a}\chi'||_{\normu2}.
    \end{split}
    \ee
\end{prop}
    \begin{prop}\label{Ria+1muL2}
    For sufficiently small $\da$ and all $|\a|\leq N_{top}-1$, it holds that
    \begin{equation*}
    (1+t)^{-1}||R_i^{\a+1}\mu||_{\normu2}\les \sum_{|\a'|\leq|\a|}(||R_i^{\a'+1}\mu||_{L^2(\Si
    _0^{\tilde{\da}})}+||\lie_{R_i}^{\a'}\chi'||_{L^2(\Si_0^{\tilde{\da}})})+\int_0^t\dfrac{1}{1+t'}\mu_m^{-\frac{1}{2}}\sqrt{\wi{E}_{1,|\a|+2}(t',u)}dt'.
    \end{equation*}
    \end{prop}
    \begin{pf}
    Commuting $R_i^{\a+1}$ with \eqref{transportmu} yields
    \bee
    LR_i^{\a+1}\mu=[L,R_i^{\a+1}]\mu+R_i^{\a+1}m+e\cdot R_i^{\a+1}\mu+\sum_{|\be_1|+|\be_2|=|\a|+1,|\be_1|>0}R_i^{\be_1}e\cdot
    R_i^{\be_2}\mu.
    \ee
    One can use similar argument to estimate $||R_i^{\be_1}e\cdot
    R_i^{\be_2}\mu||_{\normu2}$ in estimating $L^2$ norm of $\chi'$ \footnote{That is, if $\be_1\geq[\frac{N_{top}}{2}]+1=N_{\infty}-2$, we bound $R_i^{\be_1}e$ in $L^2$ norm while $R_i^{\be_2}$ can be bounded by Proposition\ref{ZiamuLinfty}; while if $\be_2\geq N_{\infty}-2$, we bounded $R_i^{\be_1}e$ by bootstrap assumptions and $||R_i^{\be_2}\mu||_{\normu2}$ can be absorbed by Gronwall inequality as shown later.}, and it holds that
    \bee
    \bes
    (1+t)^{-1}||R_i^{\be_1}e\cdot R_i^{\be_2}\mu||_{\normu2}&\les
    \dfrac{\da^{\frac{3}{2}}M}{(1+t)^2}\dfrac{1}{1+t}||R_i^{\a}\mu||_{\normu2}\\
    &+\dfrac{\da^{\frac{1}{2}}M(1+\ln(1+t))}{(1+t)^2}\left(
    \mu_m^{-\frac{1}{2}}\sqrt{\wi{E}_{1,|\a|+1}(t,u)}+||\lie_{R_i}^{\a-1}\mu||_{\normu2}\right).
    \end{split}
    \ee
    \hs As for $R_i^{\a+1}m$, it holds that
    \bee
    \bes
    |R_i^{\a+1}m|&=|\dfrac{1}{2}R_i^{\a+1}(\dfrac{dH}{dh}Th)+\dfrac{R_i^{\a+1}T\fe}{(1+t)^{\lam}}|\\
    &\les |R_i^{\a+1}T\fai_0-R_i^{\a+1}(\fai_iT\fai_i)+\dfrac{a}
    {(1+t)^{\lam}}R_i^{\a+1}T\fe|+|\dfrac{R_i^{\a+1}T\fe}{(1+t)^{\lam}}|\\
    &\les (1+t)|\sd R_i^{\a}T\fai_{\a}|+|\fai_i|\cdot|\lie_{R_i}^{\a}\chi'|,
    \end{split}
    \ee
    which implies that
    \bee
    \bes
    \dfrac{1}{1+t}||R_i^{\a+1}m||_{\normu2}&\les
    ||\sd R_i^{\a}T\fai||_{\normu2}+\dfrac{\da^{\frac{3}{2}}M}{(1+t)^2}
    ||\lie_{R_i}^{\a}\chi'||_{\normu2}\\
    &\les \dfrac{1}{1+t}\mu_m^{-\frac{1}{2}}\sqrt{\wi{E}_{1,|\a|+1}(t,u)}
    +\dfrac{\da^{\frac{3}{2}}M}{(1+t)^2}||\lie_{R_i}^{\a}\chi'||_{\normu2}.
    \end{split}
    \ee
    For the commutator, it follows from \eqref{estrgpi} that
    \bee
    \bes
    \dfrac{1}{1+t}||[L,R_i^{\a+1}]\mu||_{\normu2}&\les\dfrac{\da^{\frac{3}{2}}
    M(1+\ln(1+t))}{(1+t)^2}\dfrac{1}{1+t}||R_i^{\a}\mu||_{\normu2}\\
    &+\dfrac{\da^{\frac{1}{2}}M(1+\ln(1+t))}{1+t}||\lie_{R_i}^{\a}\chi'||
    _{\normu2}.
    \end{split}
    \ee
    Note that
    \bee
    L\left(\dfrac{1}{1+t}R_i^{\a+1}\mu\right)=-\dfrac{1}{(1+t)^2}R_i^{\a+1}\mu+\dfrac{1}{1+t}LR_i^{\a+1}\mu.
    \ee
    Then, Collecting the estimates above yields
    \bee
    \bes
    (1+t)^{-1}||R_i^{\a+1}\mu||_{\normu2}&\leq ||R_i^{\a+1}\mu||_{L^2(\Si_0^
    {\tilde{\da}})}+\int_0^t(1+t')^{-1}||LR_i^{\a+1}\mu||_{L^2(\Si_{t'}^{u})}dt'\\
    &\les ||R_i^{\a+1}\mu||_{L^2(\Si_0^{\tilde{\da}})}+\int_0^t\dfrac{\da^{\frac{3}{2}}M(1+\ln(1+t'))}{(1+t')^2}(1+t')^{-1}
    ||R_i^{\a}\mu||_{L^2(\Si_{t'}^{u})}\\
    &+\dfrac{1}{1+t'}\mu_m^{-\frac{1}{2}}\sqrt{\wi{E}_{1,|\a|+1}(t',u)}+
    \dfrac{\da^{\frac{1}{2}}M(1+\ln(1+t))}{1+t}||\lie_{R_i}^{\a}\chi'||_{L^2(\Si_{t'}^{u})}dt'.
    \end{split}
    \ee
    Applying the Gronwall inequality yields
    \bee\label{riamu1}
    \bes
    (1+t)^{-1}||R_i^{\a+1}\mu||_{\normu2}&\les||R_i^{\a+1}\mu||_{L^2(\Si_0^{\tilde{\da}})}
    +\int_0^t\dfrac{1}{1+t'}\mu_m^{-\frac{1}{2}}\sqrt{\wi{E}_{1,|\a|+1}(t',u)}dt'\\
    &+\int_0^t\dfrac{\da^{\frac{1}{2}}M(1+\ln(1+t'))}{1+t'}||\lie_{R_i}^{\a}\chi'||
    _{L^2(\Si_{t'}^{u})}dt'.
    \end{split}
    \ee
    To deal with the last integral in\eqref{riamu1}, define
    \bee
    \bes
    I_t(t')&=-\int_{t'}^t\dfrac{\da^{\frac{1}{2}}M(1+\ln(1+s))}{(1+s)^2}ds,\\
    J(t)&=\int_0^t\dfrac{\da^{\frac{3}{2}}M(1+\ln(1+t'))}{1+t'}\mu_m^{-\frac{1}{2}}
    \sqrt{\wi{E}_{1,|\a|+1}(t',u)}dt'.
    \end{split}
    \ee
    Then,
    \bee\label{Itt'}
    \bes
    &\int_0^t\dfrac{\da^{\frac{1}{2}}M(1+\ln(1+s))}{(1+s)^2}
    \int_0^s\dfrac{\da^{\frac{3}{2}}M(1+\ln(1+t'))}{1+t'}\mu_m^{-\frac{1}{2}}
    \sqrt{\wi{E}_{1,|\a|+1}(t',u)}dt'ds\\
    &=I_t(t)J(t)-I_t(0)J(0)-\int_0^tI_t(t')\dfrac{dJ(t')}{dt'}dt'.\\
    \end{split}
    \ee
    \hs Note that
    \bee
    -I_t(t')\leq\int_{t'}^{\infty}\dfrac{\da^{\frac{1}{2}}M(1+\ln(1+s))}{(1+s)^2}ds
    \leq\dfrac{1+\ln(1+t')}{1+t'}.
    \ee
    Then, it follows from Proposition\ref{Riachi'L2} and \eqref{Itt'} that
    \bee\label{riamu2}
    \bes
    &\int_0^t\dfrac{\da^{\frac{1}{2}}M(1+\ln(1+t'))}{1+t'}||\lie_{R_i}^{\a}\chi'||
    _{L^2(\Si_{t'}^{u})}dt'\\
    &\les\sum_{|\a'|\leq|\a|}||\lie_{R_i}^{\a'}\chi'||_{L^2(
    \Si_0^{\tilde{\da}})}+\int_0^t\dfrac{\da^{\frac{1}{2}}M(1+\ln(1+s))}{(1+s)^2}J(s)ds\\
    &\les\sum_{|\a'|\leq|\a|}||\lie_{R_i}^{\a}\chi'||_{L^2(\Si_0^{\tilde{\da}})}+\int_0^t\dfrac{\da^{\frac{3}{2}}M(1+\ln(1+t'))}{(1+t')^2}
    \mu_m^{-\frac{1}{2}}\sqrt{\wi{E}_{1,|\a|+1}(t',u)}dt'.
    \end{split}
    \ee
    \hs Combining the estimates \eqref{riamu1} and \eqref{riamu2} yields the desired estimates.
    \end{pf}
    \begin{remark}
    Similarly arguments and results hold for $R_i$ replaced by $Z_i\in
    \{Q, T,R_i\}$ in Proposition\ref{Riachi'L2} and \ref{Ria+1muL2}.
    \end{remark}
    To conclude, one has the following Proposition.
    \begin{prop}\label{estzichi'zimuL2}
    Under the bootstrap assumption, it holds that for all $1\leq|\a|\leq N_{top}-1$
    \begin{equation*}
    \bes
    ||\lie_{Z_i}^{\a}\chi'||_{\normu2}&\les(1+t)^{-1}\sum_{|\a'|\leq|\a|}||\lie_{Z_i}^{\a'}\chi'||_{
    L^2(\Si_0^{\tilde{\da}})}\\
    &+(1+t)^{-1}\left(\int_0^t\dfrac{\da^{\frac{3}{2}-l}M(1+\ln(1+t'))}{1+t'}
    \mu_m^{-\frac{1}{2}}\sqrt{\wi{E}_{1,|\a|+2}(t',u)}dt'\right),\\
    (1+t)^{-1}||Z_i^{\a+1}\mu||_{\normu2}&\les\sum_{|\a'|\leq|\a|}(||Z_i^{\a+1}\mu||_{L^2(\Si_0^{\tilde{\da}}
    )}+||\lie_{Z_i}^{\a}\chi'||_{L^2(\Si_0^{\tilde{\da}})})+\int_0^t\dfrac{\da^{-l}}{1+t'}\mu_m^{-\frac{1}{2}}\sqrt{\wi{E}_{1,|\a|+1}
    (t',u)}dt',
    \end{split}
    \end{equation*}
    where $l$ is the number of $T$'s in $Z_i^{\a}$. As a consequence, it holds that
    \begin{equation*}
    ||Z_i^{\a+1}y^j||_{\normu2}\les||\lie_{Z_i}^{\a}\chi'||_{\normu2},\quad ||Z_i^{\a+1}\lam_j||_{\normu2}\les(1+t)||\lie_{Z_i}^{\a}\chi'||_{
    \normu2},
    \end{equation*}
    \begin{equation*}
    \bes
    ||\lie_{Z_i}^{\a}\tpi_L||_{\normu2}&\les(1+t)^{-1}||Z_i^{\a+1}\mu||_{
    \normu2},\\
    ||\lie_{Z_i}^{\a}\tpi_{\dl}||_{\normu2}&\les\dfrac{\da^{\frac{1}{2}}
    M(1+\ln(1+t))}{1+t}||Z_i^{\a+1}\mu||_{\normu2},\\
    ||\lie_{Z_i}^{\a}\stpi||_{\normu2}&\les\da^{\frac{1}{2}}
    M(1+\ln(1+t))||\lie_{Z_i}^{\a}\chi'||_{\normu2},\\
    \end{split}
    \end{equation*}
    \begin{equation*}
    \bes
    ||\lie_{Z_i}^{\a}\rpi_L||_{\normu2}&\les||\lie_{Z_i}^{\a}\chi'||_{\normu2},\\
    ||\lie_{Z_i}^{\a}\rpi_{\dl}||_{\normu2}&\les\dfrac{\da^{\frac{1}{2}}
    M(1+\ln(1+t))}{1+t}||Z_i^{\a+1}\mu||_{\normu2}\\
    &+M(1+\ln(1+t))||\lie_{Z_i}^{\a}\chi'||_{\normu2},\\
    ||\lie_{Z_i}^{\a}\srpi||_{\normu2}&\les(1+\da^{\frac{3}{2}}
    M(1+\ln(1+t)))||\lie_{Z_i}^{\a}\chi'||_{\normu2}.\\
    \end{split}
    \end{equation*}
    \end{prop}
    \section{\textbf{Top order estimates for $\mu$ and tr$\chi$}}\label{section7}
\subsection{\textbf{Estimates for the top order angular derivatives of tr$\chi$}}
\hs Following the same framework in \cite{CZD1} section6, we regularize the transport equation for $tr\chi$ as follows.
\begin{lem}\label{squareh}
The enthalpy $h$ satisfies the following equation:
\bee\label{boxgh}
\bes
\square_gh&=\p_1+\mu^{-1}L\fai_i\dl\fai_i-\sdfai^2\\
&-\Omega^{-2}\dfrac{d\Omega}{dh}
\dfrac{1}{2\mu}\left[Lh\cdot\dl(h-\dfrac{a\fe}{(1+t)^{\lam}})+\dl h\cdot L(h-\dfrac{a\fe}{(1+t)^{\lam}})\right]\\
&+\Omega^{-2}\dfrac{d\Omega}{dh}\sd(h)\cdot\sd(h-\dfrac{a\fe}{(1+t)^{\lam}})
+a\square_g\dfrac{\fe}{(1+t)^{\lam}}:=\p,
\end{split}
\ee
\bee
\bes
\p_1&=-\dfrac{2\eta'}{\p\eta}\dfrac{a}{(1+t)^{\lam}}\dfrac{\pa\fai_j}{\pa x^j}[\fai_0-(\fai_i)^2-\frac{\lam}{1+t}\fe]\\
&+\dfrac{1}{\p^2}\dfrac{a}{(1+t)^{\lam}}\left(\dfrac{\pa}{\pa t}-\fai_j\dfrac{\pa\fai}{\pa x^j}\right)[\fai_0-\dfrac{1}{2}(\fai_i)^2-\frac{\lam}{1+t}\fe]\\
&-\dfrac{1}{\p^2}\dfrac{a\lam}{(1+t)^{\lam+1}}\left[\fai_0-\frac{1}{2}(\fai_i)^2-\frac{\lam}{1+t}\fe\right].
\end{split}
\ee
\hs Note that the right hand side of \eqref{boxgh} is of order $1$ and moreover, it holds that
\bee
\bes
|\mu\p_1|&\les\dfrac{1}{(1+t)^{\lam}}(|T\fai_j|+\mu|\sd\fai|)|[\fai_0-(\fai_i)^2-\frac{\lam}{1+t}\fe]|+
\dfrac{1}{(1+t)^{\lam}}|\dl[\fai_0-\dfrac{1}{2}(\fai_i)^2-\frac{\lam}{1+t}\fe]|\\
&+\dfrac{1}{(1+t)^{\lam+1}}|\fai_0-\dfrac{1}{2}(\fai_i)^2-\frac{\lam}{1+t}\fe|.
\end{split}
\ee
\end{lem}
\hs It follows from \eqref{decomposewavenonconformal} and Lemma\ref{squareh} that
\bee\label{musdeh}
\mu\s{\de}h=L(\dl h)+\dfrac{1}{2}(tr\dc\cdot Lh+tr\chi\cdot\dl h)+2\zeta\cdot\sd h+\mu\p.
\ee
\eqref{transporttrchi} implies
\bee\label{mutrx}
L(\mu tr\chi)=2(L\mu)tr\chi-\mu|\chi|^2-\mu tr\a.
\ee
\hs One can rewrite $\mu tr\a$ as
\bee\label{mutra}
\mu tr\a=-\dfrac{1}{2}\dfrac{dH}{dh}(L(\dl h))
-\dfrac{dH}{dh}\zeta\cdot\s{d}h-\dfrac{1}{2}\dfrac{dH}{dh}\mu\p+\mu tr\a^{[N]},
\ee
due to Lemma\ref{squareh} and \eqref{musdeh}. Then, substituting \eqref{mutra} into \eqref{mutrx} yields
\bee\label{mutrxf}
L\left(\mu tr\chi+\breve{f}\right)=2L\mu tr\chi-\mu|\chi|^2+\breve{g},
\ee
where
\bee
\breve{f}=-\dfrac{1}{2}\dfrac{dH}{dh}\dl h,\quad \breve{g}=-\dfrac{1}{2}\dfrac{d^2H}{dh^2}Lh\cdot\dl h+\dfrac{1}{2}\dfrac{dH}
{dh}\zeta\cdot\sd h+\dfrac{1}{2}\dfrac{dH}{dh}\mu\p-\mu tr\a^{[N]}.
\ee
\hs Next, set
$F_{\a}=\mu\sd R_i^{\a}tr\chi+\sd R_i^{\a}\breve{f}.$ When $|\a|=0$, it follows from \eqref{mutrxf} that
\begin{align}
LF_0&+(tr\chi-2\mu^{-1}L\mu)F_0=(\dfrac{1}{2}tr\chi-2\mu^{-1}L\mu)\sd\breve{f}
-\mu\sd|\hat{\chi}|^2+g_0,\label{F0}\\
g_0&=\sd\breve{g}+\dfrac{1}{2}tr\chi\cdot\s{d}(2L\mu+\breve{f})-(\sd\mu)
(L tr\chi+|\chi|^2).
\end{align}
Commuting $R_i^{\a}$ with \eqref{F0} yields
\bee\label{Fa}
LF_{\a}+(tr\chi-2\mu^{-1}L\mu)F_{\a}=(\dfrac{1}{2}tr\chi-2\mu^{-1}L\mu)\sd R_i^{\a}\breve{f}-\mu\sd R_i^{\a}|\hat{\chi}|^2+g_{\a},
\ee
where $g_{\a}$ is the commutation term given by
\bee\label{expressionga}
\bes
g_{\a}&=\lie_{R_i}^{\a}g_0+\sum_{|\be_1|+|\be_2|=|\a|,|\be_1|>0}
R_i^{\be_1}(\dfrac{1}{2}tr\chi-2\mu^{-1}L\mu)\sd R_i^{\be_2}\breve{f}
-\sum_{|\be_1|+|\be_2|=|\a|,|\be_1|>0}R_i^{\mu}\cdot\sd R_i^{\be_2}|\hat{\chi}|^2\\
&+\sum_{|\be_1|+|\be_2|=|\a|-1}\lie_{R_i}^{\be_1}\cdot
\lie_{\srpi_{L}}\cdot\lie_{R_i}^{\be_2}F_0-\sum_{|\be_1|+|\be_2|=|\a|,|\be_1|>0}
R_i^{\be_1}(tr\chi-2\mu^{-1}L\mu)\lie_{R_i}^{\be_{2}}F_0\\
&-\sum_{|\be_1|+|\be_2|=|\a|,|\be_1|>0}\lie_L(R_i^{\be_1}\mu\cdot\sd R_i^{\be_2}tr\chi)-(tr\chi-2\mu^{-1}L\mu)\sum_{|\be_1|+|\be_2|=|\a|,|\be_1|>0}R_i^{\be_1}\mu
\cdot \sd R_i^{\be_2}tr\chi.
\end{split}
\ee
\hs Since $|\a|=N_{top}-1$, then $F_{\a}$ is of order $N_{top}+1$. The top order terms (order of $N_{top}+1$) on the right hand side of \eqref{Fa} come from
\been
\item $\sd R_i^{\a}\breve{f}$, $\mu\sd R_i^{\a}|\hat{\chi}|^2$,
\item $g_{\a}$, which includes $\lie_{R_i}^{\be_1}\cdot
      \lie_{\srpi_{L}}\cdot\lie_{R_i}^{\be_2}F_0$, $R_i\mu\cdot\sd R_i^{\a-1}Ltr\chi$ and $\lie_{R_i}^{\a}g_0$. In $\lie_{R_i}^{\a}g_0$, the top order terms are:
      \bee\label{topordertermsRiag0}
       (\sd\mu)R_i^{\a}(Ltr\chi+|\chi|^2),\quad tr\chi\cdot \sd R_i^{\a}(2L\mu+\breve{f}),\quad \sd R_i^{\a}\breve{g}.
      \ee
\een
\hs Now we turn to the estimates for $F_{\a}$. For any $S_{t,u}$ $1-$form $\xi$, it holds that
\bee
|\xi|L|\xi|=(\xi,\lie_L\xi)-\xi\cdot\hat{\chi}\cdot\xi-\dfrac{1}{2}tr\chi|\xi|^2.
\ee
\hs Taking $\xi=F_{\a}$ yields
\bee\label{LFa}
L|F_{\a}|\leq(2\mu^{-1}L\mu-\dfrac{3}{2}tr\chi+|\hat{\chi}|)|F_{\a}|
+(\dfrac{1}{2}tr\chi+2\mu^{-1}|L\mu|)|\sd R_i^{\a}\breve{f}|+\mu|\sd R_i^{\a}|\hat{\chi}|^2|+|g_{\a}|.
\ee
\hs Define the integrating factor
\bee
B(t)=e^{-\int_0^t(2\mu^{-1}L\mu-\frac{3}{2}tr\chi+|\hat{\chi}|)dt'}=\left(\dfrac{\mu(t,u)}{\mu(0,u)}\right)^{-2}e^{\int_0^t\frac{3}{2}tr\chi}
\cdot e^{\int_0^t-|\hat{\chi}|}.
\ee
\hs Since $e^{\int_0^t\frac{3}{2}tr\chi dt'}=e^{\int_0^t\frac{3}{2}tr\chi'
+\frac{3}{1-u+t'}dt'}\leq C\left(\dfrac{1-u+t}{1-u}\right)^3$, then
\bee
B(t)\leq C\left(\dfrac{\mu(0,u)}{\mu(t,u)}\right)^2(1+t)^3,\quad \left[B(t)\right]^{-1}\leq C\left(\dfrac{\mu(t,u)}{\mu(0,u)}\right)^2\dfrac{1}
{(1+t)^3}.
\ee
Then, Multiplying $B(t)$ to \eqref{LFa} and applying the Gronwall inequality yield
\bee\label{Fal2}
\bes
||F_{\a}||_{\normu2}&\leq B(t)^{-1}||F_{\a}||_{L^2(\Si_{0}^{\tilde{\da}})}+B(t)^{-1}\int_0^tB(\tau)||
(\mu^{-1}(L\mu)_-+\dfrac{1}{2}tr\chi)\sd R_i^{\a}\breve{f}||_{L^2(\Si_{\tau}^{u})}\\
&+B(\tau)||\mu\sd R_i^{\a}|\hat{\chi}|^2||_{L^2(\Si_{\tau}^{u})}
+B(\tau)||g_{\a}||_{L^2(\Si_{\tau}^{u})}d\tau\\
&:=||F_{\a}||_{L^2(\Si_{0}^{\tilde{\da}})}+I_1+I_2+I_3,
\end{split}
\ee
where $I_1,I_2,I_3$ are the integrals given in order respectively, which will be estimated below.
\subsubsection{\textbf{Crucial lemma in the paper}}
\hs To deal with the integrals $I_1,I_2,I_3$ above, one needs the crucial lemma(Lemma\ref{crucial}) in this paper. Define:
    \beeq
    t_0&=&\inf\{t\in[0,t^{\ast})|\mu_m(t)<\dfrac{1}{10}\},\\
    M(t)&=&\max_{(u,\ta),(t,u,\ta)\in W_{shock}}|(\mu^{-1} L\mu)_-(t,u,\ta)|,\\
    I_{b}(t)&=&\int_{0}^t\mu^{-b}(\tau)M(\tau)d\tau.
    \eeq
    Then, the following crucial lemma and corollaries hold, which will be proved at the end of this subsection.
    \begin{lem}\label{crucial}
    For sufficiently large $b>2$ and all $t\in[t_0,t^{\ast})$, it holds that
    \bee
    I_b(t)\les b^{-1}\mu_m^{-b}(t).
    \ee
    \end{lem}
    \begin{cor}\label{coro1}
    For $\be>0$, sufficiently large $b$ and all $t\in[t_0,t^{\ast})$, it holds that
\[\int_{0}^t\dfrac{(1+\ln(1+t'))^{\be}}{1+t'}\mu_m^{-b-1}dt'\les
\dfrac{1}{b}(1+\ln(1+t))^{\be+1}\mu_m^{-b}(t).
\]
    \end{cor}
    \begin{cor}\label{coro2}
For sufficiently small $\da$ and fixed $b\ge2$, there exists a constant $C_0$ independent of $b$ and $\da$ such that for all $\tau\in[0,t]$, the following bound holds:
\[\mu^b_m(t)\leq C_0\mu_m^b(\tau).\]
\end{cor}
\hs Then, one could simplify the integrals $(I_1,I_2,I_3)$ as
\bee
\bes
I_1&\les\dfrac{1}{1+t}\int_0^t(1+t')(\dfrac{1}{2}|tr\chi|+2\mu^{-1}\max (L\mu)_-)||\sd R_i^{\a}\breve{f}||_{L^2(\Si_{\tau}^{u})}d\tau,\\
I_2&\les\int_0^t||\mu\sd R_i^{\a}|\hat{\chi}|^2||_{L^2(\Si_{\tau}^{u})}d\tau,\\
I_3&\les \int_0^t||g_{\a}||_{L^2(\Si_{\tau}^{u})}d\tau.
\end{split}
\ee
\subsubsection{\textbf{Elliptic estimates on $S_{t,u}$ and estimates for $I_1$, $I_2$ and $I_3$}}
\hs We start with estimating for $I_2$. To this end, rewrite \eqref{Codazzi2} as
\bee
\s{div}\hat{\chi}-\dfrac{1}{2}\sd tr\chi=\be^{\ast}-\mu^{-1}(\zeta\cdot\chi-
\zeta\cdot tr\chi):=i,
\ee
where $i$ is of order $1$. This indicates that one is able to estimate $\sd R_i^{\a}\hat{\chi}$ in terms of $\sd R_i^{\a}tr\chi$ together with some lower order terms.
\begin{lem}
Set $\pre {(i_k\cdots i_1)} {} {\hat{\chi}}:=\hat{\lie}_{R_{i_k}}\cdots\hat{\lie}_{R_{i_1}}\hat{\chi}$. Then
$\pre {(i_k\cdots i_1)} {} {\hat{\chi}}$ satisfies the following elliptic system:
\bee\label{sdivhatchi}
\s{div}\left(\pre {(i_k\cdots i_1)} {} {\hat{\chi}}\right)=
\dfrac{1}{2}\s{d}R_{i_k}\cdots R_{i_1}tr\chi+\pre {(i_k\cdots i_1)} {} {i},
\ee
where
\bee
\bes
\pre {(i_k\cdots i_1)} {} {i}&=\left(R_{i_k}+\dfrac{1}{2}tr\pre {R_{i_k}} {} {\pi}\right)\cdots\left(R_{i_1}+\dfrac{1}{2}tr\pre {R_{i_1}} {} {\pi}\right)i\\
&+\sum_{m=0}^{k-1}\left(R_{i_{k}}+\dfrac{1}{2}tr\pre {R_{i_{k}}} {} {\pi}\right)\cdots\left(R_{i_{k-m+1}}+\dfrac{1}{2}tr\pre {R_{i_{k-m+1}}} {}
{\pi}\right)\pre {(i_{k-m}\cdots i_1)} {} {q},\\
\pre {(i_{k+1}\cdots i_1)} {} {q_a}&=\dfrac{1}{4}tr\pre {R_{i_{k+1}}} {} {\pi}\cdot
\sd R_{i_k}\cdots R_{i_1}tr\chi\\
&+\dfrac{1}{2}\pre {R_{i_{k+1}}} {} {\hat{\pi}}^{BC}\left(D_B\pre {i_k\cdots i_1} {} {\hat{\chi}}_{AC}+D_C\pre {i_k\cdots i_1} {} {\hat{\chi}}_{AB}-D_A\pre {i_k\cdots i_1} {} {\hat{\chi}}_{BC}\right)\\
&+\left(\s{div}\pre {R_{i_{k+1}}} {} {\hat{\pi}}^B\cdot\pre {i_k\cdots i_1} {} {\hat{\chi}}_{AB}\right).
\end{split}
\ee
Note that $\pre {(i_k\cdots i_1)} {} {i}$ is of order $k+1$ and when $k=0$, $\pre {(i_k\cdots i_1)} {} {i}=i$.
\end{lem}
\begin{pf}
The lemma holds trivally for $k=0$. Assuming the lemma holds with $k$ replaced by $k-1$. Then, direct computations yield
\bee
\pre {(i_k\cdots i_1)} {} {i}=
R_{i_k}\pre {(i_{k-1}\cdots i_1)} {} {i}+\dfrac{1}{2}\pre {R_{i_k}} {} {\pi}
\pre {(i_k\cdots i_1)} {} {i}+\pre {(i_{k-1}\cdots i_{1})} {} {q}.
\ee
Applying Lemma\ref{recur} completes the proof.
\end{pf}
\begin{lem}
For sufficiently small $\da$ and any symmetric $S_{t,u}$ trace free $(0,2)$ tensor $\ta$, it holds that
\bee
\int_{S_{t,u}}\mu^2(|\nabla\ta|^2+|\ta|^2)d\mu_{\sg}\les\int_{S_{t,u}}\mu^2|\s{div}
\ta|^2+|\sd\mu|^2|\ta|^2d\mu_{\sg}.
\ee
\end{lem}
\begin{pf}
Let $J^A=\ta_C^B\s{D}_B\ta^{AC}-\ta^{AB}(\s{div}\ta)_B$. Obviously, $|J|\leq |\ta|(|\nabla\ta|+|\s{div}\ta|)$. 
Direct computations yield
\bee\label{bochner}
|\na\ta|^2+2K|\ta|^2=2|\s{div}\ta|^2+\s{div}J,
\ee
where $K$ is the Gauss curvature of $S_{t,u}$. It follows from \eqref{gaussta} that $K\sim 1$. Multiplying $\mu^2$ to \eqref{bochner}, then integrating over $S_{t,u}$ and applying the divergence theorem yeild
\bee
\bes
\int_{S_{t,u}}\mu^2\left(|\nabla\ta|^2+2K|\ta|^2\right)&=\int_{S_{t,u}}2\mu^2
|\s{div}\ta|^2+\mu^2\s{div}J=\int_{S_{t,u}}2\mu^2|\s{div}\ta|^2-2\mu J\cdot\sd\mu \\
&\leq\int_{S_{t,u}}2\mu^2|\s{div}\ta|^2+|\mu||\sd\mu||\ta|(|\nabla\ta|+|\s{div}
\ta|)\\
&\leq\int_{S_{t,u}}2\mu^2|\s{div}\ta|^2+\mu^2|\s{div}\ta|^2+\dfrac{1}{2}\mu^2|
\nabla\ta|^2+2|\sd\mu|^2|\ta|^2.
\end{split}
\ee
\end{pf}
Applying this lemma to \eqref{sdivhatchi} yields
\bee
\bes
||\mu\s{D}\pre {(i_k\cdots i_1)} {} {\hat{\chi}}||_{L^2(\Si_t^{u})}&\les
||\sd\mu||_{\supnormda}||\pre {(i_k\cdots i_1)} {} {\hat{\chi}}||_{\normu2}
+||\mu\s{div}\pre {(i_k\cdots i_1)} {} {\hat{\chi}}||_{\normu2}\\
&\les ||\mu \sd R_{i_k}\cdots R_{i_1}tr\chi||_{\normu2}+ (||\sd\mu||_{\supnormda}+||tr\srpi||_{\supnormda})||\pre {(i_k\cdots i_1)} {} {\hat{\chi}}||_{\normu2}\\
&+||\mu R_{i_k}\cdots R_{i_1}\fai_{\a}||_{\normu2}.
\end{split}
\ee
For $I_2$, since
\bee\label{sdriahatchi}
\mu\sd R_i^{\a}|\hat{\chi}|^2=\sum\mu\sd\left(\hat{\lie_{R_i}^{\be_1}}\hat{\chi}\cdot\hat{\lie_{R_i}
^{\be_2}}\hat{\chi}\cdot\hat{\lie_{R_i}^{\be_3}}\sg^{-1}\cdot\hat{\lie_{R_i}
^{\be_4}}\sg^{-1}\right)=2\mu\sd\hat{\lie_{R_i}^{\a}}\hat{\chi}\cdot\hat{\chi}+\text{l.o.ts},
\ee
then
\begin{equation*}
\bes
||\mu\sd R_i^{\a}|\hat{\chi}|^2||_{\normu2}&\les |\hat{\chi}|\cdot||\mu\sd R_i^{\a}tr\chi||_{\normu2}+|\hat{\chi}|\cdot(1+t)\sqrt{\wi{E}_{0,\leq|\a|+2}}\\
&+(|\sd\mu|+|tr\srpi|)|\hat{\chi}|\cdot||\pre {(i_k\cdots i_1)} {} {\hat{\chi}}||_{\normu2}+\text{l.o.ts}\\
&\les\dfrac{\da^{\frac{3}{2}}M(1+\ln(1+t))}{(1+t)^2}\left(||F_{\a}||_{\normu2}
+||\sd R_i^{\a}\breve{f}||_{\normu2}\right)+\dfrac{\da^{\frac{3}{2}}M(1+\ln(1+t))}{1+t}\sqrt{\wi{E}_{0,\leq|\a|+2}}\\
&+\dfrac{\da^2M^2(1+\ln(1+t))^2}{(1+t)^4}\left(||\lie_{R_i}^{\a}\chi'||_{
\normu2}+\int_0^t\dfrac{\da^{\frac{3}{2}}M(1+\ln(1+\tau))}{1+\tau}\mu_m^{-\frac{1}{2}}
\sqrt{\wi{E}_{1,\leq|\a|+1}}d\tau\right).
\end{split}
\end{equation*}
It follows from the definition of $\breve{f}$ that
\bee
\sd R_i^{\a}\breve{f}=\sd R_i^{\a}\left(-\dfrac{1}{2}\dfrac{dH}{dh}\dl h\right)=-\dfrac{1}{2}\dfrac{dH}{dh}\dfrac{1}{1+t}\left( \dl R_i^{\a+1}\fai_0+\fai_i\dl R_i^{\a+1}\fai_i\right)+\text{l.o.ts},
\ee
and then
\bee
||\sd R_i^{\a}\breve{f}||_{L^2(\Si_{t}^{u})}\les\dfrac{1}{1+t}
\sqrt{E_{0,\leq |\a|+2}(t,u)}.
\ee
Hence,
\bee
\bes
I_2&\les \int_0^t\dfrac{\da^{\frac{3}{2}}M(1+\ln(1+t'))}{(1+t')^2}||F_{\a}||_{
L^2(\Si_{t'}^{u})}dt'+\da^2M^2\sum_{|\a'|\leq|\a|}||\lie_{R_i}^{\a'}\chi'||_{L^2(\Si_0^{\tilde{\da}})}\\
&+\int_0^t\dfrac{\da^{\frac{3}{2}}M(1+\ln(1+t'))}{1+t'}\sqrt{\wi{E}_{0,\leq|\a|+2}
}dt'+\int_0^t\dfrac{\da^{\frac{3}{2}}M(1+\ln(1+t'))}{1+t'}\mu_m^{-\frac{1}{2}}
\sqrt{\wi{E}_{1,\leq|\a|+1}}dt',
\end{split}
\ee
\hs For $I_3$, it follows from the discussion \eqref{expressionga}-\eqref{topordertermsRiag0} that one needs only to consider the following cases:
\been
\item $||\lie_{R_i}^{\be_1}\lie_{\rpi}F_{\be_2}||_{\normu2}\leq \dfrac{\da^{
\frac{1}{2}}M(1+\ln(1+t))}{(1+t)^2}||F_{\a}||_{\normu2}$,
\item Since $Ltr\chi+|\chi|^2=e tr\chi-tr\a'+a\eta^{-1}\frac{\hat{T}^i\fai_i}{(1+t)^{\lam}}tr\chi$, it holds that
      \begin{equation*}
      \bes
      &|\sd\mu|||R_i^{\a}(Ltr\chi+|\chi|^2)||_{\normu2}
      \leq|\sd\mu|\left(|e|\cdot||R_i^{\a}tr\chi||_{\normu2}+|tr\chi|\cdot
      ||R_i^{\a}e||_{\normu2}+||\sd^2 R_i^{\a}\fai||_{\normu2}\right)\\
      &\les\dfrac{\da^2 M^2(1+\ln(1+t))}{(1+t)^4}\sum_{|\a'|\leq|\a|}||\lie_{R_i}^{\a'}\chi'||_{
      L^2(\Si_0^{\tilde{\da}})}+\dfrac{\da^{\frac{1}{2}}M(1+\ln(1+t))}{(1+t)^2}\mu_m^{-\frac{1}{2}}
      \sqrt{\wi{E}_{1,\leq|\a|+2}}\\
      &+\dfrac{\da^2 M^2(1+\ln(1+t))}{(1+t)^2}\int_0^t
      \dfrac{1+\ln(1+t')}{1+t'}\mu_m^{-\frac{1}{2}}\sqrt{\wi{E}_{1,\leq|\a|+2}
      (t',u)}dt'.\\
      \end{split}
      \end{equation*}
\item It follows from the definition of $\breve{f}$ that $\sd R_i^{\a}(2L\mu+\breve{f})=\sd R_i^{\a}
(\mu e+\dfrac{1}{2}\dfrac{dH}{dh}\eta^{-2}\mu Lh)$ and then
     \begin{equation*}
     |tr\chi|\cdot||\sd R_i^{\a}(2L\mu+\breve{f})||_{\normu2}\les
     \dfrac{\da^{\frac{3}{2}}M}{(1+t)^3}\sum_{|\a'|\leq|\a|}\left(||R_i^{\a+1}\mu||
     _{L^2(\Si_0^{\tilde{\da}})}+||\lie_{R_i}^{\a}\chi'||_{L^2(\Si_0^{\tilde{\da}})}\right)
     +\dfrac{1}{1+t}\mu_m^{-\frac{1}{2}}\sqrt{\wi{E}_{1,\leq|\a|+2}}.
     \end{equation*}
\item\footnote{The major contribution in $\sd R_i^{\a}\breve{g}$ comes from $\sd R_i^{\a}(Lh\cdot\dl h)$.} 
      \begin{equation*}
      \bes
      \sd R_i^{\a}(Lh\cdot\dl h)&=\dfrac{1}{1+t}(\dl h\cdot LR_i^{\a+1}h+Lh\cdot\dl R_i^{\a}h)+\text{l.o.ts},\\
      ||\sd R_i^{\a}\breve{g}||_{\normda2}&\les\dfrac{\da^{\frac{1}{2}}M}{(1+t)^2}
      \left(\int_0^u\wi{F}_{1,\leq|\a|+2}(t,u')du'\right)^{\frac{1}{2}}+\int_0^t\dfrac{\da^{\frac{3}{2}M}}{(1+t')^3}\sqrt{\wi{E}_{0,\leq|\a|+2}
      (t',u)}dt'.
      \end{split}
      \end{equation*}
\een
\hs Collecting the results above yields
     \bee
     \bes
     I_3&\les \int_0^t\dfrac{\da^{\frac{1}{2}}M(1+\ln(1+t'))}{(1+t')^2}||F_{\a}
     ||_{L^2(\Si_{t'}^{u})}dt'+\dfrac{\da^{\frac{3}{2}}M}{(1+t)^3}\left(
     ||R_i^{\a+1}\mu||_{L^2(\Si_0^{\tilde{\da}})}+||\lie_{R_i}^{\a}\chi'||_{L^2
     (\Si_0^{\tilde{\da}})}\right)\\
     &+\dfrac{\da^{\frac{3}{2}}M(1+\ln(1+t))}{(1+t)^2}\mu_m^{-\frac{1}{2}}\sqrt{
     \wi{E}_{1,\leq|\a|+2}(t,u)}\\
     &+\int_0^t\dfrac{\da^{\frac{3}{2}}M}{(1+t')^3}\sqrt{\wi{E}_{0,\leq|\a|+2}
      (t',u)}dt'+\int_0^t\dfrac{\da^2M^2(1+\ln(1+t'))}{1+t'}
     \mu_m^{-\frac{1}{2}}\sqrt{\wi{E}_{1,\leq|\a|+2}(t',u)}dt'\\
     &+\dfrac{\da^{\frac{1}{2}}M}{(1+t)^2}
      \left(\int_0^u\wi{F}_{1,\leq|\a|+2}(t,u')du'\right)^{\frac{1}{2}}.
     \end{split}
     \ee
\hs Applying Lemma\ref{crucial} to $I_1$ and taking $b=\ba2$ yield
    \bee
    \bes
    I_1&\les \int_{0}^t(\mu^{-1}|L\mu|+\dfrac{1}{1+t})||\sd R_i^{\a}\breve{f}||_{L^2(\Si_{t'}^{u})}dt'\\
    &\les\int_0^t\dfrac{1}{(1+t')^2}\sqrt{\wi{E}_{0,\leq|\a|+2}(t',u)}dt'
    +\int_0^t2\mu^{-1}|L\mu|\dfrac{1}{1+t}\sqrt{\wi{E}_{0,\leq|\a|+2}(t',u)}
    dt'\\
    &\les\mu_m^{-\ba2+1}(1+\ln(1+t))^p\sqrt{\ol{E}_{0,\leq|\a|+2}}.
    \end{split}
    \ee
    Collecting the estimates for $I_1$, $I_2$ and $I_3$ and applying the Gronwall inequality yield
    \bee\label{boundssdRiatrchi}
    \bes
    ||F_{\a}||_{\normu2}&\les\dfrac{\da^{\frac{3}{2}}M}{(1+t)^3}
    \sum_{|\a'|\leq|\a|}\left(||R_i^{\a'+1}\mu||_{L^2(\Si_0^{\tilde{\da}})}+||\lie_{R_i}^{\a'}\chi'||
    _{L^2(\Si_0^{\tilde{\da}})}\right)\\
    &+\dfrac{1}{1+t}\int_0^t\da^{\frac{3}{2}}M(1+\ln(1+t'))\sqrt{\wi{E}_{0,\leq|\a|
    +2}}dt'\\
    &+\dfrac{1}{1+t}\int_0^t\dfrac{\da^{\frac{1}{2}}M(1+\ln(1+t'))}{1+t'}
    \mu_m^{-\frac{1}{2}}\sqrt{\wi{E}_{1,\leq|\a|+2}}dt'\\
    &+\dfrac{C}{\ba2}\mu_m^{-\ba2}(1+\ln(1+t))^p\sqrt{\ol{E}_{0,\leq|\a|+2}}+\dfrac{\da^{\frac{1}{2}}M}{(1+t)^2}
    \left(\int_0^t\wi{F}_{1,\leq|\a|+2}du'\right)^{\frac{1}{2}}.
    \end{split}
    \ee
    \hs One can obtain the same bounds for $||\mu\sd R_i^{\a}tr\chi||_{\normu2}$ due to $||\mu\sd R_i^{\a}tr\chi||_{\normu2}\leq||\sd R_i^{\a}\breve{f}||_{\normu2}+||F_{\a}||_{\normu2}$. It remains to prove the Lemma\ref{crucial} and the key step is to bound $\mu$ accurately from both below and above.
    \begin{pf}\footnote{Here we give the proof for Lemma\ref{crucial}. The proof for Corollary\ref{coro1} is the same as Lemma\ref{crucial} and the proof for Corollary\ref{coro2} can be found in\cite{CZD1} Corollary 6.2.}
    Note that in the shock region, $L\mu<0$. Thus, denote $-\eta_m(\tau):=\min_{(u,\ta)}L\mu(\tau,u,\ta)$ for $\tau\in[t_0,t^{\ast})$ where $\eta_m(\tau)>0$, which can be achieved at some $(u_m,\ta_m)$. Moreover, for any $t$, one can choose $(u_t,\ta_t)$ such that $\mu_m(t)=\mu(t,u_t,\ta_t)$. \textbf{Fix} $s\in [t_0,t]$. It follows from
    \bee
    \bes
    &(1+t)A(t)L\mu(t,u_s,\ta_s)=(1+s)A(s) L\mu(s,\mu_s,\ta_s)+I\\
    &=-(1+s)A(s)\eta_m(s)+\left[(1+s)A(s)L\mu(s,u_s,\ta_s)-(1+s)A(s)L\mu(s,u_m,\ta_m)\right]+I\\
    &:=-(1+s)A(s)\eta_m(s)+II,
    \end{split}
    \ee
    and Proposition\ref{accruatemu1} that $I$ can be bounded by $O(\da^{\frac{3}{2}} M)$. Furthermore, it follows from the argument of Proposition\ref{accruatemu1} that
    \begin{equation*}
    \bes
    1+O(\da)&=\mu(0,u_s,\ta_s)=\mu(s,u_s,\ta_s)+\int_s^0\frac{1}{(1+\tau)A(\tau)}d\tau[(1+s)A(s)L\mu(s,u_s,\ta_s)+O(\da^{\frac{3}{2}} M)],\\
    1+O(\da)&=\mu(0,u_m,\ta_m)=\mu(s,u_m,\ta_m)+\int_s^0\frac{1}{(1+\tau)A(\tau)}d\tau
    [(1+s)A(s)L\mu(s,u_m,\ta_m)+O(\da^{\frac{3}{2}} M)].
    \end{split}
    \end{equation*}
    Thus, this shows that $\left|(1+s)A(s)(L\mu(s,u_s,\ta_s)-L\mu(s,u_m,\ta_m))\right|\les\da^{\frac{3}{2}} M$, which implies $|II|\les\da^{\frac{3}{2}} M$. Therefore,
    \bee
    \bes
    \mu_m(t)&=\mu_m(s,u_t,\ta_t)+\int_s^tL\mu(\tau,u_t,\ta_t)d\tau\\
    &\geq \mu_m(s)+\int_s^t\dfrac{1}{(1+\tau)A(\tau)}d\tau((1+s)A(s)L\mu(s,u_t,\ta_t)+I)\\
    &\geq\mu_m(s)+\int_s^t\frac{1}{(1+\tau)A(\tau)}d\tau\left(-(1+s)A(s)\eta_m(s)-\dfrac{1}{b}
    \right),
    \end{split}
    \ee
    provided that $\da^{\frac{3}{2}} M<\dfrac{1}{b}$. On the other hand, it holds that
    \bee\label{upper}
    \bes
    \mu_m(t)&\leq \mu_m(s,u_s,\ta_s)+\int_s^t\frac{(1+\tau)A(\tau)}{(1+\tau)A(\tau)}L\mu(\tau,u_s,\ta_s)d\tau\\
    &=\mu_m(s)+\int_s^t\dfrac{1}{(1+\tau)A(\tau)}d\tau(-(1+s)A(s)\eta_m(s)+II)\\
    &\leq \mu_m(s)+\left(-(1+s)A(s)\eta_m(s)+\dfrac{1}{b}\right)\int_s^t\dfrac{1}{(1+\tau)A(\tau)}d\tau.
    \end{split}
    \ee
    Hence, it follows from Proposition\ref{keymu1} that\footnote{$x_1=\int_0^t\frac{1}{(1+y)A(y)}dy$.}
    \begin{equation*}
    \bes
    I_b(t)&\leq \int_0^t\left[\mu_m(s)+\int_s^{\tau}\dfrac{1}{(1+y)A(y)}dy
    \left(-(1+s)A(s)\eta_m(s)-\dfrac{1}{b}\right)\right]^{-b-1}\dfrac{1}{(1+\tau)A(\tau)}\left(\int_0^{\tau}\frac{1}{(1+y)A(y)}\right)^{-1}
    d\tau\\
    &\leq\int_0^{x_1}\left[\mu_m(s)+[x-\int_0^s\dfrac{1}{(1+y)A(y)}dy]
    \left(-(1+s)A(s)\eta_m(s)-\dfrac{1}{b}\right)\right]^{-b-1}dx\\
    &\les\dfrac{1}{b}\left[\mu_m(s)+\int_s^t\frac{1}{(1+\tau)A(\tau)}d\tau
    \left(-(1+s)A(s)\eta_m(s)-\dfrac{1}{b}\right)\right]^{-b}.
    \end{split}
    \end{equation*}
    Then,
    \bee
    \bes
    I_b(t)&\leq \dfrac{1}{b}\left[\dfrac{\mu_m(s)+
    \left(-(1+s)A(s)\eta_m(s)-\dfrac{1}{b}\right)\int_s^t\frac{1}{(1+\tau)A(\tau)}d\tau}{\mu_m(s)+
    \left(-(1+s)A(s)\eta_m(s)+\dfrac{1}{b}\right)\int_s^t\frac{1}{(1+\tau)A(\tau)}d\tau}\right]
    ^{-b}\mu_m^{-b}(t)\\
    &\les\dfrac{1}{b}\left(\dfrac{-(1+s)A(s)\eta_m(s)-\frac{1}{b}}{-(1+s)A(s)\eta_m(s)+\frac{1}
    {b}}\right)^{-b}\mu_m^{-b}(t)\les \dfrac{1}{b}\mu_m^{-b}(t).
    \end{split}
    \ee
    \end{pf}
    \subsection{\textbf{Estimates for the top order spatial derivatives of $\mu$}}
    \hs In this subsection, we derive the estimates for $R_i^{\a'}T^l\s{\de}\mu$ by making use of the basic transport equation \eqref{transportmu}. 
    First, commuting $\s{\de}$ with \eqref{transportmu} yields:
    \bee\label{lmusdemu}
    L(\mu\s{\de}\mu)=(m+2\mu e)\s{\de}\mu+\mu[L,\s{\de}]\mu+\mu\s{\de}m+\mu^2\s{\de}e.
    \ee
    \hs Note that
    \begin{equation*}
    \mu\s{\de}m=L\left(\breve{f}_0\right)+\dfrac{1}{2}tr\chi\breve{f}_0+\dfrac{1}{2}\dfrac{dH}{dh}\eta^{-2}\mu^2\sd tr\chi\cdot\sd h+m\s{\de}\mu+\frac{a\mu^2}{(1+t)^{\lam}\eta^2}(\sd x^i)\fai_i\cdot\sd tr\chi+m_0,
    \end{equation*}
    where $\breve{f}_0=\dfrac{1}{2}\dfrac{dH}{dh}T\dl h$ and $m_0$ doesn't contain any acoustical terms of order 2. Similarly,
    \begin{equation*}
    \mu^2\s{\de}\mu=L(\mu\breve{f}_1)+\dfrac{1}{2}tr\chi\breve{f}_1+\mu\left[\dfrac{1}{\eta}\dfrac{d\eta}{dh}\sd h+
    \dfrac{1}{\eta}\hat{T}^i\sd\fai_i-\dfrac{1}{\eta^2}(\sd h-a\dfrac{T\fe}{(1+t)^{\lam}}-\eta\hat{T}^i\sd\fai_i)\right]\mu\sd tr\chi+e_0,
    \end{equation*}
    where $\breve{f}_1=\dfrac{1}{\eta}\dfrac{d\eta}{dh}L\dl h+\dfrac{1}{\eta}\hat{T}^iL\dl\fai_i$ and $e_0$ doesn't contain any acoustical terms of order 2. Set $\breve{f}=\breve{f}_0+\mu\breve{f}_1$. Then, 
    \bee\label{lf0mu}
    L(\mu\s{\de}\mu)=(2\mu^{-1}L\mu-tr\chi)\mu\s{\de}\mu+L\breve{f}+\dfrac{1}{2}
    tr\chi\cdot\breve{f}-2\mu\hat{\chi}\cdot\hat{\s{D}}^2\mu+\breve{g},
    \ee
    where the principle acoustical part of $\breve{g}$ is given by
    \bee
    \left[\breve{g}\right]_{P.A.}=
    \left(-\sd\mu+\dfrac{2\mu}{\eta}\hat{T}^i\sd\fai_i-\dfrac{2\mu}{\eta^2}\sd h+\dfrac{a\mu}{(1+t)^{\lam}\eta^2}(T\fe+(\sd x^i)\fai_i)\right)\mu\sd tr\chi.
    \ee
     \hs Define $F_{\a',l}=\mu R_i^{\a'}T^l\s{\de}\mu-R_i^{\a'}T^l\breve{f}$. Then, $F_{\a',l}$ satisfies the following transport equation:
    \bee\label{lfal}
    \bes
    LF_{\a',l}&+(tr\chi-2\mu^{-1}L\mu)F_{\a',l}=(-\dfrac{1}{2}tr\chi+2\mu^{-1}L\mu)
    R_i^{\a'}T^l\breve{f}-2\mu\hat{\chi}\cdot R_i^{\a}T^l(\hat{\s{D}}^2\mu)+g_{\a',l},\\
    g_{\a',l}&=R_i^{\a'}g_{0,l}+\sum_{k=0}^{\a'-1}R_i^{\a'}\srpi_LF_{\a'-k-1,l}
    +\sum_{k=0}^{\a'-1}R_i^{\a'}y_{\a'-k-1,l},
    \end{split}
    \ee
    where $y_{\a',l}$ is given by
    \bee
    \bes
    y_{\a',l}&=\left[2(R_iL\mu)-\mu R_itr\chi\right]R_i^{\a}T^l\s{\de}\mu+
    \dfrac{1}{2}R_i(tr\chi)R_i^{\a}T^l\breve{f}\\
    &-2R_i(\mu\hat{\chi})\cdot R_i^{\a'}T^l(\hat{\s{D}}^2\mu)-4tr\rpi\cdot
    \mu\hat{\chi}\cdot R_i^{\a'}T^l(\hat{\s{D}}^2\mu).
    \end{split}
    \ee
    By the expression of $g_{\a',l}$, 
    the top order acoustical terms in $g_{\a',l}$ come from: $R_i^{\a'}g_{0,l}$ and $R_i^k\srpi_L F_{\a'-k-1,l}$, where the latter one is equalient to $R_i^{\a'+1}T^{l-1}\s{\de}\mu$. In $R_i^{\a'}g_{0,l}$, the top order acoustical terms come from:
    \been
    \item  $R_i^{\a'}T^k\Lambda F_{0,l-k-1}$, which is equivalent to $\mu R_i^{\a'+1}T^{l-1}\s{\de}\mu$.
    \item $R_i^{\a'}T^l\breve{g}$, which is equivalent to $\mu R_i^{\a'}T^l\sd
    tr\chi$ and the following two cases hold:
          \bee\label{discussionga'l}
          \mu\sd R_i^{\a}tr\chi,\hs \text{if}\ l=0;\quad \mu \sd R_i^{\a'}T^{l-1}\s{\de}\mu,\hs \text{if}\ l\geq 1.
          \ee

    \een
    It follows from \eqref{lfal} that
    \begin{equation*}
    L|F_{\a',l}|\leq (2\mu^{-1}L\mu-\dfrac{3}{2}tr\chi+|\hat{\chi}|)|F_{\a',l}|+(\dfrac{1}{2}
    tr\chi-2\mu^{-1}L\mu)|R_i^{\a'}T^l\breve{f}|+2\mu|\hat{\chi}|\cdot|R_i^{\a'}T^l(\hat{\s{D}}^2\mu)|+|g_{\a',l}|.
    \end{equation*}
    Multiplying $B(t)$ on the both sides and applying the Gronwall inequality yield
    \bee
    \bes
    ||F_{\a',l}||_{\normu2}&\leq B(t)^{-1}||F_{\a',l}||_{L^2(\Si_0^{\tilde{\da}})}+B(t)^{-1}\int_{0}^tB(\tau)\left(\dfrac{1}{2}tr\chi+
    2\mu^{-1}|L\mu|\right)||R_i^{\a'}T^l\breve{f}||_{L^2(\Si_{\tau}^{u})}\\
    &+B(\tau)\cdot2\mu|\hat{\chi}|\cdot||R_i^{\a'}T^l\hat{\s{D}}^2\mu||_{L^2(\Si_{\tau}^{u})
    }+B(\tau)||g_{\a',l}||_{L^2(\Si_{\tau}^{u})}d\tau\\
    &=||F_{\a',l}||_{L^2(\Si_{0}^{\tilde{\da}})}+I_1+I_2+I_3,
    \end{split}
    \ee
    where $I_1,I_2,I_3$ are the integrals given in order respectively, which will be estimated below.
    \subsubsection{\textbf{Elliptic estimates for $\mu$ on $S_{t,u}$ and estimates for $I_1$, $I_2$, $I_3$}}
    \hs Note that $I_2$ involves $\hat{\s{D}}^2\mu=\s{D}^2\mu-\dfrac{1}{2}(\s{\de}\mu) g$. The following lemma implies that one could bound $\s{D}^2\mu$ by $\s{\de}\mu$.
    \begin{lem}\label{ellipticmu}
    For any $S_{t,u}$ function $\fe$, it holds that
    \bee
    \int_{S_{t,u}}\mu^2\left(\dfrac{1}{2}|\nabla^2\fe|^2+K|\sd\fe|^2\right)d
    \mu_{\sg}\leq 2\int_{S_{t,u}}\mu^2|\s{\de}\fe|^2d\mu_{\sg}+3\int_{S_{t,u}}
    |\sd\mu|^2\cdot|\sd\fe|^2d\mu_{\sg},
    \ee
    where $K$ is the Gauss curvature of $S_{t,u}$.
    \end{lem}
    \begin{pf}
    Define $1-$form $\xi=\sd\fe$. Then,
    $\nabla_a\nabla^b\xi_b-\nabla^b\nabla_a\xi_b=-K\xi_a$, i.e.
    \bee\label{nablaasdefe}
    \nabla_a(\s{\de}\fe)=\nabla^b(\nabla_a\xi_b)-K\xi_a,
    \ee
    where $K$ is the Gauss curvature of $S_{t,u}$. Multiplying $-\mu^2\xi$ to \eqref{nablaasdefe}, then integrating over $S_{t,u}$ and using the divergence theorem yield
    \bee
    \bes
    \int_{S_{t,u}}-\mu^2\xi\nabla_a(\s{\de}\fe)&=\int_{S_{t,u}}\left(
    2\mu\sd\mu\cdot\sd\fe+\mu^2\s{\de}\fe\right)\s{\de}\fe,\\
    \int_{S_{t,u}}-\mu^2\xi\nabla^b(\nabla_a\xi_b)&=\int_{S_{t,u}}\left(
    2\mu\sd\mu\sd\fe+\mu^2\nabla^2\fe\right)\cdot\nabla^2\fe.
    \end{split}
    \ee
    Hence,
    \bee
    \bes
    &\int_{S_{t,u}}\mu^2\left(\dfrac{1}{2}|\nabla^2\fe|^2+K|\sd\fe|^2\right)d
    \mu_{\sg}=\int_{S_{t,u}}\left(
    2\mu\sd\mu\cdot\sd\fe+\mu^2\s{\de}\fe\right)\s{\de}\fe-\left(
    2\mu\sd\mu\sd\fe\right)\cdot\nabla^2\fe d\mu_{\sg}\\
    &\leq2\int_{S_{t,u}}\mu^2|\s{\de}\fe|^2+3
    |\sd\mu|^2\cdot|\sd\fe|^2d\mu_{\sg}-\dfrac{1}{2}\mu^2|\nabla^2\fe|^2d\mu_{\sg}.
    \end{split}
    \ee
    \end{pf}
    \hs Applying this lemma to $I_2$ yields
    \bee\label{riatld2mu}
    \bes
    ||\mu R_i^{\a'}T^l(\hat{\s{D}}^2\mu)||_{\normu2}&\leq
    ||\mu R_i^{\a'}T^l\s{\de}\mu||_{\normu2}+||\mu\sd\mu\sd R_i^{\a'}T^l\mu||_{\normu2}\\
    &\leq ||F_{\a',l}||_{\normu2}+||R_i^{\a'}T^l\breve{f}||_{\normu2}
    +\dfrac{1}{(1+t)}||R_i^{\a'+1}T^l\mu||_{L^2(\Si_{0}^{\tilde{\da}})}\\
    &+\dfrac{1}{1+t}\left(\da^{-l}||\lie_{R_i}^{\a'+1}\lie_{T}^l\chi'||_{L^2(\Si_
    0^{\da})}+\int_0^t\da^{-l}\dfrac{1}{1+t'}\mu_m^{-\frac{1}{2}}\sqrt{
    \wi{E}_{1,\leq|\a|+1}}dt'\right).
    \end{split}
    \ee
    \hs It follows from the definition of $\breve{f}$ that
    \bee
    ||R_i^{\a'}T^l\breve{f}||_{\normu2}\les\da^{-(l+1)}
    \sqrt{\wi{E}_{0,\leq|\a|+2}(t,u)}+\dfrac{\da^{\frac{3}{2}}M}{(1+t)^2}||R_i^{\a'}T^l
    \mu||_{\normu2}.
    \ee
    Therefore,
    \bee
    \bes
    \da^{(l+1)}I_2&\leq\int_0^t\dfrac{\da^{\frac{3}{2}}M(1+\ln(1+t'))}{(1+t')^2}
    \da^{l+1}||F_{\a',l}||_{L^2(\Si_{t'}^{u})}dt'\\
    &+\da^{l+2}\left(||R_i^{\a'+1}T^l\mu||_{L^2(\Si_0^{\tilde{\da}})}+\lVert\lie_{R_i}^{\a'}
    \lie_{T}^l\chi'\rVert_{L^2(\Si_0^{\tilde{\da}})}\right)\\
    &+\int_0^t\dfrac{\da^{\frac{3}{2}}M(1+\ln(1+t'))}{(1+t')^2}\sqrt{\wi{E}_{0,\leq
    |\a|+2}(t',u)}+\int_0^t\dfrac{\da^{\frac{3}{2}}M}{(1+t')}\mu_m^{-\frac{1}{2}}
    \sqrt{\wi{E}_{1,\leq|\a|+1}(t',u)}dt'.
    \end{split}
    \ee
    \hs For $I_3$, it follows from the discussion around \eqref{discussionga'l} and \eqref{boundssdRiatrchi} that
    \begin{equation*}
    \bes
    \da^{(l+1)}I_3&
    \les\dfrac{\da^{\frac{3}{2}}M(1+\ln(1+t'))}{1+t'}
    \|\mu\sd R_i^{\a}tr\chi\|_{L^2(\Si_{t'}^{u})}dt'\\
    &+\int_0^t\dfrac{\da^{\frac{1}{2}}M(1+\ln(1+t'))}{(1+t')^2}\da^{l+1}\|F_{\a',l}
    \|_{L^2(\Si_{t'}^{u})}dt'+\int_0^t\dfrac{\da^{\frac{1}{2}}M(1+\ln(1+t'))}
    {(1+t')^2}\sqrt{\wi{E}_{0,\leq|\a|+2}(t',u)}dt'.
    \end{split}
    \end{equation*}
    \hs For $I_1$, it follows from Lemma\ref{crucial} that
    \bee
    \bes
    I_1&\leq\int_{0}^t(tr\chi+2\mu^{-1}|L\mu|)||R_i^{\a'}T^l\breve{f}||_{L^2(\Si_{
    t'}^{u})}dt'\les\int_0^t\dfrac{1}{1+t'}\|R_i^{\a'}T^l\breve{f}\|_{L^2(\Si_{t'}^{u})}dt'\\
    &+\int_0^t\mu^{-1}|L\mu|\|R_i^{\a'}T^l\breve{f}\|_{L^2(\Si_{t'}^{u})}dt'\\
    &\les\int_0^t\dfrac{\da^{\frac{1}{2}}M(1+\ln(1+t'))}{1+t'}\sqrt{\wi{E}_
    {0,\leq|\a|+2}}dt'
    +\dfrac{C}{\ba2}\mu_m^{-\ba2}(1+\ln(1+t))^p\sqrt{\ol{E}_{0,\leq|\a|+2}}.
    \end{split}
    \ee
    Collecting the estimates for $I_1$, $I_2$, $I_3$ and applying the Gronwall inequality yield
    \bee\label{boundsFa'l}
    \bes
    \da^{l+1}||F_{\a',l}||_{\normu2}&\les\da^{\frac{3}{2}}M\left(\da^{l+1}
    \|R_i^{\a+1}T^l\s{\de}\mu\|_{L^2(\Si_0^{\tilde{\da}})}+\da^{l+1}\|\lie_{R_i}^{\a'}
    \lie_T^l\chi'\|_{L^2(\Si_0^{\tilde{\da}})}\right)\\
    &+\dfrac{C}{\ba2}\mu_m^{-\ba2}\sqrt{\ol{E}_{0,\leq|\a|+2}}(1+\ln(1+t))^p\\
    &+\dfrac{1}{1+t}\int_0^t\da^{\frac{1}{2}}M(1+\ln(1+t'))(\da\sqrt{\wi{E}_{1,\leq
    |\a|+2}}+\sqrt{\wi{E}_{0,\leq|\a|+2}})dt'.
    \end{split}
    \ee
    \hs As a corollary, the same estimate holds for $\mu R_i^{\a'}T^l\s{\de}\mu$ due to
    \bee\label{boundsRia'Tlsdemu}
    \da^{l+1}||\mu R_i^{\a'}T^l\s{\de}\mu||_{\normu2}\leq
    \da^{l+1}||F_{\a',l}||_{\normu2}+\da^{l+1}||R_i^{\a'}T^l\breve{f}||_{\normu2}.
    \ee
    \section{\textbf{Top order energy estimates}}\label{section8}
\subsection{\textbf{Contributions of the top order acoustical terms to the error integrals}}
\hs To complete top order energy estimates, we need to estimate the following error integrals:
\bee\label{topordererrorintegral}
\bes
&-\int_{W_t^u}\til{\p}_{|\a|+2}\cdot K_0R_i^{\a+1}\fai/K_0R_i^{\a'}T^{l+1}\fai
dt'du'd\mu_g,\\
&-\int_{W_t^u}\til{\p}_{|\a|+2}\cdot (K_1+2(1+t))R_i^{\a+1}\fai/(K_1+2(1+t))R_i^{\a'}T^{l+1}\fai
dt'du'd\mu_g.
\end{split}
\ee
It has been shown that the contributions of top order acoustical terms to the error integrals are
\[
T\fai\cdot \sd R_i^{\a}tr\chi,\hs T\fai\cdot R_i^{\a'}T^l\s{\de}\mu.
\]
\hs Since $K_0=(1+\mu)L+\dl$, $K_1=\frac{2(1+t)}{\frac{1}{2}\til{tr}\chi}L$, which implies $(1+\mu)L$ has a decay factor $\frac{1+\ln(1+t)}{(1+t)^2}$ compared with $K_1$, then for the error integrals associated to $K_0$, it suffices to consider the contribution from $\dl$ in $K_0$.
\subsubsection{\textbf{Contribution associated with $K_0$}}
      We first consider the following space-time integral:
      \bee\label{sdRiatrchi}
      \int_{W_t^u}T\fai\cdot\sd R_i^{\a}tr\chi\cdot\dl R_i^{\a+1}\fai\ dt'du'd\mu_{\sg},
      \ee
      which can be bounded as
      \bee\label{K0main}
      \int_{0}^t\sup(\mu^{-1}|T\fai|)\|\mu\sd R_i^{\a}tr\chi\|_{L^2(\Si_{t'}^{u})}\|\dl R_i^{\a+1}\fai\|_{L^2(\Si_{t'}^{u})}dt'.
      \ee

      \hs For $\dl R_i^{\a+1}\fai$, it follows from the definition of $E_0$ that
      \bee
      \|\dl R_i^{\a+1}\fai\|_{L^2(\Si_{t}^{\da})}\leq\mu_m^{-\ba2}N(t)^p
      \sqrt{\ol{E}_{0,\leq|\a|+2}},
      \ee
      where $N(t)=1+\ln(1+t)$. For the term $\mu\sd R_i^{\a}tr\chi$, it follows from \eqref{boundssdRiatrchi} that
      \bee\label{K0trx}
      \bes
      \|\mu\sd R_i^{\a}tr\chi\|_{\normu2}&\leq\dfrac{1}{1+t}\int_0^t\da^{\frac{3}{2}}M\cdot N(t')\mu_m^{-\ba2}
      N(t')^p\sqrt{\ol{E}_{0,\leq|\a|+2}(t',u)}\ dt'\\
      &+\dfrac{C}{\ba2}\mu_m^{-\ba2}N(t')^p\sqrt{\ol{E}_{0,\leq|\a|+2}(t,u)}\\
      &+\dfrac{1}{1+t}\int_0^t\dfrac{\da^{\frac{1}{2}}M\cdot N(t')}{1+t'}
      \mu_m^{-\ba2-\frac{1}{2}}N(t')^q\sqrt{\ol{E}_{1,\leq|\a|+2}}\ dt'\\
      &+\dfrac{\da^{\frac{1}{2}}M}{(1+t)^2}\mu_m^{-\ba2}N(t)^q\left(
      \int_0^u\ol{F}_{1,\leq|\a|+2}\ du'\right)^{\frac{1}{2}}.
      \end{split}
      \ee
      \hs To deal with the integral involving $\mu^{-1}$ above, we have to consider it in the shock region and in the non-shock region. It follows from Lemma\ref{crucial} that for $t\geq t_0$
      \bee
      \int_{-2}^{t_0}\mu_m^{-\ba2-1}\leq\mu_m^{-\ba2},\quad \int_{-2}^t\mu_m^{-\ba2-1}\leq\dfrac{C}{\ba2}\mu_m^{-\ba2}.
      \ee
      Hence, it follows from \eqref{K0main}, \eqref{K0trx} and Lemma\ref{crucial} that
      \bee\label{K0trchiL2final1}
      \bes
      \eqref{K0main}&\leq \dfrac{CM\da^{\frac{1}{2}}}{2\ba2-\frac{1}{2}}\mu_m^{-2\ba2}\left(N(t)^{2p}
      \ol{E}_{0,\leq|\a|+2}+N(t)^{2q}\ol{E}_{1,\leq|\a|+2}\right)\\
      &+\int_0^t\dfrac{\da^{\frac{3}{2}}M(1+\ln(1+t'))^2}{(1+t')^2}
      \mu_m^{-2\ba2}(N(t')^{2p}\ol{E}_{0,\leq|\a|+2}+N(t')^{2q}\ol{E}_{1,\leq|\a|+2})\ dt'\\
      &+\dfrac{CM\da^{\frac{1}{2}}}{(1+t)^2}N(t)^{2q}\int_0^u
      \ol{F}_{0,\leq|\a|+2}\ du'.
      \end{split}
      \ee
      Later we will choose $\ba2$ suitably large.\\
     \hs Next we consider the space-time integral:
     \bee\label{Ria'TlsdemuK0}
     \da^{2l+2}\int_{W_t^u}T\fai\cdot R_i^{\a'}T^l\s{\de}\mu\cdot\dl R_i^{\a'} T^{l+1}\fai\ dt'du'd\mu_{\sg},
     \ee
     which can be bounded as
      \bee\label{K0main2}
      \da^{2l+2}\int_{0}^t\sup(\mu^{-1}|T\fai|)\cdot\|\mu R_i^{\a'}T^l\s{\de}\mu\|_{L^2(\Si_{t'}^u)}\cdot\|\dl R_i^{\a'}T^{l+1}\fai\|_{L^2(\Si_{t'}^u)}dt'.
      \ee
      Similarly,
      \bee\da^{l+1}\|\dl R_i^{\a'}T^{l+1}\fai\|_{\normu2}\leq
      \mu_m^{-\ba2}N(t)^p\sqrt{\ol{E}_{0,\leq|\a|+2}}.
      \ee
       For $\mu R_i^{\a'}T^l\s{\de}\mu$, it follows from \eqref{boundsFa'l} and \eqref{boundsRia'Tlsdemu} that
      \begin{equation*}
      \bes
      \da^{l+1}\|\mu R_i^{\a'}T^l\s{\de}\mu\|_{\normu2}&\les\dfrac{1}{1+t}\int_0^t\da^{\frac{1}{2}}M\mu_m^{-\ba2} N(t')^{p+2}\sqrt
      {\ol{E}_{0,\leq|\a|+2}}\ dt'+\mu_m^{-\ba2}N(t)^p\sqrt{\ol{E}_{0,\leq|\a|+2}}\\
      &+\dfrac{1}{1+t}\int_0^t\da^{\frac{3}{2}}M\mu_m^{-\ba2} N(t')^{q+2}
      \sqrt{\ol{E}_{1,\leq|\a|+2}}\ dt'\\
      &+\dfrac{\da^{\frac{3}{2}}M}{1+t}\mu_m^{-\ba2}N(t)^q\left(\int_0^u\ol{F}_{1,\leq
      |\a|+2}du'\right)^{\frac{1}{2}}.
      \end{split}
      \end{equation*}
      \hs This together with \eqref{K0main2} and \eqref{boundsRia'Tlsdemu} shows
      \bee\label{K0muL2final1}
      \bes
      \eqref{Ria'TlsdemuK0}&\leq
      \dfrac{C\da^{\frac{1}{2}}M}{2\ba2}\mu_m^{-2\ba2}N(t)^{2p}\ol{E}_{0,\leq|\a|+2}+\dfrac{C\da^2M^2}{2\ba2}\mu_m^{-2\ba2}N(t)^{2q}\left(\ol{E}_
      {1,\leq|\a|+2}+\int_0^u\ol{F}_{1,\leq|\a|+2}\right)\\
      &+\int_0^t\dfrac{\da M^2(1+\ln(1+t'))^3}{(1+t')^2}\mu_m^{-2\ba2}
      (N(t')^{2p}\ol{E}_{0,\leq|\a|+2}+N(t')^{2q}\ol{E}_{1,\leq|\a|+2})\ dt'.
      \end{split}
      \ee

\subsubsection{\textbf{Contribution associated with $K_1$}}
We need the following elementary lemma of calculus on the manifold.
\begin{lem}\label{change}
Let $f,g$ be arbitrary functions defined on $S_{t,u}$ and $X$ be an $S_{t,u}$ tangent vector field. Then, it holds that
\begin{align}
\int_{S_{t,u}}f(Xg)d\mu_{\sg}&=-\int_{S_{t,u}}\left\{g(Xf)+\dfrac{1}{2}(tr\pre X {} {\s{\pi}})fg\right\}d\mu_{\sg},\\
\dfrac{\pa}{\pa t}\int_{S_{t,u}}fd\mu_{\sg}&=\int_{S_{t,u}}(L+tr\chi)fd\mu_{\sg}.
\end{align}
\end{lem}

We start with the space-time integral:
      \bee\label{K1main}
      \int_{W_t^u} T\fai\cdot\sd R_i^{\a}tr\chi\cdot (K_1+2(1+t'))R_i^{\a+1}\fai\ dt'du'd\mu_{\sg}.
      \ee
      \hs Due to Lemma\ref{change}, the space-time integral \eqref{K1main} can be rewritten as
      \bee
      \bes
      &\int_{W_t^u}\dfrac{2(1+t')}{\frac{1}{2}\til{tr}\chi}\cdot
      T\fai\cdot \sd R_i^{\a}tr\chi\cdot(L+\dfrac{1}{2}\til{tr}\chi)(R_i^{\a+1}\fai)\
      dt'du'd\mu_{\sg}\\
      &=\int_{W_t^u}(L+\til{tr}\chi)\left[\wv\cdot T\fai\cdot\sd R_i^{\a}tr\chi\cdot R_i^{\a+1}\fai\right]\ dt'du'd\mu_{\sg}\\
      &-\int_{W_t^u}(L+\dfrac{1}{2}\til{tr}\chi)\left[\wv\cdot T\fai\cdot\sd R_i^{\a}tr\chi\right]\cdot R_i^{\a+1}\fai\ dt'du'd\mu_{\sg}\\
      &=\int_{\Si_t^u}\wv\cdot T\fai\cdot\sd R_i^{\a}tr\chi\cdot R_i^{\a+1}\fai-
      \int_{\Si_0^u}\wv\cdot T\fai\cdot\sd R_i^{\a}tr\chi\cdot R_i^{\a+1}\fai\\
      &-\int_{W_t^u}(L+\dfrac{1}{2}\til{tr}\chi)\left[\wv\cdot T\fai\cdot\sd R_i^{\a}tr\chi\right]\cdot R_i^{\a+1}\fai\ dt'du'd\mu_{\sg}\\
      &:=I_0-I_1-I_2,
      \end{split}
      \ee
      where $I_1$ can be bounded by the "initial energy". It follows from Lemma\ref{change} and Proposition\ref{relationangular} that
      \bee
      \bes
      I_0&\les \int_{\Si_t^u}\wv\cdot\dfrac{1}{1+t'} T\fai\cdot R_i^{\a+1}tr\chi\cdot R_i^{\a+1}\fai\\
      &=-\int_{\Si_t^u}\dfrac{4}{\til{tr}\chi}R_i^{\a}tr\chi\left[
      R_i(T\fai\cdot R_i^{\a+1}\fai)+T\fai\cdot R_i^{\a+1}\fai\cdot\dfrac{1}{2}
      tr\srpi\right]\\
      &=-\int_{\Si_t^u}\dfrac{4}{\til{tr}\chi}T\fai\cdot R_i^{\a}tr\chi\cdot R_i^{\a+2}\fai-
      \int_{\Si_t^u}\dfrac{4}{\til{tr}\chi}R_i^{\a+1}\fai\cdot R_i^{\a}tr\chi\left(
      R_i T\fai+T\fai\cdot\dfrac{1}{2}tr\srpi\right)\\
      &:=-H_0-H_1.
      \end{split}
      \ee
      \hs Since $|tr\srpi|\les\dfrac{\da^{\frac{1}{2}}M(1+\ln(1+t))}{1+t}$, then it suffices to estimates $H_0$. Since the term $R_i^{\a}tr\chi$ in $H_0$ is not a top order term, then it follows from Proposition\ref{Riachi'L2} that
      \bee
      H_0\les \da^{\frac{1}{2}}M(1+t)\|R_i^{\a}tr\chi\|_{\normu2}\cdot\|\sd R_i^{\a+1}\fai\|_{\normu2}\les\dfrac{C\da^{\frac{1}{2}}M}{\ba2}\mu_m^{-2\ba2}N(t)^{2q}
      \ol{E}_{1,\leq|\a|+2}.
      \ee
      \hs For $I_2$, it holds that
      \bee
      \bes
      I_2&=\int_{W_t^u}(L+\til{tr}\chi)(\sd R_i^{\a}tr\chi)
      \wv\cdot T\fai\cdot R_i^{\a+1}\fai\ dt'du'd\mu_{\sg}\\
      &+\int_{W_t^u}(L-\dfrac{1}{2}\til{tr}\chi)\left(\wv\cdot T\fai\right)
      \sd R_i^{\a}tr\chi\cdot R_i^{\a+1}\fai\ dt'du'd\mu_{\sg}\\
      &:=A_1+A_2.
      \end{split}
      \ee
      \hs We first estimate $A_2$ and rewrite the term $(L-\dfrac{1}{2}\til{tr}\chi)(\wv\cdot T\fai)$ as
      \bee\label{Tfaiwv}
      \left(L+\dfrac{1}{2}\til{tr}\chi\right)T\fai\cdot\wv+
      (L-\til{tr}\chi)(\wv)\cdot T\fai.
      \ee
      \hs It follows from the proof of Proposition\ref{accruatemu1} that
      \bee
      \left|(L+\dfrac{1}{2}\til{tr}\chi)T\fai\right|\les\max\left\{
      \dfrac{\da^{\frac{1}{2}}M}{(1+t)^{\lam+1}},\dfrac{\da^{\frac{3}{2}}M}{
      (1+t)^3}\right\}\leq\dfrac{\da^{\frac{1}{2}}M}{(1+t)^{2+\ep}},
      \ee
      where $0<\ep<\lam-1$ is a constant. For the second term in \eqref{Tfaiwv}, one can rewrite this term as
      \bee
      \bes
      &(L-\til{tr}\chi)\left(\wv\right)=\dfrac{1+t}{(\frac{1}{2}\til{tr}\chi)^2}
      \left(\dfrac{\til{tr}\chi}{1+t}-L\til{tr}\chi-(\til{tr}\chi)^2\right)^2\\
      =&\dfrac{1+t}{(\frac{1}{2}\til{tr}\chi)^2}\cdot\left(\dfrac{-2u}{(1-u+t)^2(1+t)}+(\dfrac{1}{1+t}-\dfrac{4}{1-u+t})
      \til{tr}\chi'-(\til{tr}\chi')^2-L\til{tr}\chi'\right),
      \end{split}
      \ee
      which can be bounded by $\da^{\frac{1}{2}}M
      (1+\ln(1+t))$. Hence, it follows that
      \bee\label{1+tepdaM}
      \left|(L-\dfrac{1}{2}\til{tr}\chi)\left(\wv\cdot T\fai\right)\right|\leq
      \dfrac{\da^{\frac{1}{2}}M}{(1+t)^{\ep}}.
      \ee
      \hs Thus,
      \bee
      \bes
      |A_2|&\les\int_0^t\dfrac{\da^{\frac{1}{2}}M}{(1+t')^{\ep}}\mu^{-1}
      \|\mu\sd R_i^{\a}tr\chi\|_{L^2(\Si_{t'}^u)}\cdot\|R_i^{\a+1}\fai\|_
      {L^2(\Si_{t'}^u)}\ dt'\\
      &\leq\dfrac{C\da^{\frac{1}{2}}M}{2\ba2}\mu_m^{-2\ba2}N(t)^{2p}\ol{E}_{0,\leq|\a|+2}+\dfrac{C\da^\frac{1}{2}M}{2\ba2}\mu_m^{-2\ba2}N(t)^{2q}\left(\ol{E}_
      {1,\leq|\a|+2}+\int_0^u\ol{F}_{1,\leq|\a|+2}\right)\\
      &+\int_0^t\dfrac{\da M^2(1+\ln(1+t'))^2}{(1+t')^{1+\ep}}\mu_m^{-2\ba2}
      (N(t')^{2p}\ol{E}_{0,\leq|\a|+2}+N(t')^{2q}\ol{E}_{1,\leq|\a|+2})\ dt'.
      \end{split}
      \ee
      \hs For $A_1$, note that
      \bee
      (L+\til{tr}\chi)\sd R_i^{\a}tr\chi=\sd(L+\til{tr}\chi) R_i^{\a}tr\chi-(\sd \til{tr}\chi)(R_i^{\a}tr\chi).
      \ee
      \hs Since $|\sd \til{tr}\chi|\les\dfrac{\da^{\frac{3}{2}}M(1+\ln(1+t))}
      {(1+t)^3}$, it suffices to estimate the contribution from $\sd(L+\til{tr}\chi) R_i^{\a}tr\chi$. It follows from Lemma\ref{change} that
      \bee
      \bes
      &\int_{W_t^u}\sd(L+\til{tr}\chi) R_i^{\a}tr\chi\cdot T\fai\cdot R_i^{\a+1}\fai\ dt'du'd\mu_{\sg}\\
      \les&\int_{W_t^u}(1+t)^{-1}R_i(L+\til{tr}\chi)R_i^{\a}tr\chi\cdot T\fai\cdot R_i^{\a+1}\fai\ dt'du'd\mu_{\sg}\\
      =&-\int_{W_t^u}(1+t)^{-1}(L+\til{tr}\chi)R_i^{\a}tr\chi\cdot R_i^{\a+2}\fai\cdot T\fai\ dt'du'd\mu_{\sg}\\
      -&
      \int_{W_t^u}(1+t)^{-1}(L+\til{tr}\chi)R_i^{\a}tr\chi\cdot R_i^{\a+1}\fai\cdot\left(T\fai\cdot\dfrac{1}{2}tr\srpi+RT\fai\right)\ dt'du'd\mu_{\sg}\\
      :=&-A_{11}-A_{12}.
      \end{split}
      \ee
      It suffices to estimate $A_{11}$ due to $|tr\srpi|\les\dfrac{\da^{\frac{1}{2}}M(1+\ln(1+t))}{1+t}$.\\
      \hs Note that
      \bee
      (L+\til{tr}\chi)R_i^{\a}\chi=\left(L+\dfrac{2}{1-u+t}\right)
      R_i^{\a}tr\chi'+\til{tr}\chi'R_i^{\a}tr\chi,
      \ee
      and
      \bee
      \bes
      (L+\dfrac{2}{1-u+t})tr\chi'&=
      etr\chi'+|\chi'|^2-tr\a'+a\frac{\eta^{-1}\hat{T}\fe}{(1+t)^{\lam}}tr\chi'\\
      &+\dfrac{2}{1-u+t}\left(1+e+a\frac{\eta^{-1}\hat{T}\fe}{(1+t)^{\lam}}\right)tr\chi':=\p_0,\\
      (L+\dfrac{2}{1-u+t})R_i^{\a}tr\chi'&=R_i^{\a}\p_0
      +\sum_{|\be_1|+|\be_2|=|\a|-1}
      R_i^{\be_1}\srpi_LR_i^{\be_2}tr\chi',
      \end{split}
      \ee
      where
      \begin{equation*}
      \|R_i^{\be_1}\srpi_LR_i^{\be_2}tr\chi'\|_{\normu2}\les
      \dfrac{\da^{\frac{1}{2}}M(1+\ln(1+t))}{(1+t)^3}
      \int_0^t\dfrac{\da^{\frac{3}{2}}M(1+\ln(1+t'))}{1+t'}\mu_m^{-\frac{1}{2}}\sqrt{\wi{E}_{1,|\a|+2}(t',u)}dt',
      \end{equation*}
       due to Proposition\ref{Riachi'L2}, similar for the terms $R_i^{\a}(|\chi|^2), R_i^{\be_1}e\cdot R_i^{\be_2}tr\chi',R_i^{\a}tr\chi'$ in $R_i^{\a}tr\a'$ except the top order term $R_i^{\a}\s{\de}h.$ Indeed, one can bound this term as
     \bee
     \|R_i^{\a}\s{\de}h\|_{\normu2}\les\|\sd^2R_i^{\a}\fai\|_{\normu2}\les\dfrac{1}{(1+t)^2}\mu_m^{-\ba2}N(t)^q\sqrt{\ol{E}_{1,\leq|\a|+2}}.
     \ee
     Therefore,
     \bee
     |A_{11}|\les\dfrac{C\da^{\frac{1}{2}}M}{2\ba2-1}\mu_m^{-2\ba2}N(t)^{2q}
     \ol{E}_{1,\leq|\a|+2}+\int_0^t\dfrac{\da^{\frac{1}{2}}M}{(1+t)^2}\mu_m^{-2\ba2}N(t)^{2q}
     \ol{E}_{1,\leq|\a|+2}\ dt'.
     \ee
     \hs Collecting the results above yields
     \bee\label{Ria+1trchiK1final}
     \bes
     \eqref{K1main}&\leq\dfrac{C\da^{\frac{1}{2}}M}{2\ba2}\mu_m^{-2\ba2}N(t)^{2p}\ol{E}_{0,\leq|\a|+2}+\dfrac{C\da^\frac{1}{2}M}{2\ba2}\mu_m^{-2\ba2}N(t)^{2q}\left(\ol{E}_
      {1,\leq|\a|+2}+\int_0^u\ol{F}_{1,\leq|\a|+2}\right)\\
      &+\int_0^t\dfrac{\da M^2(1+\ln(1+t'))^2}{(1+t')^{1+\ep}}\mu_m^{-2\ba2}
      (N(t')^{2p}\ol{E}_{0,\leq|\a|+2}+N(t')^{2q}\ol{E}_{1,\leq|\a|+2})\ dt'.
     \end{split}
     \ee
 Finally, we consider the space-time integral:
      \bee\label{Ria'TlsdemuK1}
      \da^{2l+2}\int_{W_t^u}T\fai\cdot R_i^{\a'}T^l\s{\de}\mu\cdot (K_1+2(1+t'))R_i^{\a'}T^{l+1}\fai\ dt'du'd\mu_{\sg},
      \ee
      which can be written similarly as
      \bee\label{K1main2}
      \bes
      &\da^{2l+2}\int_{W_t^u}(L+\til{tr}\chi)\left[\wv\cdot T\fai\cdot R_i^{\a'}T^l\s{\de}\mu\cdot R_i^{\a'}T^{l+1}\fai\right]\ dt'du'd\mu_{\sg}\\
      -&\da^{2l+2}\int_{W_t^u}(L+\frac{1}{2}\til{tr}\chi)
      \left[\wv\cdot T\fai\cdot R_i^{\a'}T^l\s{\de}\mu\right]R_i^{\a'}T^{l+1}\fai\ dt'du'd\mu_{\sg}\\
      =&\da^{2l+2}\int_{\Si_t^u}\wv\cdot T\fai\cdot R_i^{\a'}T^l\s{\de}\mu\cdot R_i^{\a'}T^{l+1}\fai\\
      -&\da^{2l+1}\int_{\Si_{0}^u}\wv\cdot T\fai\cdot R_i^{\a'}T^l\s{\de}\mu
      \cdot R_i^{\a'}T^{l+1}\fai\\
      -&\da^{2l+2}\int_{W_t^u}(L+\frac{1}{2}\til{tr}\chi)
      \left[\wv\cdot T\fai\cdot R_i^{\a'}T^l\s{\de}\mu\right]R_i^{\a'}T^{l+1}\fai\ dt'du'd\mu_{\sg}\\
      :=&J_0-J_1-J_2,
      \end{split}
      \ee
      where $J_1$ can be bounded by the "initial energy". Similar as before, it follows from Lemma\ref{change} that
      \bee
      \bes
      |J_0|&\les\int_{\Si_t^u}(1+t')^2T\fai\cdot R_i^{\a'}T^l\s{\de}\mu
      \cdot R_i^{\a'}T^{l+1}\fai\ dt'du'd\mu_{\sg}\\
      &=-\da^{2l+2}\int_{\Si_t^u}(1+t')^2T\fai\cdot R_i^{\a'-1}T^l\s{\de}\mu\cdot R_i^{\a'+1}
      T^{l+1}\fai\\
      &-\da^{2l+2}\int_{\Si_t^u}(1+t')^2R_i^{\a'-1}T^l\s{\de}\mu\cdot R_i^{\a'}T^{l+1}\fai
      \left(R_iT\fai+T\fai\dfrac{1}{2}tr\srpi\right)\\
      &:=-H_1-H_2,
      \end{split}
      \ee
      where it suffices to estimate $H_1$. Since the term $R_i^{\a'-1}T^l\s{\de}\mu$ is not a top order term, then it follows from Proposition\ref{estzichi'zimuL2} that
      \bee
      |H_1|\les\dfrac{\da^{\frac{1}{2}}M}{1+t}\da^{l+1}\|R_i^{\a'+1}T^l\mu\|_{L^2(\Si_t^{u})}
      \cdot\da^{l+1}\|R_i^{\a'+1}T^{l+1}\fai\|_{\normu2}\leq\dfrac{C\da^{\frac{1}{2}}M}{\ba2-\frac{1}{2}}N(t)^{2q}
      \ol{E}_{1,\leq|\a|+2}.
      \ee
      \hs For the space-time integral $J_2$, it holds that
      \bee
      \bes
      J_2&=-\int_{W_t^u}\wv\cdot(L+\til{tr}\chi)R_i^{\a'}T^l\s{\de}\mu\cdot
      T\fai\cdot R_i^{\a'}T^{l+1}\fai\ dt'du'd\mu_{\sg}\\
      &-\int_{W_t^u}R_i^{\a'}T^l\s{\de}\mu\cdot R_i^{\a'}T^{l+1}\fai
      \left[(L-\frac{1}{2}\til{tr}\chi)\left(T\fai\cdot\wv\right)\right]\ dt'du'd\mu_{\sg}\\
      &:=-B_1-B_2.
      \end{split}
      \ee
      It follows from \eqref{1+tepdaM} that
      \bee
      \bes
      |B_2|&\leq\int_0^t\dfrac{\da^{\frac{1}{2}}M}{(1+t)^{\ep}}\mu^{-1}
      \|R_i^{\a'}T^{l+1}\fai\|_{L^2(\Si_{t'}^u)}\cdot\|\mu R_i^{\a'}T^l\s{\de}\mu\|_{L^2(\Si_{t'}^u)}\ dt'\\
      &\leq\dfrac{CM\da^{\frac{1}{2}}}{2\ba2}\mu_m^{-2\ba2}N(t)^{2p}
      \ol{E}_{0,\leq|\a|+2}+\dfrac{C\da^2M^2}{2\ba2}\mu_m^{-2\ba2}N(t)^{2q}\left[\ol{E}_{1,\leq
      |\a|+2}+\int_0^u\ol{F}_{1,\leq|\a|+2}du'\right]\\
      &+\int_0^t\dfrac{\da M^2(1+\ln(1+t'))^3}{(1+t')^{1+\ep}}
      \mu_m^{-2\ba2}(N(t')^{2p}\ol{E}_{0,\leq|\a|+2}+N(t')^{2q}\ol{E}_{1,\leq|\a|+2})\ dt'.
      \end{split}
      \ee
      \hs Note that
      \bee
      \bes
      (L+\til{tr}\chi)R_i^{\a'}T^l\s{\de}\mu&=R_i(L+\til{tr}\chi)R_i^{\a¡®-1}T^l\s{\de}\mu-(\srpi_L+R_i\til{tr}\chi)R_i^{\a'-1}T^l\s{\de}\mu.
      \end{split}
      \ee
      \hs Since $|\srpi|+|R_i\til{tr}\chi|\les\dfrac{\da^{\frac{1}{2}}M(1+\ln(1+t))}{(1+t)^2}$, it suffices to consider the contribution from $R_i(L+\til{tr}\chi)R_i^{\a¡®-1}T^l\s{\de}\mu$. It holds that
      \bee
      \bes
      &\da^{2l+2}\int_{W_t^u}(1+t')^2R_i(L+\til{tr}\chi)R_i^{\a-1}T^l\s{\de}
      \mu\cdot T\fai\cdot R_i^{\a'}T^{l+1}\fai\ dt'du'd\mu_{\sg}\\
      =&-\da^{2l+2}\int_{W_t^u}(1+t')^2(L+\til{tr}\chi)R_i^{\a'-1}T^l\s{\de}\mu
      \cdot R_i^{\a'+1}T^{l+1}\fai\cdot T\fai\ dt'du'd\mu_{\sg}\\
      &-\da^{2l+2}\int_{W_t^u}(1+t')^2(L+\til{tr}\chi)R_i^{\a'-1}T^l\s{\de}\mu
      \cdot R_i^{\a'}T^{l+1}\fai\cdot \left[R_iT\fai+T\fai\dfrac{1}{2}tr\srpi\right]\ dt'du'd\mu_{\sg}\\
      :=&-B_{11}-B_{12},
      \end{split}
      \ee
      where it suffices to estimate $B_{11}$. Decompose $B_{11}$ as:
      \bee
      \bes
      B_{11}&=\da^{2l+2}\int_{W_t^u}(1+t')^2LR_i^{\a'-1}T^l\s{\de}\mu\cdot R_i^{\a'+1}T^{l+1}\fai\cdot T\fai\ dt'du'd\mu_{\sg}\\
      &+\da^{2l+2}\int_{W_t^u}(1+t')^2(\til{tr}\chi)R_i^{\a'-1}T^l\s{\de}\mu\cdot R_i^{\a'+1}T^{l+1}\fai\cdot T\fai\ dt'du'd\mu_{\sg}\\
      &:=B_{111}+B_{112},
      \end{split}
      \ee
      where $B_{112}$ has been estimated before since $R_i^{\a'-1}T^l\s{\de}\mu$ is not a top order term. It remains to estimate $B_{111}$. Commuting $R_i^{\a'-1}T^l\s{\de}$ with \eqref{transportmu} yields
      \bee
      LR_i^{\a'-1}T^l\s{\de}\mu=\dfrac{1}{2}\dfrac{dH}{dh}R_i^{\a'-1}T^{l+1}\s{\de}h
      +\dfrac{1}{\eta}\dfrac{d\eta}{dh}LR_i^{\a'-1}T^l\s{\de}h+\dfrac{1}{\eta}\hat{T}^iLR_i^{\a'-1}T^l\s{\de}\fai_i+\text{l.o.ts},
      \ee
      where
      \bee
      \da^{l+1}\|R_i^{\a'-1}T^{l+1}\s{\de}h\|_{\normu2}+\da^{l+1}\|LR_i^{\a'-1}T^l\s{\de}h\|_{\normu2}\les
      \dfrac{\da}{(1+t)^2}\mu_m^{-\ba2-\frac{1}{2}}N(t)^q
      \sqrt{\ol{E}_{1,\leq|\a|+2}}.
      \ee
      \hs Hence, $B_{11}$ can be bounded as
      \begin{equation*}
      \bes
      |B_{11}|&\leq\int_0^t(1+t')\da^{\frac{1}{2}}\|(L+\til{tr}\chi)
      R_i^{\a'-1}T^l\s{\de}\mu\|_{L^2(\Si_{t'}^u)}\cdot\|R_i^{\a'+1}T^{l+1}\fai
      \|_{L^2(\Si_{t'}^u)}\ dt'\\
      &\les\int_0^t(1+t')\da^{\frac{1}{2}}M\mu_m^{-\ba2-\frac{1}{2}}N(t')^q
      \sqrt{\ol{E}_{1,\leq|\a|+2}}\\
      &\cdot\left(\dfrac{1}{(1+t')^2}\mu_m^{-\ba2-\frac{1}{2}}N(t')^q\sqrt{\ol{E}
      _{1,\leq|\a|+2}}+\dfrac{1}{(1+t')^2}\int_0^{t'}\dfrac{\da}{1+\tau}
      \mu_m^{-\ba2-\frac{1}{2}}N(\tau)^q\sqrt{\ol{E}_{1,\leq|\a|+1}}d\tau\right)dt'\\
      &\leq\dfrac{C\da^{\frac{1}{2}}M}{2\ba2}\mu_m^{-2\ba2}N(t)^{2q}
      \ol{E}_{1,\leq|\a|+2}.
      \end{split}
      \end{equation*}
      \hs Consequently,
      \bee\label{Ria'TlsdemuK1final}
      \bes
      \eqref{Ria'TlsdemuK1}&\leq\dfrac{CM\da^{\frac{1}{2}}}{2\ba2}\mu_m^{-2\ba2}N(t)^{2p}
      \ol{E}_{0,\leq|\a|+2}+\dfrac{C\da^2M^2}{2\ba2}\mu_m^{-2\ba2}N(t)^{2q}\left[\ol{E}_{1,\leq
      |\a|+2}+\int_0^u\ol{F}_{1,\leq|\a|+2}du'\right]\\
      &+\int_0^t\dfrac{\da M^2(1+\ln(1+t'))^3}{(1+t')^{1+\ep}}
      \mu_m^{-2\ba2}(N(t')^{2p}\ol{E}_{0,\leq|\a|+2}+N(t')^{2q}\ol{E}_{1,\leq|\a|+2})\ dt',
      \end{split}
      \ee
where the contributions from the initial data are omitted.
\subsection{\textbf{Top order energy estimates}}
\hs Following the same argument in \cite{CZD1} section7.2, one can complete the top order energy estimates associated to $K_0$ and $K_1$ as follows.
\subsubsection{\textbf{Estimates associated with $K_1$}}
\hs Let $Z_i\in\{Q,T,R_i\}$. It follows from \eqref{energy1},\eqref{energy2} and Remark\ref{energyremark} that
\bee\label{toporderK11}
\bes
&\sum_{|\a'|\leq|\a|}\da^{2l'}\left(\til{E}_1(Z^{\a'+1}_i\fai)(t,u)
+\til{F}_1(Z^{\a'+1}_i\fai)(t,u)+K(Z^{\a'+1}_i\fai)(t,u)\right)\\
&\leq C\sum_{|\a'|\leq|\a|}\da^{2l'}(\til{E}_1(Z^{\a'+1}_i\fai)(0,u)+\til{E}_0(Z^{\a'+1}_i\fai)(0,u))+
C\sum_{|\a'|\leq|\a|}\da^{2l'}\int_{W^t_{u}}Q_{1}^{|\a|+2},
\end{split}
\ee
where $l'$ is the number of $T$'s in $Z_i^{\a'+1}$ and $Q_1^{|\a|+2}:=(-\p_{|\a|+2}K_1(Z_i^{\a+1}\fai)+\frac{1}{2}T^{\a\be}(Z_i^{\a+1}\fai)\pi_{1,\a\be})$. The error integrals $\int_{W_t^{u}}Q_{1}^{|\a|+2}$ can be estimated as follows.\\
 \hs The contributions from the top order acoustical terms in $\wi{\p}_{|\a|+2}$ have been estimated in last subsection, see\eqref{Ria+1trchiK1final} and \eqref{Ria'TlsdemuK1final}; while the contributions from the lower order terms can be bounded as\footnote{For detailed analysis, see (7.45)-(7.49) in\cite{CZD1}.}:
    \bee
    \bes
    &\mu_m^{-2\ba2}\int_0^t\dfrac{\da^{\frac{1}{2}}M\cdot N(t')^q}{(1+t')^2}
    \sqrt{\ol{E}_{1,\leq|\a|+2}}\left(N(t')^p\sqrt{\ol{E}_{0,\leq|\a|+2}}+
    N(t')^q\left(\int_0^u\ol{F}_{1,\leq|\a|+2}\right)^{\frac{1}{2}}\right)\ dt'\\
    +&\mu_m^{-2\ba2}\da^{\frac{1}{2}}M\cdot N(t)^{2q}K_{\leq|\a|+2}(t,u)\\
    \leq&\mu_m^{-2\ba2}\int_0^t\dfrac{\da^{\frac{1}{2}}M}{(1+t')^2}N(t')^{2p}
    \ol{E}_{0,\leq|\a|+2}\ dt'+\mu_m^{-2\ba2}\da^{\frac{1}{2}}M\cdot N(t)^{2q}K_{\leq|\a|+2}(t,u)\\
    +&\mu_m^{-2\ba2}\int_0^t\dfrac{\da^{\frac{1}{2}}M}{(1+t')^2}N(t')^{2q}
    \left(\ol{E}_{1,\leq|\a|+2}+\int_0^u\ol{F}_{1,\leq|\a|+2}\right)\ dt',
    \end{split}
    \ee
    where $K_{\leq|\a|+2}$ is defined as
    \[
    K_{\leq|\a|+2}=\sup_{\tau\in [0,t]}\{(1+\ln(1+t'))^{-2q}\mu^{ u }_m(\tau)^{2b_{k+1}}\sum_{\fai}\sum_{|\a|=k-1}\da^{2l}K(Z^{\a+1}\fai)\}.
    \]
\hs Set $q=p+2$. Then, it follows that
    \bee\label{toporderK12}
    \bes
    &\mu_m^{2\ba2}\sum_{|\a'|\leq||\a|}\da^{2l'}\left(\til{E}_1(Z^{\a'+1}_i\fai)(t,u)
     +\til{F}_1(Z^{\a'+1}_i\fai)(t,u)+K(Z^{\a'+1}_i\fai)(t,u)\right)\\
     \leq&C\sum_{|\a'|\leq|\a|}\da^{2l'}(\til{E}_1(Z^{\a'+1}_i\fai)(0,u)+\til{E}_0(Z^{\a'+1}_i\fai)(0,u))+\dfrac{CM\da^{\frac{1}{2}}}{2\ba2}\mu_m^{-2\ba2}N(t)^{2p}
      \ol{E}_{0,\leq|\a|+2}\\
      &+\dfrac{C\da^2M^2}{2\ba2}\mu_m^{-2\ba2}N(t)^{2q}\left[\ol{E}_{1,\leq
      |\a|+2}+\int_0^u\ol{F}_{1,\leq|\a|+2}du'\right]\\
      &+\int_0^t\dfrac{\da M^2(1+\ln(1+t'))^3}{(1+t')^{1+\ep}}
      \mu_m^{-2\ba2}(N(t')^{2p}\ol{E}_{0,\leq|\a|+2}+N(t')^{2q}\ol{E}_{1,\leq|\a|+2})\ dt'.
    \end{split}
    \ee
   \hs Note that the right hand side above is increasing in $t$. Thus, \eqref{toporderK12} implies
    \bee\label{toporderK13}
    \bes
    &\ol{E}_{1,\leq|\a|+2}(t,u)+\ol{F}_{1,\leq|\a|+2}(t,u)+K_{\leq|\a|+2}(t,u)\leq\text{contributions from the initial data}\\
    &+\dfrac{C\da^2M^2}{2\ba2}\ol{E}_{1,\leq|\a|+2}+\dfrac{C\da^{\frac{1}{2}}M}{2\ba2}
    N(t')^{-4}\ol{E}_{0,\leq|\a|+2}+\dfrac{C\da^2M^2}{2\ba2}\int_0^u\ol{F}_{1,\leq|\a|+2}\ du'\\
    &+\int_0^t\dfrac{\da M^2(1+\ln(1+t'))^3}{(1+t')^{1+\ep}}
    \left(\ol{E}_{1,\leq|\a|+2}+N(t')^{-4}\ol{E}_{0,\leq|\a|+2}\right)dt'+CM\da^{\frac{1}{2}}
    K_{\leq|\a|+2},
    \end{split}
    \ee
    where the contributions from the initial data are $C\sum_{|\a'|\leq|\a|}\da^{2l'}(\til{E}_1(Z_i^{\a'+1}\fai)+\til{E}_0(Z_i^{\a'+1}\fai))(0,u)$+other initial data, which can be omitted\footnote{See remark7.1 in\cite{CZD1}.}. Choosing $\ba2$ such that $\dfrac{CM}{\ba2}\leq\dfrac{1}{4}$ and keeping only $\ol{E}_1$ on the left hand side of \eqref{toporderK13} yield that for sufficiently small $\da$
    \bee
    \bes
    \ol{E}_{1,\leq|\a|+2}(t,u)&\leq C\sum_{|\a'|\leq|\a|}\da^{2l'}(\til{E}_1(Z^{\a'+1}_i\fai)(0,u)+\til{E}_0(Z^{\a'+1}_i\fai)(0,u))\\
    &+\int_0^t\dfrac{\da M^2(1+\ln(1+t'))^3}{(1+t')^{1+\ep}}\ol{E}_{1,\leq|\a|+2}dt'+\dfrac{C\da^{\frac{1}{2}}M}{2\ba2}N(t')^{-4}\ol{E}_{0,\leq|\a|+2}\\
    &+\int_0^t\dfrac{\da M^2(1+\ln(1+t'))^3}{(1+t')^{1+\ep}}N(t')^{-4}\ol{E}_{0,\leq|\a|+2}dt'+\dfrac{C\da^2M^2}{2\ba2}\int_0^u\ol{F}_{1,\leq|\a|+2}\ du'.
    \end{split}
    \ee
    \hs Applying the Gronwall inequality yields
    \begin{equation*}\label{toporderK14}
    \bes
   \ol{E}_{1,\leq|\a|+2}(t,u)&\leq C\sum_{|\a'|\leq|\a|}\da^{2l'}(\til{E}_1(Z^{\a'+1}_i\fai)(0,u)+\til{E}_0(Z^{\a'+1}_i\fai)(0,u))+\dfrac{C\da^{\frac{1}{2}}M}{2\ba2}N(t')^{-4}\ol{E}_{0,\leq|\a|+2}\\
    &+\dfrac{C\da^2M^2}{2\ba2}\int_0^u\ol{F}_{1,\leq|\a|+2}\ du'+\int_0^t\dfrac{\da M^2(1+\ln(1+t'))^3}{(1+t')^{1+\ep}}N(t')^{-4}\ol{E}_{0,\leq|\a|+2}dt'.
    \end{split}
    \end{equation*}
    \hs Applying the same argument to $\ol{F}_1$ yields
    \bee\label{toporderK1final}
    \bes
    &\ol{E}_{1,\leq|\a|+2}(t,u)+\ol{F}_{1,\leq|\a|+2}(t,u)+K_{\leq|\a|+2}(t,u)\leq C\sum_{|\a'|\leq|\a|}\da^{2l'}(\til{E}_1(Z^{\a'+1}_i\fai)(0,u)+\til{E}_0(Z^{\a'+1}_i\fai)(0,u))\\
    &+C\sum_{
    |\a'|\leq|\a|+2}\da^{2l'}\|\lie_{Z_i}^{\a'}\chi'\|_{L^2(\Si_{0}^{\tilde{\da}})}
    +C\sum_{|\a'|+l'\leq|\a|+2}\da^{2l'}\|Z_i^{\a'}T^l\mu\|_{
    L^2(\Si_{0}^{\tilde{\da}})}\\
    &+\dfrac{C\da^{\frac{1}{2}}M}{2\ba2}N(t')^{-4}\ol{E}_{0,\leq|\a|+2}
    +\int_0^t\dfrac{\da M^2(1+\ln(1+t'))^3}{(1+t')^{1+\ep}}N(t')^{-4}\ol{E}_{0,\leq|\a|+2}dt'.
    \end{split}
    \ee
    This completes the top order energy estimates associated to $K_1$.
    \subsubsection{\textbf{Estimates associated with $K_0$}}
    \hs Similarly, one has
\bee\label{toporderK01}
\bes
&\sum_{|\a'|\leq||\a|}\da^{2l'}\left(\til{E}_0(Z^{\a'+1}_i\fai)(t, u)+\til{F}_0(Z^{\a'+1}_i\fai)(t, u)\right)\\
\leq&C\sum_{|\a'|\leq|\a|}\da^{2l'}(\til{E}_0(Z^{\a'+1}_i\fai)(0, u)+\til{E}_1(Z^{\a'+1}_i\fai)(0,u))+C\sum_{|\a'|\leq|\a|}\da^{2l'}\int_{W_t^{ u}}Q_{0}^{|\a|+2},
\end{split}
\ee
where $l'$ is the number of $T$'s in $Z_i^{\a'+1}$ and $Q_0^{|\a|+2}:=(-\p_{|\a|+2}K_0(Z_i^{\a+1}\fai)+\frac{1}{2}T^{\a\be}(Z_i^{\a+1}\fai)\pi_{0,\a\be})$. The error integrals $\int_{W_t^{u}}Q_{0}^{|\a|+2}$ can be bounded as due to \eqref{toporderK1final}
    \bee\label{contributionslotK0}
    \bes
    &\dfrac{CM\da^{\frac{1}{2}}}{2\ba2}\mu_m^{-2\ba2}\left[N(t)^{2p}
      \ol{E}_{0,\leq|\a|+2}+N(t)^{2q}\left(\ol{E}_
      {1,\leq|\a|+2}+\int_0^u\ol{F}_{1,\leq|\a|+2}\right)\right]\\
      +&\int_0^t\dfrac{\da M^2(1+\ln(1+t'))^3}{(1+t')^2}\mu_m^{-2\ba2}
      (N(t')^{2p}\ol{E}_{0,\leq|\a|+2}+N(t')^{2q}\ol{E}_{1,\leq|\a|+2})\ dt'\\
    +&\mu_m^{-2\ba2}\int_0^t\dfrac{\da^{\frac{1}{2}}M}{(1+t')^2}
    N(t')^{2p}\ol{E}_{0,\leq|\a|+2}+\mu_m^{-2\ba2}\da^{\frac{1}{2}}M\cdot N(t)^{2p}\ol{E}_{0,\leq|\a|+2}.
    \end{split}
    \ee
     \hs Substituting \eqref{contributionslotK0} and the estimates established in the previous subsections into \eqref{toporderK01} yields
    \bee
    \bes
    &N(t)^{-2p}\mu_m^{2\ba2}\sum_{|\a'|\leq||\a|}\da^{2l'}
    \left(\til{E}_0(Z^{\a'+1}_i\fai)(t, u)+\til{F}_0(Z^{\a'+1}_i\fai)(t, u)\right)\\
\leq&C\sum_{|\a'|\leq|\a|}\da^{2l'}\left(E_0(Z^{\a'+1}_i\fai)(0, u)+E_1(Z_i^{\a'+1}\fai)(0,u)\right)\\
    +&C\left(\da^{\frac{1}{2}}M+\dfrac{\da^{\frac{1}{2}}M}{2\ba2}\right)
    \ol{E}_{0,\leq|\a|+2}+\int_0^t\dfrac{\da M^2(1+\ln(1+t'))^3}{(1+t')^2}\ol{E}_{0,\leq|\a|+2}\ dt'.
    \end{split}
    \ee
    \hs Applying the same argument as before yields
    \bee\label{toporderK02}
    \ol{E}_{0,\leq|\a|+2}(t,u)+\ol{F}_{0,\leq|\a|+2}(t,u)\\
\leq C\sum_{|\a'|\leq|\a|}\da^{2l'}\left(\til{E}_0(Z^{\a'+1}_i\fai)(0, u)+\til{E}_1(Z_i^{\a'+1}\fai)(0,u)\right).
    \ee
    Define a quantity depending on initial data as\footnote{Indeed, the derivatives involving $\chi'$ and $\mu$ can be absorbed by the former two due to the definitions.}:
    \bee
    \bes
    \mathcal{D}_{|\a|+2}&=\sum_{|\a'|\leq|\a|}\da^{2l'}\left(\til{E}_0(Z^{\a'+1}_i\fai)
    (0, u)+\til{E}_1(Z_i^{\a'+1}\fai)(0,u)\right)\\
    &+\sum_{
    |\a'|\leq|\a|+1}\da^{2l'}\|\lie_{Z_i}^{\a'}\chi'\|_{L^2(\Si_{0}^{\tilde{\da}})}
    +\sum_{|\a'|+l'\leq|\a|+2}\da^{2l'}\|Z_i^{\a'}T^l\mu\|_{
    L^2(\Si_{0}^{\tilde{\da}})}.
    \end{split}
    \ee
    Then, it follows from \eqref{toporderK02} and \eqref{toporderK1final} that
    \bee\label{toporderen}
    \bes
    \ol{E}_{1,\leq|\a|+2}(t,u)+\ol{F}_{1,\leq|\a|+2}(t,u)+K_{\leq|\a|+2}(t,u)
    &\leq C\mathcal{D}_{|\a|+2},\\
    \ol{E}_{0,\leq|\a|+2}(t,u)+\ol{F}_{0,\leq|\a|+2}(t,u)&\leq
    C\mathcal{D}_{|\a|+2}.
    \end{split}
    \ee
    \section{\textbf{The decent scheme}}\label{section9}
    \hs The following lemma refines Corollary\ref{coro2} whose proof is given in \cite{CZD1}.
\begin{lem}\label{refine}
Let $k$ and $c$ be positive constants with $k>1$. Then, if $\da$ is sufficiently small depending on an upper bound for $c$ and a lower bound for $k$, it holds that, for any $\tau\in[0,t]$,
\[\mu_m^c(t)\leq k\mu_m^c(\tau).
\]
\end{lem}
\subsection{\textbf{The next-to-top order estimates}}
\hs Recall that the top order is $|\a|+2=N_{top}+1$ and set $\ba2=b_{|\a|+1}+1$. The following computations are similar to the top order energy estimates with $|\a|+2=N_{top}+1$ replaced by $|\a|+1=N_{top}$.\\
\hs Corresponding to \eqref{K1main}, we consider the space-time integral
      \bee\label{nttK11}
      \bes
      &\int_{W_t^u}R_i^{\a}tr\chi\cdot T\fai\cdot (K_1+2(1+t))R_i^{\a}\fai\ dt'du'd\mu_{\sg}\\
      =&\int_{W_t^u}\dfrac{2(1+t')}{\frac{1}{2}\til{tr}\chi}R_i^{\a}tr\chi\cdot
      T\fai\cdot(L+\frac{1}{2}\til{tr}\chi)R_i^{\a}\fai\ dt'du'd\mu_{\sg}\\
      \leq C \da M&\left(\int_0^t\|R_i^{\a}tr\chi\|^2_{L^2(\Si_{t'}^u)}\ dt'
      \right)^{\frac{1}{2}}\cdot\left(\int_0^u\|(1+t)(L+\frac{1}{2}\til{tr}\chi)
      R_i^{\a}\fai\|^2_{L^2(C_{u'}^t)}\ du'\right)^{\frac{1}{2}}.
      \end{split}
      \ee
      \hs It follows from the definition of $\ol{F}_1$, \eqref{toporderen} and Proposition\ref{Riachi'L2} that
      \begin{equation*}
      \bes
      &\int_{W_t^u}R_i^{\a}tr\chi\cdot T\fai\cdot (K_1+2(1+t))R_i^{\a}\fai\ dt'du'd\mu_{\sg}\les\da^2M^2\mu_m^{-2\nba}N(t)^{4}\int_0^u\ol{F}_{1,
      \leq|\a|+1}\\
      &+\da^2M^2\cdot N(t)^4\ol{E}_{1,\leq|\a|+2}\int_0^t\dfrac
      {(1+\ln(1+t'))^6}{(1+t')^2}\mu_m^{-2\ba2+1}\ dt'\\
      &\les\da^2M^2\mu_m^{-2\nba}(1+\ln(1+t))^{4}\left(\mathcal{D}_{|\a|+2}+
      \int_0^u\ol{F}_{1,
      \leq|\a|+1}\right).
      \end{split}
      \end{equation*}
      \hs Next, corresponding to \eqref{Ria'TlsdemuK1}, we consider the contribution from $R_i^{\a'-1}T^l\s{\de}\mu$:
      \bee\label{nttK12}
      \bes
      &\da^{2l+2}\int_{W_t^u}R_i^{\a'-1}T^l\s{\de}\mu\cdot T\fai\cdot
      (K_1+2(1+t))R_i^{\a'-1}
      T^{l+1}\fai\ dt'du'd\mu_{\sg}\\
      \les &\da^{\frac{1}{2}}M\da^{l+1}
      \left(\int_{0}^t\|R_i^{\a'-1}T^l\s{\de}\mu\|^2_{L^2(\Si_{t'}^{u}
      )}\right)^{\frac{1}{2}}\cdot\da^{l+1}\left(\int_0^u\|(1+t)(L+\frac{1}{2}\til{tr}\chi
      )R_i^{\a'-1}T^{l+1}\fai\|
      ^2_{L^2(C_{u'}^t)}\right)^{\frac{1}{2}}.
      \end{split}
      \ee
      \hs Similarly, one can bound it as due to Proposition\ref{estzichi'zimuL2}
      \bee
      \bes
      &\da^{2l+2}\int_{W_t^u}R_i^{\a'-1}T^l\s{\de}\mu\cdot T\fai\cdot
      (K_1+2(1+t))R_i^{\a'-1}
      T^{l+1}\fai\ dt'du'd\mu_{\sg}\\
      &\leq C\da^{\frac{3}{2}}M\mu_m^{-2\nba}N(t)^4\int_0^u
      \ol{F}_{1,\leq|\a|+1}\ du'+C\da^{\frac{3}{2}}M\mu_m^{-2\nba}\ol{E}_{0,\leq
      |\a|+1}\\
      &+C\da^{\frac{3}{2}}M\cdot N(t)^4\ol{E}_{1,\leq|\a|+2}\int_0^t
      \dfrac{(1+\ln(1+t'))^4}{(1+t')^2}\mu_m^{-2\ba2+1}\ dt'\\
      &\leq C\da^{\frac{3}{2}}M\mu_m^{-2\nba}N(t)^4\left(\mathcal{D}_{|\a|+2}+
      \int_0^u
      \ol{F}_{1,\leq|\a|+1}\ du'\right)+C\da^{\frac{3}{2}}M\mu_m^{-2\nba}\ol{E}_{0,\leq
      |\a|+1}.
      \end{split}
      \ee
     \hs Corresponding to \eqref{sdRiatrchi}, we consider the space-time integral:
      \bee\label{nttK01}
      \int_{W_t^u}R_i^{\a}tr\chi\cdot T\fai\cdot\dl R_i^{\a}\fai\ dt'du'd\mu_{\sg},
      \ee
      which can be bounded as due to Proposition\ref{Riachi'L2}
      \bee
      \bes
       &\int_0^t|T\fai|\cdot\|R_i^{\a}tr\chi\|_{L^2(\Si_{t'}^u)}
      \cdot\|\dl R_i^{\a}\fai\|_{L^2(\Si_{t'}^u)}\ dt'\\
      \les&\da^2 M^2\sqrt{\mathcal{D}_{|\a|+2}}\sqrt{\ol{E}_{0,\leq
      |\a|+1}}\int_0^t\dfrac{(1+\ln(1+t'))^5}{(1+t')^2}\mu_m^{-\ba2+\frac{3}{2}}\ dt'\\
      \leq& C\da^2M^2\mu_m^{-2\nba}\mathcal{D}_{|\a|+2}+C\da^2M^2\mu_m^{-2\nba}\ol{E}_{0,
      \leq|\a|+1}.
      \end{split}
      \ee
      \hs Finally, corresponding to \eqref{Ria'TlsdemuK0}, we consider the space-time integral:
      \bee\label{nttK02}
      \bes
      &\da^{2l+2}\int_{W_t^u}R_i^{\a'-1}T^l\s{\de}\mu\cdot T\fai\cdot\dl R_i^{\a'-1}T^{l+1}\fai\ dt'du'd\mu_{\sg}\\
      \leq&\int_0^t|T\fai|\cdot\da^{l+1}\|R_i^{\a'-1}T^l\s{\de}\mu\|_{
      L^2(\Si_{t'}^u)}\cdot\da^{l+1}\|\dl R_i^{\a'-1}T^{l+1}\fai\|_{
      L^2(\Si_{t'}^u)}\ dt'.
      \end{split}
      \ee
      It follows from the definition of $\ol{E}_0$ and Proposition\ref{estzichi'zimuL2} that
      \bee
      \bes
      &\da^{2l+2}\int_{W_t^u}R_i^{\a'-1}T^l\s{\de}\mu\cdot T\fai\cdot\dl R_i^{\a'-1}T^{l+1}\fai\ dt'du'd\mu_{\sg}\\
      &\leq C\da^{\frac{3}{2}}M\sqrt{\mathcal{D}_{|\a|+2}}\sqrt
      {\ol{E}_{0,\leq|\a|+1}}\int_0^t\dfrac{(1+\ln(1+t'))^4}{(1+t')^2}
      \mu_m^{-2\ba2+\frac{3}{2}}\ dt'\\
      &\leq C\da^{\frac{3}{2}}M\mu_m^{-2\nba}
      \left(\mathcal{D}_{|\a|+2}+\ol{E}_{0,\leq|\a|+1}\right).
      \end{split}
      \ee
      We then turn to the next-to-top order energy estimates.\\
      \hs Similar to \eqref{toporderK11}, it holds that
      \bee
      \bes
      &\mu_m^{2\nba}N(t)^{-4}\sum_{|\a'|\leq|\a|-1}\da^{2l'}\left(\til{E}_1(Z_i^{\a'+1}\fai)+\til{F}_1(Z_i^{\a'+1}\fai)
      +K(Z_i^{\a'+1})\right)(t,u)\\
      \leq&C\sum_{|\a'|\leq|\a|-1}\da^{2l'}(E_{1}(Z_i^{\a'+1}\fai)(0,u)+E_{0}(Z_i^{\a'+1}\fai)(0,u))
      +C\mathcal{D}_{|\a|+2}+C\da M^2\int_0^u\ol{F}_{1,\leq|\a|+1}\ du'\\
      +&C\da M\int_{0}^t\dfrac{(1+\ln(1+t'))^2}{(1+t')^2}
      \left(\ol{E}_{1,\leq|\a|+2}+N(t')^{-4}\ol{E}_{0,\leq|\a|+1}\right)\ dt'\\
      +&C\da M\wi{E}_{1,\leq|\a|+1}+C\da M\ol{E}_{0,\leq|\a|+1}+C\da^{\frac{1}{2}}K_{\leq|\a|+1}.
      \end{split}
      \ee
      Then,
      \bee\label{ntteK11}
      \bes
      &\left(\ol{E}_{1,\leq|\a|+1}+\ol{F}_{1,\leq|\a|+1}+K_{\leq|\a|+1}\right)(t,u)\leq C\mathcal{D}_{|\a|+2}+C\da M\int_0^u\ol{F}_{1,\leq|\a|+1}\ du'+C\da M\wi{E}_{0,\leq|\a|+1}\\
      +&C\int_{0}^t\dfrac{(1+\ln(1+t'))^2}{(1+t')^2}
      \left(\ol{E}_{1,\leq|\a|+2}+N(t')^{-4}\ol{E}_{0,\leq|\a|+1}\right)\ dt'.
      \end{split}
      \ee
      \hs Applying the same argument as before and using Lemma\ref{refine} yield
      \bee\label{ntte1}
      \bes
      &\ol{E}_{1,\leq|\a|+1}+\ol{F}_{1,\leq|\a|+1}+K_{\leq|\a|+1}\leq C\mathcal{D}_{|\a|+2}
      +C\da M\ol{E}_{0,\leq|\a|+1}\\
      +&C\da M\int_{0}^t\dfrac{(1+\ln(1+t'))^2}{(1+t')^2}N(t')^{-4}\ol{E}_{0,\leq|\a|+1}\ dt'.
      \end{split}
      \ee
    \hs Similar to \eqref{toporderK01}, the following estimate holds for $E_0$ and $F_0$
      \bee
      \bes
      &\mu_m^{2\nba}\sum_{|\a'|\leq|\a|-1}\da^{2l'}\left(
      E_0(Z_i^{\a'+1}\fai)+F_0(Z_i^{\a'+1}\fai)\right)(t,u)\\
      \leq&C\sum_{|\a'|\leq|\a|-1}\da^{2l'}({E}_0(Z_i^{\a'+1}\fai)(0,u)+{E}_1(Z_i^{\a'+1}\fai)(0,u))
      +C\mathcal{D}_{|\a|+2}+C\da M\ol{E}_{0,\leq|\a|+1}\\
      +&C\da^{\frac{1}{2}}K_{\leq|\a|+1}(t,u)+C\da(1+\ln(1+t))^4\ol{E}_{1,\leq|\a|+1}+
      C(1+\ln(1+t))^4\int_0^u\ol{F}_{1,\leq|\a|+1}\ du'\\
      \leq&C\mathcal{D}_{|\a|+2}+C\da M\ol{E}_{0,\leq|\a|+1}+C\int_{0}^t
      \dfrac{(1+\ln(1+t'))^2}{(1+t')^2}\ol{E}_{0,\leq|\a|+1}\ dt'.
      \end{split}
      \ee
      \hs Applying the same argument as before yields
      \bee\label{ntte2}
      \left(\ol{E}_{0,\leq|\a|+1}+\ol{F}_{0,\leq|\a|+1}\right)(t,u)\leq
      C\mathcal{D}_{|\a|+2}+C\da M\int_{0}^t\dfrac{(1+\ln(1+t'))^2}{(1+t')^2}\ol{E}_{0,\leq|\a|+1}\ dt'.
      \ee
      \hs Applying the Gronwall inequality and noticing \eqref{ntte1} and \eqref{ntte2} yield
      \bee\label{nttefinal}
      \bes
      \ol{E}_{0,\leq|\a|+1}(t,u)+\ol{F}_{0,\leq|\a|+1}(t,u)&\leq C\mathcal{D}_{|\a|+2},\\
      \ol{E}_{1,\leq|\a|+1}(t,u)+\ol{F}_{1,\leq|\a|+1}(t,u)+K_{\leq|\a|+1}(t,u)&\leq
      C\mathcal{D}_{|\a|+2}.
      \end{split}
      \ee
\subsection{\textbf{The decent scheme}}
\hs Set
    \bee
    b_{|\a|+2-n}=\ba2-n,\quad b_{|\a|+1-n}=\ba2-1-n, \quad \ba2=[\ba2]+\dfrac{3}{4}.
    \ee
    Let the argument for the next-to-top order estimates be the $1-$st step and $n-$th step be the corresponding energy estimates with $\ba2$ and $\nba$ replaced by $b_{|\a|+2-n}$ and $b_{|\a|+1-n}$, respectively. Then, \textbf{as long as $b_{|\a|+1-n}>0$}, i.e. $n\leq[\nba]$, the $n-$th step can be proceeded exactly in the same way as the $1-$st step.\\
    \hs In the $n-$th step, one has to the following integral:
    \bee\label{nthmu}
    \int_{0}^t\dfrac{(1+\ln(1+t'))^6}{(1+t')^2}\mu_m^{-2b_{|\a|+2-n}+1}\ dt'.
    \ee
    \hs Since $n\leq[\nba]$, then $-2b_{|\a|+2-n}+1\leq-\dfrac{5}{2}$. Thus, the integral \eqref{nthmu} can be bounded as $C\mu_m^{-2b_{|\a|+2-n}+2}=C\mu_m^{-2
    b_{|\a|+1-n}}$ due to Lemma\ref{refine} and then the same argument in the $1-$st step can be applied to the $n-$th step. Hence, it holds that
    \bee
    \bes
    \ol{E}_{0,\leq|\a|+1-n}+\ol{F}_{0,\leq|\a|+1-n}&\leq C\mathcal{D}_{|\a|+2},\\
    \ol{E}_{1,\leq|\a|+1-n}+\ol{F}_{1,\leq|\a|+1-n}+K_{\leq|\a|+1-n}&\leq C\mathcal
    {D}_{|\a|+2},
    \end{split}
    \ee
    for $n=0,1,\cdots,[\nba]$.\\
    \hs Up to now, $\ol{E}_{0,\leq|\a|+1-n}=\sup\{\mu_m^{-2b_{|\a|+1-n}}E_{0,|\a|+1-n}\}$, the power of $\mu$ is still not eliminated. Hence, we consider the next step. Set $n=[\ba2]$ to be the final step. Then, $b_{|\a|+2-n}=\dfrac{3}{4}$, $b_{|\a|+1-n}=-\dfrac{1}{4}$, and the following integral needs to be considered:
    \bee
    \int_{0}^t\dfrac{(1+\ln(1+t'))^6}{(1+t')^2}\mu_m^{-\frac{1}{2}}\ dt'.
    \ee
    It follows from the proof of Lemma\ref{crucial} that
\begin{equation*}
\bes
&\int_{0}^t\dfrac{(1+\ln(1+t'))^6}{(1+t')^2}\mu_m^{-\frac{1}{2}}\ dt'\leq
\dfrac{(1+\ln(1+t))^2}{1+t}\int_0^t\dfrac{\mu_m^{-\frac{1}{2}}}{1+t'}\ dt'\\
&\leq\dfrac{(1+\ln(1+t))^2}{1+t}\dfrac{1}{-(1+s)A(s)\eta_m(s)+\frac{1}{4}}
\cdot\left[\mu_m(s)+\int_s^t\frac{1}{(1+\tau)A(\tau)}d\tau(-(1+s)A(s)\eta_m(s)+\frac{1}{4})\right]^{\frac{1}{2}}\\
&\leq C,
\end{split}
\end{equation*}
for any \textbf{fixed} $s\in[0,t]$.\\
\hs Set $b_{|\a|+1-n}=0$. Then, the final step can be proceeded in the same way as the previous steps and the following estimates hold for $n=[\ba2]$:
    \bee\label{desiredenergy}
    \bes
    \ol{E}_{0,\leq|\a|+1-n}+\ol{F}_{0,\leq|\a|+1-n}&\leq C\mathcal{D}_{|\a|+2},\\
    \ol{E}_{1,\leq|\a|+1-n}+\ol{F}_{1,\leq|\a|+1-n}+K_{\leq|\a|+1-n}&\leq C\mathcal
    {D}_{|\a|+2},
    \end{split}
    \ee
    which are the desired energy estimates since 
    the power of $\mu$ has been eliminated.
    \section{\textbf{Recovery of the bootstrap assumptions and completion of the proof}}\label{section10}
\begin{lem}\label{JA}(c.f. Lemma 17.1 in\cite{christodoulou2014compressible})
There exists a numerical constant $C$ 
such that
\bee
\|\fai\|_{L^{\infty}(S_{t,u})}\leq C\cdot(1+t)^{-1} S^{\frac{1}{2}}_{[2]}(\fai),
\ee
where
\bee
S_{[2]}(\fai)=\int_{S_{t,u}}\sum_{i,j}(|\fai|^2+|R_i\fai|^2+|R_iR_j\fai|^2)\ d\mu_{\sg}.
\ee
\end{lem}
\hs Denote $S_n(t,u)$ to be
\bee
S_n(t,u)=\sum_{|\a'|\leq n}\int_{S_{t,u}}|\da^{l'}Z_i^{\a'}\fai_{\gamma}|^2\ d\mu_{\sg},
\ee
which is the sum of the integrals on $S_{t,u}$ of square of all the variations of $\fai$ up to order $n$ where $l'$ is the number of $T'$s in $Z_i^{\a'}$. Then, it follows from Lemma\ref{faibound} that
\bee
S_{|\a|-[\ba2]+1}\leq C\tilde{\da}(\ol{E}_{0,\leq|\a|+1-[\ba2]}+\ol{E}_{1,\leq|\a|+1-[\ba2]}),
\ee
for all $(t, u)\in[0,t^{\ast})\times[0,\tilde{\da}]$. The constant $C$ here depends on $\dfrac{\sqrt{\det \sg(t,u)}}{\sqrt{\det\sg(t,0)}}$ and is a numerical constant. It follows from \eqref{desiredenergy} that
\[
\ol{E}_{0,\leq|\a|+1-[\ba2]}+\ol{E}_{1,\leq|\a|+1-[\ba2]}\leq C\cdot \mathcal{D}_{|\a|+2}.
\]
\hs Now for any variation $\fai$ up to order $|\a|+1-[\ba2]-2$, it holds that
\[S_2(\fai)\leq S_{|\a|-[\ba2]+1}(t, u).
\]
\hs Thus, by Lemma\ref{JA}, the following estimate holds for all $|\a'|\leq|\a|-1-[\ba2]$:
\bee
\da^{l'}\sup|Z_i^{\a'}\fai|\leq C_0(1+t)^{-1}\da^{\frac{3}{2}}\sqrt{\mathcal{D}_{|\a|+2}},
\ee
for all $(t,u)\in[0,t^{\ast})\times[0,\tilde{\da}]$. The constant $C_0$ here depends on some numerical constants.\\
\hs Hence, by choosing $M$ suitably large such that $C_0\sqrt{\mathcal{D}_{|\a|+2}}<M$, one can recover the bootstrap assumptions.\\
\hs To finish the proof, it suffices to prove $t^{\ast}=s^{\ast}$, i.e. either $t^{\ast}=+\infty$, then the smooth solution exists globally, or $t^{\ast}<\infty$, then at least at one point on $\Si_{t^{\ast}}^{\tilde{\da}}$ such that $\mu_m(t^{\ast})=0$, and a shock forms in finite time.\\
\hs If $t^{\ast}< s^{\ast}$, then \textbf{$\mu_m$ is positive on $\Si_{t^{\ast}}^{\tilde{\da}}$}. Hence, the \textbf{Jacobian $\de$ of the transformation between the acoustical coordinates and the rectangular coordinates never vanishes on $\Si_{t^{\ast}}^{\tilde{\da}}$}, i.e. the transformation between two coordinates is regular on $\Si_{t^{\ast}}^{\tilde{\da}}$. Moreover, in the acoustical coordinates, $\fai$ and its derivatives $\fai_{\gamma}$ are regular on $\Si_{t^{\ast}}^{\tilde{\da}}$ due to the bootstrap assumptions. 
Therefore, in the rectangular coordinates, $\fai$ and its derivatives $\fai_{\gamma}$ are regular on $\Si_{t^{\ast}}^{\tilde{\da}}$ which belongs to the Sobolev space $H^3$. By the standard local well-posedness theory, one can obtain an extension of the solution to some $t_1>t^{\ast}$, which is a contradiction!
\bibliographystyle{plain}
\bibliography{ref} 
\end{document}